\documentclass[a4paper,10pt]{article}
\usepackage[paper=a4paper,left=30mm,right=30mm,top=30mm,bottom=30mm]{geometry} 

\usepackage{amsmath, amsthm, amsfonts, amssymb,amscd}
\usepackage{bbm}
\usepackage{fancyhdr}

\usepackage{pstricks,pstricks-add,pst-math,pst-xkey}

\usepackage{graphicx}
\usepackage{psfrag}
\usepackage{wrapfig}
\usepackage{color}
\usepackage[abs]{overpic}
\definecolor{shadecolor}{gray}{0.90}
\usepackage{framed}
\usepackage[colorlinks]{hyperref}
\usepackage{url}
\usepackage{cancel}

\usepackage{soul}
\setstcolor{blue}

  \usepackage{lineno}
  \usepackage{setspace}

\theoremstyle{plain}

\theoremstyle{definition}

\theoremstyle{remark}
\newtheorem*{rem}{Remark}

\title{A structure-preserving split finite element discretization of the split 1D wave equations}
\author{Werner Bauer\footnote{Department of Mathematics, Imperial College London, London, United Kingdom; w.bauer@imperial.ac.uk}\ , \quad
       J\"orn Behrens\footnote{Department of Mathematics/CEN-Center for Earth System Research and 
       Sustainability, Universit\"at Hamburg, Hamburg, Germany}}
\date{}
\begin{document}

\maketitle


\begin{abstract}
 
  We introduce a new finite element (FE) discretization framework applicable 
  for covariant split equations. The introduction of additional differential forms (DF) 
  that form pairs with the original ones permits the splitting of the equations 
  into topological momentum and continuity equations and metric-dependent closure 
  equations that apply the Hodge-star operator. Our discretization framework conserves this 
  geometrical structure and provides for all DFs proper FE spaces such that the 
  differential operators (here gradient and divergence) hold in strong form. 
  We introduce lowest possible order discretizations of the split 1D wave equations, in 
  which the discrete momentum and continuity equations follow by trivial projections onto piecewise 
  constant FE spaces, omitting partial integrations. Approximating the Hodge-star by nontrivial Galerkin 
  projections (GP), the two discrete metric equations follow by projections onto either 
  the piecewise constant (GP0) or piecewise linear (GP1) space. 
  
  Out of the four possible realizations, our framework gives us three schemes with 
  significantly different behavior. 
  The split scheme using twice GP1 is unstable and shares the dispersion relation 
  with the P1--P1 FE scheme that approximates both variables by piecewise linear 
  spaces (P1). The split schemes that apply a mixture of GP1 and GP0 share the dispersion 
  relation with the stable P1--P0 FE scheme that applies piecewise linear and piecewise
  constant (P0) spaces. However, the split schemes exhibit second order convergence 
  for both quantities of interest. 
  For the split scheme applying twice GP0, we are not aware of a corresponding standard 
  formulation to compare with. Though it does not provide a satisfactory 
  approximation of the dispersion relation as short waves are propagated much too
  fast, the discovery of the new scheme illustrates the potential of
  our discretization framework as a toolbox to study and find FE schemes by 
  new combinations of FE spaces.

\end{abstract}



  \section{Introduction}
  \label{intro}

  The Finite Element (FE) method provides a powerful framework to discretize 
  partial differential equations (PDEs) and includes methods to prove the 
  discrete models' convergence, stability, and accuracy properties (see e.g. \cite{BoffiBrezzi2013,FenicsBook}).
  By offering flexibility in the choice of computational (unstructured, h/p-adapted) meshes (cf. \cite{Behrens2006}) 
  while providing an approximation of the continuous PDEs with the required order of accuracy, 
  FE discretizations are nowadays appreciated in all research areas that apply numerical modeling.
  
  Discretizations using finite element methods provide one important advantage over other methods: 
  Starting from a variational formulation the discretization follows simply by substituting the continuous 
  by discrete function spaces (Galerkin methods) while the differential operators remain unchanged.
  There exist a large variety of different suitable FE spaces to choose from. 
  However, not all choices lead to well-behaving schemes. In particular mixed FE schemes suffer from 
  this problem, where different variables of the PDE system are represented by different FE spaces.
  In such schemes, certain combinations of FE spaces lead to instabilities that exhibit spurious modes, 
  rendering the solutions useless, in particular when studying nonlinear phenomena. A famous example 
  for an unstable scheme is given by an approximation of both velocity and height fields of the 
  1D shallow-water equations by piecewise linear functions (cf. \cite{WaltersCarey1983} and Sect.~\ref{sec_standard}),
  where it is well known that equal order FE pairs are always unstable \cite{Cotterlbb2009}. 
 
  In order to avoid unsuitable choices, the Finite Element Exterior Calculus (FEEC) method \cite{ArnoldFEEC2006,ArnoldFEEC2010} 
  provides means for choosing a suitable pair of FE spaces that is guaranteed to lead to 
  a stable mixed discretization. In particular, FEEC puts geometrical constraints on the 
  FE spaces such that geometric properties, like the Helmholtz decomposition of vector-fields,
  are preserved in the discrete case. As a result, FEEC pairs of spaces always satisfy 
  the inf-sup condition \cite{Arnold2002} while combinations of FE spaces that are {\em not} stable
  are ruled out. For the above mentioned 1D wave equations, approximating the velocity with piecewise linear and the height field
  with piecewise constant spaces satisfies the requirements of FEEC and gives indeed a stable scheme
  (cf. again \cite{WaltersCarey1983} and Sect.~\ref{sec_standard}).
  
  Although providing a very general mathematical framework, naturally there are issues for which 
  FEEC yields no satisfying answers. Let us consider, for instance, problems in geophysical 
  fluid dynamics (GFD), in which an additional Coriolis term in the equations considers effects 
  caused by the earth's rotation \cite{PedloskyGFD1979}. For an atmosphere in rest, the Coriolis 
  force that depends on the velocity and the gradient of the pressure (or height) are in 
  {\em geostrophic balance}. To maintain this balance in the discretization, the 
  pressure (or height) field should be represented discretely at one order of consistency higher than that of the 
  velocity field. Unfortunately, this contradicts the requirement imposed by FEEC on this FE pair. 
  In \cite{Cotterlbb2009}, this issue could be resolved by applying a combination of FE and 
  Discontinuous Galerkin (DG) spaces.
  
  Moreover, in order to meet the regularity requirement of the chosen FE pairs, FEEC requires the PDEs 
  to be written in weak variational form, in which partial integration has been performed. As pointed 
  out recently in \cite{LeeWinther2016}, the conventional mixed (weak) form of the equations causes 
  certain operators, such as the co-derivative, to be non-local (global) operators. As a consequence, 
  such FEEC methods are not locally volume preserving, which reduces the quality of the local representation 
  of the quantities of interest (cf. \cite{LeeWinther2016}). 
  
  In this manuscript our main goal is to introduce a FE discretization framework that provides 
  an alternative methodology to avoid mentioned unsuitable FE choices with GFD in mind.
  More specifically, we develop a framework that applies two FEEC pairs 
  instead of one, therefore providing a larger variety of different combinations of FE spaces, 
  in which both derivatives and co-derivatives are local operators. This framework is based on formulating 
  the PDEs in split form, as introduced in \cite{BauerPhD:2013,Bauer2016} for the GFD equations.
  The split equations consist of a topological and metric part while employing 
  straight and twisted differential forms to adequately model the physical quantities 
  of interest. The FE discretization framework translates this geometrical structure 
  from the continuous to similarly structured discrete equations (cf. Sect.~\ref{sec_structure_waveeq}).

  Our approach shares some basic ideas with other discretization techniques, 
  in particular {\em mimetic discretizations} (see e.g. \cite{BossavitGFD2001,Bochev2006,Desbrun_Hirani_Leok_Marsden2005}, 
  and \cite{PalhaGerritsma2014} for a historical overview). There, the PDEs are also formulated by differential forms 
  and a clear distinction between purely topological and metric terms is achieved. 
  Applying algebraic topology as discrete counterpart to differential geometry, the discrete equations mimic 
  the underlying geometrical structure and are therefore denoted as structure-preserving (cf. \cite{KreefPalhaGerritsma2010}).
  Similar ideas of distinguishing between metric-dependent and metric-free terms in a GFD related context 
  can also be found in \cite{CotterThuburn2012} introducing FEEC discretizations of the nonlinear rotating 
  shallow-water equations. In spite of these similarities, none of the schemes associate
  a proper FE space to each variable, as suggested by our framework. 
  
  For the sake of a clear exposition, we focus on a simple example, the split 1D linear shallow-water set of equations. 
  Extending our framework to treat also more practically relevant equations, such as the nonlinear rotating shallow-water 
  equations, is subject of ongoing and future work. 
  By investigating structure-preserving methods that apply lowest order (piecewise linear and constant) FE spaces
  to keep computational costs low, we address the requirements of GFD in developing schemes that satisfy 
  first principle conservation laws (i.e mass and momentum conservation) and that are suited for simulations with 
  integration times in the order of years and longer.

  We structure the manuscript as follows. In Sect.~\ref{sec_structure_waveeq}, we introduce 
  the split set of 1D wave equations and motivate the use of the split form of equations as principle 
  formulation for their discretization. Recalling in Sect.~\ref{sec_standard} two 
  low-order mixed FE schemes, namely the unstable P1--P1 and the stable P1--P0 pairs, 
  we introduce in Sect.~\ref{sec_split} the new discretization framework, referred to as {\em split FE method}.
  We suggest a solving algorithm and present the schemes'
  discrete dispersion relations. Comparing them if possible with the 
  conventional mixed schemes, we perform in Sect.~\ref{sec_num_analysis} numerical
  simulations to investigate conservation behavior, convergence rates,
  and accuracy of all schemes. Finally in Sect.~\ref{sec_summary}, we draw conclusions and 
  provide an outlook for ongoing and future work.

 \section{Hierarchical structuring of the wave equation}
 \label{sec_structure_waveeq} 
 
 \newcommand{\LO}{\color{black}     \Lambda^0               \color{black}}  
  \newcommand{\LtwO}{\color{black}   \widetilde\Lambda^0     \color{black}} 
  \newcommand{\LOh}{\color{black}    \Lambda^0_h             \color{black}}
  \newcommand{\LtwOh}{\color{black}  \widetilde\Lambda^0_h   \color{black}}

  \newcommand{\LI}{\color{black}     \Lambda^1               \color{black}}  
  \newcommand{\LtwI}{\color{black}   \widetilde\Lambda^1     \color{black}} 
  \newcommand{\LIh}{\color{black}    \Lambda^1_h             \color{black}}
  \newcommand{\LtwIh}{\color{black}  \widetilde\Lambda^1_h   \color{black}}

  \newcommand{\LOhtest}{\color{black}    \hat\Lambda^0_h             \color{black}}
  \newcommand{\LIhtest}{\color{black}    \hat\Lambda^1_h             \color{black}}
 \newcommand{\LOtest}{\color{black}    \hat\Lambda^0             \color{black}}
  \newcommand{\LItest}{\color{black}    \hat\Lambda^1             \color{black}}
  \newcommand{\LtwIhtest}{\color{black}  \widehat{\widetilde\Lambda^1_h}   \color{black}}

 In order to convey the basic ideas to the reader most clearly and comprehensively, 
 we develop our FE discretization framework on a simple
 example, the set of 1D wave equations. As the split form of the equations (cf. \cite{Bauer2016})
 is an essential ingredient of our approach, we first introduce the corresponding
 split 1D wave equations and compare them with conventional formulations.

 \paragraph{A hierarchical structuring of the wave equation.}
  
 Compared to conventional formulations, the split 1D wave equations 
 arrange themselves in the following hierarchy of linear wave equations:
 \begin{enumerate}
  \item {\bf one} $2^{\rm nd}$-order equation: 
    \begin{equation}\label{equ_one_wave}     
      \frac{\partial^2 h(x,t)}{\partial t^2} - c^2 \frac{\partial^2 h(x,t)}{\partial x^2}  = 0, 
    \end{equation}
    with phase velocity $c = \omega/k$ for wave frequency $\omega$ and wave number $k$
    and with boundary condition (BC) $h(0,t) = h(L,t)\, \forall t$; 
  \item {\bf two} $1^{\rm st}$-order equations:
    \begin{equation}\label{equ_standard_wave}
      \frac{\partial u(x,t)}{\partial t} + g \frac{\partial h(x,t)}{\partial x} = 0, \qquad 
       \frac{\partial h(x,t)}{\partial t} + H \frac{\partial u(x,t)}{\partial x}  = 0, 
    \end{equation}
    with phase velocity $c = \sqrt{gH}$ and BCs $h(0,t) = h(L,t) , \, u(0,t) = u(L,t)\, \forall t$;
  \item {\bf two} $1^{\rm st}$-order equations and {\bf two} closure equations:
    \begin{equation}\label{equ_split_wave}
      \begin{split}
       \frac{\partial u^{(1)} }{\partial t}  + g {\bf d} h^{(0)}  = 0, & \qquad 
       \frac{\partial \widetilde h^{(1)} }{\partial t} + H {\bf d} \widetilde u^{(0)}  = 0, \\ 
     \widetilde u^{(0)}  = \widetilde\star u^{(1)} , & \qquad \widetilde h^{(1)}  = \widetilde\star h^{(0)} ,
    \end{split}      
    \end{equation}
    with phase velocity $c = \sqrt{gH}$ and BCs $ h^{(0)}(0,t) = h^{(0)}(L,t) , \, \widetilde u^{(0)}(0,t) = \widetilde u^{(0)}(L,t)\, \forall t$; 
 \end{enumerate}
 for a periodic domain $x \in [0,L]$. The smooth function $h(x,t)$ 
 denotes the height elevation with mean fluid height $H$ while assuming a trivial bottom topography.
 The smooth function $u(x,t)$ denotes the velocity of some fluid parcel in direction of $x$ 
 at time $t\in [0,T] \subset \mathbbm R$. $g$ is the gravitational acceleration. 
 Describing the momentum and height of a shallow water column, 
 equations~\eqref{equ_standard_wave} are usually referred to as shallow-water equations. 
 Throughout this manuscript however, we will refer to the preceding formulations 
 as wave equations. 
 

 The splitting of equations~\eqref{equ_split_wave} into topological (first line) and metric parts 
 (second line) is based on their formulation in terms of differential forms: 
 $u^{(1)}  \in (\LI;0,T)$ denotes a time-dependent straight 1-form and $h^{(0)} \in (\LO;0,T)$ a time-dependent straight 0-form 
 (straight function) -- straight forms do not change their signs when the orientation 
 of the manifold changes; $\widetilde u^{(0)}  \in (\LtwO;0,T)$ denotes a time-dependent twisted 0-form and 
 $\widetilde h^{(1)} \in (\LtwI;0,T)$ a time-dependent twisted 1-form -- twisted forms compensate the change in 
 signs which would result from a change in orientation. The latter property is shared by 
 the twisted Hodge-star operator $\widetilde \star: \Lambda^k \rightarrow \widetilde \Lambda^{(1 - k)}$ 
 (resp. $ \widetilde \Lambda^k \rightarrow  \Lambda^{(1 - k)}$) mapping from 
 straight (resp. twisted) $k$-forms to twisted (resp. straight) $(1-k)$-forms,
 here for $k=0,1,$ as we are in one dimension. The index $^{(k)}$ denotes the degree, 
 and $\Lambda^k,\widetilde \Lambda^k$ the space of all $k$-forms.
 
 The exterior derivative is defined as the map ${\bf d}: \LO \rightarrow \LI$ 
 (or ${\bf d}: \LtwO \rightarrow \LtwI$ for twisted forms, cf. diagram~\eqref{equ_spaces_cd}).
 It reduces in one dimension to the total derivative of a smooth function $f(x)$,
 i.e. for $f^{(0)} \in \LO, {\bf d}f^{(0)} = \partial_x f(x) dx \in \LI$.  
 The map of $\bf d$ on twisted forms is defined analogously. 
 
 We refer the reader to Sec.~\ref{sec_varform_spliteq} for an explicit representation of
 the split wave equations~\eqref{equ_split_wave} in the local coordinate $x \in [0,L]$, and to \cite{Bauer2016} for an elaborated discussion 
 about straight and twisted forms, the twisted Hodge-star operators, 
 and the derivation of $n$-dimensional split equations of GFD.

 \paragraph{The split wave equations:}
 
 We argued in detail in \cite{Bauer2016} why it is favorable to represent the equations of GFD 
 in terms of differential forms rather than vector-calculus notation, and we discussed advantages 
 of the split form in particular. Here, the split representation of the wave equations fits naturally 
 into the arrangement of wave equations presented above, in which the saddle point (mixed) form~\eqref{equ_standard_wave} 
 results from a splitting of the standard form~\eqref{equ_one_wave}, and the split form~\eqref{equ_split_wave} from 
 a further splitting of the saddle point form~\eqref{equ_standard_wave}.
 This somehow suggests to study in more detail possible benefits of a discretization based on 
 the split form \eqref{equ_split_wave} rather than on the mixed form \eqref{equ_standard_wave}, 
 just as there exist cases in which discretizations based on equations~\eqref{equ_standard_wave}, 
 usually discretized by {\em mixed} FE methods, perform better than schemes 
 resulting from {\em standard} FE approaches that rely on formulation~\eqref{equ_one_wave}.

 Let us first elaborate on the latter point, for which we present two examples where 
 standard FE schemes even fail to provide correct solutions. As pointed out in \cite{ArnoldFEEC2010}, 
 the standard FE solution to the vector Poisson equation on non-convex polyhedral domains will 
 converge to a false solution of the problem for almost all forcing functions $f$.
 The same happens when calculating standard FE solutions for the vector Laplacian on an annulus. 
 However, correct solutions can be obtained when using mixed formulations, in particular those based 
 on FEEC \cite{ArnoldFEEC2010}.
 In addition, mixed formulations, in which two or even more FE spaces are used to approximate separate 
 variables, are well suited for saddle point problems. These arise, for instance, in constrained minimization 
 problems (Stokes equations or Darcy flow), in which the elimination of one of the variables is not 
 possible (cf. \cite{FenicsBook}, Chapter 36). In this latter case, equations in form~\eqref{equ_one_wave} 
 do even not exist. 
 There are also some physical and numerical reasons in favor of mixed FE methods. According to \cite{LeVeque2002}, 
 it is the first-order system, such as \eqref{equ_standard_wave}, that follows from physics 
 (i.e. from first conservation principles of mass and momentum) and not the second-order equation, 
 such as \eqref{equ_one_wave}. In particular, efficient numerical schemes can be derived more easily
 from the first-order system.

 Following in particular the latter line of argumentation, the split form more accurately models the 
 physical properties than the first-order system, because system~\eqref{equ_split_wave} 
 provides for each variable the adequate (straight or twisted) differential form that suits in dimension 
 and orientation the corresponding property of a real fluid (cf. \cite{Bauer2016}, Section 7). 
 It seems appropriate during the discretization process to provide for each variable a suitable 
 FE space too. As a proof of concept, one aim of this manuscript is to exploit this additional freedom 
 in the choice of FE spaces and study, in terms of convergence, accuracy, stability and dispersion relation, 
 possible advantages, but also drawbacks, that come along with this generalized FE method.

 More precisely, by using the split form of the equations, we approximate each variable by a FE space 
 such that the discrete version of the exterior derivative $\bf d$ satisfies the mappings 
 $\Lambda_h^0 \rightarrow \Lambda_h^1$ and $\widetilde \Lambda_h^0 \rightarrow \widetilde\Lambda_h^1$. 
 $ \Lambda_h^i,\widetilde \Lambda_h^i, i = 0,1$, are FE spaces that approximate the corresponding 
 continuous spaces. As pointed out later in more detail, this approach guarantees correct discrete 
 topological equations up to projection error caused by projecting the continuous equations into the FE spaces. 
 Any additional errors are caused by the metric-dependent closure equations that are realized by nontrivial 
 projections between straight and twisted spaces (cf. Sect.~\ref{sec_split}). Such a clear separation 
 between projection and additional errors, the latter caused by partial integrations in the weak form, 
 is not obvious in mixed FE methods. Part of future work will be to study how these different error sources relate.

 \newcommand{\indxVOtest} {  {\color{black}  {l'}      \color{black}} } 
 \newcommand{\indxVOtrial}{  {\color{black}  {l}     \color{black}} } 
 \newcommand{\indxVItest} {  {\color{black}  {m'}      \color{black}} } 
 \newcommand{\indxVItrial}{  {\color{black}  {m}     \color{black}} }

 \newcommand{\VO}{\color{black}   V^1      \color{black}} 
 \newcommand{\VOh}{\color{black}  V^1_h    \color{black}}
 \newcommand{\VI}{\color{black}   V^0      \color{black}} 
 \newcommand{\VIh}{\color{black}  V^0_h    \color{black}}
 \newcommand{\VtwO}{\color{black}   \widetilde V^1      \color{black}} 
 \newcommand{\VtwOh}{\color{black}  \widetilde V^1_h    \color{black}}
 \newcommand{\VtwI}{\color{black}   \widetilde V^0      \color{black}} 
 \newcommand{\VtwIh}{\color{black}  \widetilde V^0_h    \color{black}}
 \newcommand{\VOtest}{\color{black}   \hat V^1      \color{black}} 
 \newcommand{\VOhtest}{\color{black}  \hat V^1_h    \color{black}}
 \newcommand{\VItest}{\color{black}   \hat V^0      \color{black}} 
 \newcommand{\VIhtest}{\color{black}  \hat V^0_h    \color{black}}

 \newcommand{\fctVOtest}{\color{black}  \phi_{l'}(x)     \color{black}} 
 \newcommand{\fctVOtrial}{\color{black}  \phi_{l}(x)     \color{black}}
 \newcommand{\fctVItest}{\color{black}  \chi_{m'}(x)     \color{black}} 
 \newcommand{\fctVItrial}{\color{black}  \chi_{m}(x)     \color{black}}
 
 \newcommand{\uVOtrial}{\color{black}  u_{l}     \color{black}} 
 \newcommand{\hVOtrial}{\color{black}  h_{l}     \color{black}} 
 
 \newcommand{\uVItrial}{\color{black}  u_{m}     \color{black}} 
 \newcommand{\hVItrial}{\color{black}  h_{m}     \color{black}} 
 
 \newcommand{\utwVOtrial}{\color{black}  \widetilde{u}_{l}     \color{black}} 
 \newcommand{\htwVOtrial}{\color{black}  \widetilde{h}_{l}     \color{black}} 
 
 \newcommand{\utwVItrial}{\color{black}  \widetilde{u}_{m}     \color{black}} 
 \newcommand{\htwVItrial}{\color{black}  \widetilde{h}_{m}     \color{black}}

 \newcommand{\massmxnn}{ {\bf M}^{nn} } 
 \newcommand{\stiffmxnn}{ {\bf D}^{nn} }
 \newcommand{\massmxee}{ {\bf M}^{ee} } 
 
 \newcommand{\stiffmxen}{ {\bf D}^{en} } 
 \newcommand{\stiffmxne}{ {\bf D}^{ne} }
 
  \newcommand{\massmxne}{ {\bf M}^{ne}  } 
  \newcommand{\massmxen}{ {\bf M}^{en}  }  
 
 \newcommand{\prjcmxne}{{\bf P}^{ne} }

 \newcommand{\velvctn}{\color{black}{\bf u}_{n} \color{black}}
 \newcommand{\hvctn}{\color{black}{\bf h}_{n} \color{black}} 
 \newcommand{\velvcte}{\color{black}{\bf u}_{e} \color{black}}
 \newcommand{\hvcte}{\color{black}{\bf h}_{e} \color{black}}
 \newcommand{\metrivcte}{\color{black} {\bf \Delta x}_{e} \color{black}}

 \section{The mixed finite element method}
 \label{sec_standard}

 \begin{figure}[t]\centering
      \includegraphics[scale=0.9]{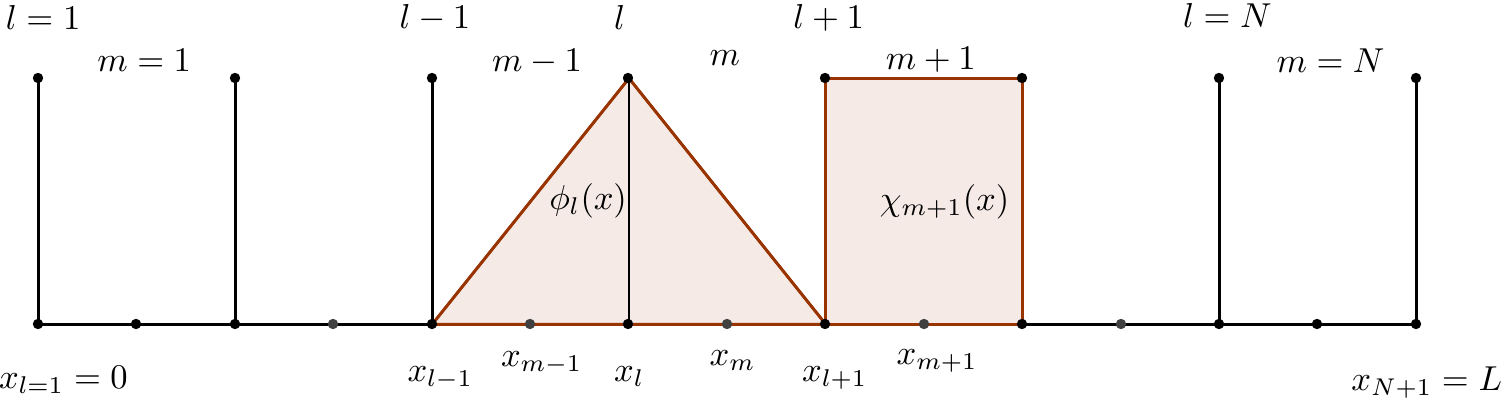}
  \caption{Discretization of the periodic domain $x \in [0,L]$ by $N$ subintervals (elements) 
          and $N$ independent nodes; the former are denoted by $m$, the latter by $l$ for $m,l = 1,..,N$.
          $ \fctVOtrial$ are the basis functions of $\VOh$ and  $ \fctVItrial$ of $\VIh$ (cf. 
          Appendix~\ref{append_int_values} for their definition).
          }  
   \label{fig_mesh}
 \end{figure}

 In this section, we briefly review two examples of the general idea behind mixed FE methods 
 and we introduce notations and definitions required for the remainder of the paper. 
 First, we observe that there are two equations for two prognostic variables in 
 equation~\eqref{equ_standard_wave}. In order to discretize them one has the choice of mixing 
 the corresponding FE spaces, which gives rise to the naming convention. Here we consider (i) 
 similar $(\VOh,\VOh,)$ or (ii) mixed FE spaces $(\VOh,\VIh)$, where $\VOh$ refers to a linear 
 and $\VIh$ to a constant FE space.


 In more detail, we subdivide the domain $[0,L]$ into $N$ subintervals $[x_l,x_{l+1}]$, 
 for $l= 1,\dots N,$ as shown in Figure~\ref{fig_mesh}. We denote the points $x_l$ as nodes $l$ 
 and the subintervals as elements $m$ with centers $x_m:= \frac{1}{2}(x_l + x_{l+1})$ and element sizes 
 $\Delta x_m = (x_{l+1} - x_l)$, for $l,m=1,\dots N$. 
 Basis functions and coefficients at nodes carry the subindices $l$ and $l'$ and those at element centers 
 the subindices $m$ and $m'$.
 To impose periodic boundary conditions, we identify the nodes at $x=0$ and at $x=L$ 
 with each other and consider only 1 Degree of Freedom (DoF) associated to this node. 
 We have thus $N$ independent DoF for both nodes and elements. 
 
 Throughout this manuscript we will use the Lagrange functions $\{ \fctVOtrial  \}_{\indxVOtrial = 1}^{N}$
 as basis for the piecewise linear FE space $\VOh$ and the 
 step functions $\{ \fctVItrial  \}_{\indxVItrial = 1}^{N}$ as basis for the piecewise
 constant FE space $\VIh$ (cf. Appendix~\ref{append_int_values} for their definitions).

 \subsection{The unstable P1--P1 finite element discretization}
  \label{sec_standard00}

  In order to find approximate solutions of the wave equations~\eqref{equ_standard_wave},
  we apply Galerkin's method. To this end, we formulate the wave equations in variational form 
  and approximate $u(x,t)$ and $h(x,t)$ by suitable discrete functions.
  Being aware of the problems of spurious modes (discussed below), we start with 
  a naive choice of a pair of FE spaces approximating both $u(x,t)$ and $h(x,t)$ 
  by piecewise linear (P1) functions $u_h(x,t)\in (\VOh;0,T)$ and $h_h(x,t)\in(\VOh;0,T)$.
  We refer to this scheme as P1--P1 scheme.

  The discrete variational form of the wave equations is then given by: 
  find $(u_h(x,t), h_h(x,t)) \in (\VOh;0,T) \times (\VOh;0,T)$, such that
  \begin{equation}
   \begin{split}\label{equ_var_form00}
    \int_L  \left( \frac{\partial u_h(x,t)}{\partial t} + g \frac{\partial h_h(x,t)}{\partial x}  \right) \hat\phi(x) dx 
    & = 0 , \quad \forall  \hat\phi(x) \in \VOhtest \subset \VO ,  \qquad \\
     \int_L  \left( \frac{\partial h_h(x,t)}{\partial t} + H \frac{\partial u_h(x,t)}{\partial x}  \right) \hat\phi(x) dx 
    & = 0 , \quad \forall  \hat\phi(x)  \in \VOhtest \subset \VO ,
   \end{split}
  \end{equation}
  for any test function $\hat\phi(x)$ in the test space $\VOhtest$. 
  Expanding the variables in terms of the basis of the trial space $\VOh$, 
  \begin{equation}
   {u}_h(x,t) \approx  \sum_{\indxVOtrial = 1}^{N} \uVOtrial(t) \fctVOtrial \, \quad \text{and} \quad 
   {h}_h(x,t) \approx  \sum_{\indxVOtrial = 1}^{N} \hVOtrial(t) \fctVOtrial \, ,
  \end{equation}
  with time dependent coefficients $\uVOtrial(t)$ and $\hVOtrial(t)$, and varying the test 
  functions $\hat\phi(x)$ over the basis functions $\{ \fctVOtest  \}_{\indxVOtest = 1}^{N}$ 
  of the test space $\VOhtest$, we obtain
  \begin{equation}\label{equ_stnd_1}
   \begin{split}
           \sum_{\indxVOtrial = 1}^{N} \partial_t  \uVOtrial(t) \int_L \fctVOtrial \fctVOtest  dx 
       + g \sum_{\indxVOtrial = 1}^{N} \hVOtrial(t) \int_L \frac{d\fctVOtrial}{dx} \fctVOtest dx 
       =     0    ,    \quad \text{\ for\ } \indxVOtest = 1,\dots,N  , \\
           \sum_{\indxVOtrial = 1}^{N} \partial_t  \hVOtrial(t) \int_L \fctVOtrial  \fctVOtest  dx 
       + H \sum_{\indxVOtrial = 1}^{N} \uVOtrial(t) \int_L \frac{d\fctVOtrial}{dx} \fctVOtest dx 
       =     0   ,     \quad \text{\ for\ }\indxVOtest  = 1,\dots,N  . 
   \end{split}
  \end{equation}

 This discrete system can analogously be written in matrix-vector formulation. Defining the 
 coefficient vectors $\velvctn = (u_1(t),...,u_l(t),...,u_{N}(t))$ and $\hvctn = (h_1(t),...,h_l(t),...,h_{N}(t))$, 
 with subindex $n$ for nodes, we obtain the following linear system of algebraic equations:
 \begin{align}\label{equ_matrix_inst}
 \massmxnn  \frac{\partial \velvctn }{\partial t}  + g \stiffmxnn \hvctn   = 0 , \quad 
 \massmxnn  \frac{\partial \hvctn   }{\partial t}  + H \stiffmxnn \velvctn = 0 ,
 \end{align}
 in which $\massmxnn$ is a $(N\times N)$ mass-matrix with metric-dependent coefficient 
 \begin{equation}\label{equ_mssmxnn_coeff}
  M_{\indxVOtrial\indxVOtest} = \int_L \fctVOtrial \fctVOtest  dx ,  
 \end{equation}
 and $\stiffmxnn $ a $(N\times N)$ stiffness-matrix with metric-independent coefficients
 \begin{equation}\label{equ_stiffxnn_coeff}
  D_{\indxVOtrial\indxVOtest} = \int_L \frac{d \fctVOtrial  }{dx} \fctVOtest dx \, .
 \end{equation}
 In Appendix~\ref{append_int_values}, we provide explicit 
 representations of these matrices for the mesh of 
 Figure~\ref{fig_mesh} subject to periodic boundary conditions.

 \paragraph{Time discretization.} 
 To discretize the semi-discrete equations~\eqref{equ_matrix_inst} in 
 time, we use the Crank-Nicolson (CN) scheme, a symmetric, implicit time discretization method. 
 Then, for the time $t \in [0,T]$ and a time step size of $\Delta t$, the full discretized 
 matrix-vector equations read
 \begin{align}\label{equ_matrix_inst_CN}
      \massmxnn \velvctn^{t+1} &  = \massmxnn  \velvctn^{t} +  \frac{1}{2}\Delta t g \stiffmxnn \hvctn^{t}   + \frac{1}{2}\Delta t g \stiffmxnn \hvctn^{t+1} ,\\
      \massmxnn \hvctn^{t+1}   &  = \massmxnn  \hvctn^{t}   +  \frac{1}{2}\Delta t H \stiffmxnn \velvctn^{t} + \frac{1}{2}\Delta t H \stiffmxnn \velvctn^{t+1} .
 \end{align}
 We will solve these equations by fixed point iteration, even though they could be solved 
 directly due to the linearity. The motivation for this inefficient algorithmic choice lies 
 in future work on non-linear equations and the possibility to treat possible kernels for 
 the discrete Hodge-star operators in the split schemes (cf. \eqref{equ_augmentedsystem}).
 For more details on the solver, we refer to Sect.~\ref{sec_split_CN} in which we introduce 
 a corresponding algorithm for the split finite element schemes.

 The choice of using an iterative algorithm does not alter the results discussed in this paper. 
 Only the time step size has to be restricted, in spite of the unconditional stability of the CN scheme.
 Similarly to the CFL number in explicit schemes, $\Delta t$ must not exceed a certain 
 threshold $\mu$ in order to guarantee that even the fastest waves are sufficiently well 
 resolved and that the iterative solver converges. Defining the CFL number for the 1D wave equations by 
 \begin{equation}\label{cflnumber}
  \mu = \sqrt{g H} \frac{\Delta t}{\Delta x} \, ,
 \end{equation}
 we restrict $\Delta t$ for a given $\Delta x$ such that $\mu$ is smaller that some constant. 
 This constant, determined empirically for all test cases and mesh resolutions applied in this manuscript, 
 is given for the unstable scheme by $\mu\leq 1.15$ (see Sect.~\ref{sec_split_CN} for more details).

 \paragraph{Discussion about spatial stability.} 
 
 Throughout this paper, we study the stability properties of each 
 discrete scheme by investigating whether it supports {\em spurious modes}.  
 Such spurious modes are free modes, i.e. unphysical, very oscillatory small-scale waves on the grid scale 
 that are not constrained by the discrete equations. In case of nonlinear 
 equations, these modes get coupled to the smooth solution and grow fast, 
 making the scheme unusable. Besides studying the inf-sup condition (see, e.g. \cite{BoffiBrezzi2013}), 
 one way of investigating the occurrence of such 
 spurious modes is to determine the dispersion relation of the discrete scheme. 
 
 For the unstable P1--P1 FE scheme, we associate spurious modes with those 
 modes that show zero group velocity $\omega$, for some wave numbers $k$. These modes are 
 trapped to the grid scale and may lead to the just described instabilities.
 More precisely, the angular frequency $\omega = \omega(k)$ satisfies the dispersion 
 relation (see Appendix~\ref{app_disc_disp_rel})
 \begin{equation}\label{equ_disprel_P1P1}
  c_d = \frac{\omega}{k} =  \pm \sqrt{gH} \frac{\sin (k\Delta x)}{k \Delta x} \left[ \frac{3  }{ 2 + \cos (k\Delta x)  } \right].
 \end{equation}
 For small $k\rightarrow 0$, the discrete dispersion relation converges to the 
 continuous one with $c_d \rightarrow c = \sqrt{gH}$. 
 However, besides the correct zero frequency at $k =0$, this dispersion relation has 
 at the shortest wave length $k = \frac{\pi}{\Delta x}$ a second zero solution (spurious mode)
 with zero wave speed leading to an unphysical standing wave. 
 

 \subsection{The stable P1--P0 finite element discretization}
  \label{sec_standard01}
 
 The instability demonstrated in the previous subsection, caused by a choice of similar FE spaces 
 for both variables, is a general feature of methods with same FE spaces for approximating $u$ and $h$ 
 (cf.~\cite{CotterMcRae2014}). To avoid this problem, we present a stable mixed FE method using 
 different FE spaces to approximate velocity and height. 
 
 As proposed by other authors (e.g.~\cite{LeRoux2007}), we approximate $u(x,t)$ by the 
 piecewise linear (P1) function $u_h(x,t)\in (\VOh;0,T)$ and $h(x,t)$ by the piecewise constant (P0)
 function  $h_h(x,t) \in (\VIh;0,T)$ and we apply Galerkin's method to find these approximations. 
 We refer to this scheme as P1--P0 scheme.
 
 The formulation of the wave equations~\eqref{equ_standard_wave} in variational form implies 
 multiplication of the momentum equation by a test function $\hat\phi(x) \in \VOhtest$ and of 
 the continuity equation by a test function $\hat\chi(x) \in \VIhtest$. 
 Because the derivative of $h_h(x,t)$ is not well-defined globally, the corresponding term in the momentum 
 equation requires integration by parts, whereas for the remaining terms 
 trivial projections into the corresponding spaces are suitable. 
 
 The discrete variational/weak form of the wave equations is thus given by:
 find $(u_h(x,t), h_h(x,t)) \in (\VOh;0,T) \times (\VIh;0,T)$, such that
 \begin{equation} 
 \begin{split}\label{equ_var_form01} 
   \int_L  \left( \frac{\partial u_h(x,t)}{\partial t} \hat\phi(x) - g h_h(x,t) \frac{\partial \hat\phi(x) }{\partial x}  \right) dx 
   & = 0, \quad \forall \hat\phi(x) \in \VOhtest \subset \VO , \\
   \int_L  \left( \frac{\partial h_h(x,t)}{\partial t}      +  H \frac{\partial u_h(x,t)}{\partial x}  \right) \hat\chi(x) dx 
   & = 0, \quad \forall \hat\chi(x) \in \VIhtest \subset \VI. 
 \end{split}
 \end{equation}
 The momentum equation is given in {\em weak} form, due to the integration by parts, while the continuity 
 equations remains in variational form. The boundary terms from integration by parts in the weak momentum 
 equation vanish because of the periodic boundary conditions and continuity conditions at cell interfaces. 
 This special treatment is a source of additional 
 error, while the trivially projected terms do not introduce further errors beyond those caused 
 by the approximation of the initial conditions by the FE spaces \cite{CotterMcRae2014}.

 Expanding the variables in terms of the corresponding bases yields
 \begin{equation}
     {u}_h(x,t) \approx \sum_{\indxVOtrial = 1}^{N} \uVOtrial(t) \fctVOtrial \quad \text{and} \quad 
     {h}_h(x,t) \approx \sum_{\indxVItrial = 1}^{N} \hVItrial(t) \fctVItrial.
 \end{equation}
 Substituting these expansions into equations~\eqref{equ_var_form01} gives 
 \begin{equation}\label{equ_P0P1}
	\begin{split}
	 \sum_{\indxVOtrial = 1}^{N} \partial_t \uVOtrial(t) \int_L \fctVOtrial \fctVOtest  dx 
     - g \sum_{\indxVItrial = 1}^{N} \hVItrial(t)  \int_L \fctVItrial \frac{d\fctVOtest}{dx}  dx 
     & =   0  ,    \quad \text{\ for\ } \indxVOtest = 1,\dots,N , \\
	 \sum_{\indxVItrial = 1}^{N} \partial_t \hVItrial(t)  \int_L \fctVItrial \fctVItest  dx 
     + H \sum_{\indxVOtrial = 1}^{N} \uVOtrial(t) \int_L \frac{d \fctVOtrial}{dx} \fctVItest   dx  
     & =     0  ,   \quad \text{\ for\ } \indxVItest = 1,\dots,N .
	\end{split}
      \end{equation}

 Analogously to equation~(\ref{equ_stnd_1}) this discrete system can be written in matrix-vector formulation. For the 
 coefficient vectors $\velvctn = (u_1(t),...,u_l(t),...,u_{N}(t))$ with subindex $n$ for nodes 
 and $\hvcte = (h_1(t),...,h_m(t),...,h_{N}(t))$ with subindex $e$ for elements,
 we obtain the following linear system 
  \begin{align}\label{equ_matrix_stbl}
 \massmxnn   \frac{\partial \velvctn }{\partial t}  - g \stiffmxne  \hvcte   = 0 , \quad 
 \massmxee   \frac{\partial \hvcte   }{\partial t}  + H \stiffmxen  \velvctn = 0 .
 \end{align}
 $\massmxnn$ is the mass-matrix from above with metric-dependent coefficients~\eqref{equ_mssmxnn_coeff}, 
 $\massmxee$ is a $(N\times N)$ mass-matrix with metric-dependent coefficients
 \begin{equation}\label{equ_mssmxee_coeff}
  M_{\indxVItrial\indxVItest} = \int_L \fctVItrial \fctVItest  dx ,
 \end{equation}
 and $\stiffmxne $ and $\stiffmxen  $ are $(N\times N)$ stiffness-matrices. 
 Recall that in 1D with periodic BC, the number of nodes and the number of 
 elements coincide, therefore all matrices are of size $N \times N$. 
 The metric-independent coefficients of $\stiffmxen$ are defined by 
 \begin{equation}\label{equ_stiffxne_coeff}
  D_{\indxVOtrial\indxVItest} = \int_L \frac{d \fctVOtrial  }{dx} \fctVItest dx \, ,
 \end{equation}
 and those of $\stiffmxne $ via the property $\stiffmxne = (\stiffmxen)^T$, in which $()^T$ denotes 
 the transpose of a matrix. We refer to Appendix~\ref{append_int_values} to see the full
 matrices for the example given in Figure~\ref{fig_mesh}.

 \paragraph{Time discretization.} 
 
 To discretize the semi-discrete equations~\eqref{equ_matrix_stbl} in 
 time, we use the Crank-Nicolson scheme again. For the time $t \in [0,T]$ 
 and a time step size of $\Delta t$, the fully discretized 
 matrix-vector equations read
 \begin{align}\label{equ_matrix_stbl_CN}
      \massmxnn  \velvctn^{t+1} &  =  \massmxnn  \velvctn^{t} -  \frac{1}{2}\Delta t g \stiffmxne \hvcte  ^{t} - \frac{1}{2}\Delta t g \stiffmxne \hvcte  ^{t+1} ,\\
      \massmxee   \hvcte^{t+1} &  =   \massmxee  \hvcte  ^{t} +  \frac{1}{2}\Delta t H \stiffmxen \velvctn^{t} + \frac{1}{2}\Delta t H \stiffmxen \velvctn^{t+1} .
  \end{align}
 We solve this implicit system of equations iteratively by fixed point iteration (cf. Sect.~\ref{sec_split_CN} for more details).

 Here, we restrict the time step size such that for all meshes and test cases studied
 the CFL number~\eqref{cflnumber} satisfies $\mu\leq 0.57 $. This is only half the value 
 compared to the unstable P1--P1 scheme, since the lower order representation 
 of the height field leads effectively to a coarsened spatial resolution.

 \paragraph{Discussion about spatial stability.}
 
 We briefly discuss the stability of the method by studying it's discrete dispersion relation. 
 As shown in Appendix~\ref{app_disc_disp_rel}, $\omega = \omega(k)$
 satisfies the dispersion relation
    \begin{equation}\label{equ_disprel_P0P1}
       c_d =  \frac{\omega}{k} = \pm \sqrt{gH}  \frac{\sin (k \frac{\Delta x}{ 2} ) } { k \frac{\Delta x}{ 2}} \left[  \frac{ 3}{ 2 + \cos (k\Delta x)  } \right]^{1/2} \, .
    \end{equation}
 The discrete wave speed $c_d$ converges for $k \rightarrow 0$ to the analytical 
 wave speed $c$ also for the P1--P0 case.
 However, this relation has only one root for $\omega$ at $k=0$, 
 while it shows a good approximation to the continuous dispersion relation
 at the shortest wave length $k = \frac{\pi}{\Delta x}$, in contrast to the unstable P1--P1 scheme 
 (cf. Fig.~\ref{fig_disp_relation}). Hence, this scheme exhibits no spurious modes and is therefore stable.

 \section{The split finite element method}
 \label{sec_split}

  In this section, we introduce a new finite element discretization method that uses the 
  split formulation of the wave equations~\eqref{equ_split_wave} (cf. \cite{Bauer2016} for 
  a more general formulation of such split equations). It is 
  referred to in the following as {\em split finite element (FE) method}.
  
  The split FE method is based on the following general ideas; it 
  \begin{itemize}
    \item discretizes the split form of the equations and keeps the splitting into 
          topological and metric equations preserved during the discretization process; 
    \item treats straight and twisted variables independently and provides for 
          each of them suitable FE spaces such that the topological momentum and continuity equations 
          can be written in variational form without partial integration; 
    \item provides FE spaces that form a pair of straight and twisted 
          cohomology chains, in which $\mathbf{d}$ maps between either straight or twisted chain elements
          according to diagram~\eqref{equ_spaces_cd}; 
    \item represents the metric closure equations (Hodge-star operators) as
          projections between the straight and twisted FE spaces (cf. diagram~\eqref{equ_spaces_cd}). 
  \end{itemize}
  Following these algorithmic ideas, the discretization of the topological equations is a trivial projection 
  into the corresponding FE spaces with only the projection error occurring. 
  Additional errors are introduced by metric-dependent closure equations, 
  in which the discretization of the continuous Hodge-star operators is 
  performed by non-trivial projections between straight and twisted spaces.

  \subsection{Variational formulation of the split wave equation}
  \label{sec_varform_spliteq}
     
  Starting from the split wave equations~\eqref{equ_split_wave}, we seek a discretization 
  by standard FE techniques. In particular,
  we want to apply Galerkin's method, which implies to formulate the split wave equations 
  in variational form.
  
  Employing the notation of differential forms, we notice that the time-dependent differential 
  forms in \eqref{equ_split_wave} can be represented in the 
  local coordinate $x \in [0,L]$ with time parameter $t \in [0,T]$. Denoting with $dx$ the dual basis to 
  the tangential basis $\partial_x:= \frac{\partial}{\partial x}$, 
  the 1-forms are given by $u^{(1)} = u^{(1)} (x,t) dx \in (\LI;0,T)$ and $\widetilde h^{(1)} = \widetilde h^{(1)} (x,t) dx \in (\LtwI;0,T)$.
  The 0-forms can be written as $h^{(0)} = h^{(0)} (x,t) \in (\LO;0,T)$ and $\widetilde u^{(0)} = \widetilde u^{(0)} (x,t) \in (\LtwO;0,T)$. 
  To distinguish the coordinate functions from the respective forms, the former carry the notation $(x,t)$. 
  The exterior derivative ${\bf d}:\LO \rightarrow \LI$ and ${\bf d}:\LtwO \rightarrow \LtwI$ 
  maps the corresponding 0-forms to the 1-forms ${\bf d} h_h^{(0)} = \partial_x h_h^{(0)} (x,t) dx$ and 
  ${\bf d}\widetilde u_h^{(0)} = \partial_x \widetilde u_h^{(0)} (x,t) dx$. 
  
  Furthermore, this coordinate representation allows us to describe the action
  of $\widetilde \star$ more precisely. It is defined via its action on the dual basis $dx$ and the constant
  function $1$ describing the unit volume. 
  In other words, denoting a choice of orientation as direct frame 'Or' and it's opposite 
  orientation as skew frame '-Or', then $\widetilde \star$ maps the straight 1-form 
  $dx \in \LI$ to the twisted 0-from $\widetilde \star dx = \{1 \ \text{in Or}, - 1 \ \text{in -Or}\} \in \LtwI$
  or the straight function $1\in \LI$ to the twisted 1-form 
  $\widetilde \star 1 = \{dx \ \text{in Or}, - dx \ \text{in -Or}\} \in \LtwI$. 
  For both the straight and the twisted Hodge-star operators a self-adjoint property holds 
  (e.g.\ $\widetilde \star \widetilde\star =  Id$. cf. \cite{Bauer2016} for more details).

  The idea behind the definition of twisted forms and the twisted Hodge-star is to guarantee that 
  the equations remain the same in both direct and skew frames.   
  To enhance readability, we will therefore assume for the remainder of the manuscript 
  to be in the direct frame 'Or' in which $\widetilde \star$ maps to the positive valued forms 
  and skip an analogous derivation for the case of '-Or'. However, we will keep the notation $\widetilde\,$
  to distinguish between the various spaces. Using this local representation, 
  the metric closure equations in \eqref{equ_split_wave} read:

  $ \widetilde u^{(0)} = \widetilde\star u^{(1)} = u^{(1)}(x,t) \widetilde\star dx =  u^{(1)}(x,t) $ \ and \ 
  $ \widetilde h^{(1)} = \widetilde\star h^{(0)} = h^{(0)}(x,t) \widetilde\star 1 =   h^{(0)}(x,t)dx$.

  \paragraph{Split variational form of the split 1D wave equations.}
  Using these local representations in \eqref{equ_split_wave} and multiplying the resulting 
  equations by a test function $\hat\chi(x) \in \LItest$, we obtain the variational form for the 
  topological equations: 
  find $(u^{(1)}(x,t), h^{(0)}(x,t),\widetilde h^{(1)}(x,t),\widetilde u^{(0)}(x,t)) \in (\LI;0,T) \times (\LO;0,T) \times (\LtwI;0,T) \times (\LtwO;0,T) $
  such that
  \begin{equation}
    \begin{split}\label{equ_varsplit_topo}
        \int_L \left( \frac{\partial u^{(1)}(x,t)}{\partial t} + g \frac{\partial h^{(0)}(x,t)}{\partial x} \right) \hat\chi(x) dx 
    & = 0 \, , \quad  \forall \hat\chi(x) \in \LItest ,  \qquad \\
    \int_L \left( \frac{\partial \widetilde h^{(1)}(x,t)}{\partial t} + H \frac{\partial \widetilde u^{(0)}(x,t)}{\partial x}  \right) \hat\chi(x) dx 
    & = 0 \, , \quad  \forall \hat\chi(x) \in \LItest ,
   \end{split}
   \end{equation}
  subject to the variational metric equations
  \begin{equation}\label{equ_varsplit_metri} 
  \begin{split}
         \int_L     u^{(0)}(x,t)  \hat\tau^i(x)  dx & =  \int_L   u^{(1)}(x,t) \hat\tau^i(x)  dx  \, , 
         \quad        \forall \hat\tau^i(x)  \in \hat\Lambda^{i}, i = 0,1,  \\ 
          \int_L    \widetilde h^{(1)}(x,t) \hat\tau^{j}(x)   dx & =  \int_L    h^{(0)}(x,t) \hat\tau^{j}(x)    dx \, ,  
          \quad        \forall \hat\tau^{j}(x)  \in \hat\Lambda^{{j}}, {j} = 0,1. 
  \end{split}  
  \end{equation}
  The latter equations follow by multiplying the local representation of the metric equations in \eqref{equ_split_wave}
  with the test functions $\hat\tau^{i,j}(x), {i,j} = 0,1,$ that can be elements of either $\LOtest$ or $\LItest$. 
  
  Having four equations for four unknowns (the two straight and two twisted ones),
  the system is closed. We will refer to this set of equations as 
  {\em split variational form} of the split 1D wave equations. 
  
  \newcommand{\starULih}{\color{black}\widetilde\star^{\, u}_{1-i}\color{black}}  
  \newcommand{\starHLjh}{\color{black}\widetilde\star^{\, h}_{1-j}\color{black}}
  \newcommand{\starULOh}{\color{black}\widetilde\star^{\,u}_{0}\color{black}}  
  \newcommand{\starHLOh}{\color{black}\widetilde\star^{\,h}_{0}\color{black}}
  \newcommand{\starULIh}{\color{black}\widetilde\star^{\,u}_{1}\color{black}}
  \newcommand{\starHLIh}{\color{black}\widetilde\star^{\,h}_{1}\color{black}}
  \renewcommand{\starULOh}{\color{black}\widetilde\star^{\,u}_{1}\color{black}}  
  \renewcommand{\starHLOh}{\color{black}\widetilde\star^{\,h}_{1}\color{black}}
  \renewcommand{\starULIh}{\color{black}\widetilde\star^{\,u}_{0}\color{black}}
  \renewcommand{\starHLIh}{\color{black}\widetilde\star^{\,h}_{0}\color{black}}

  \subsection{Discrete split variational formulation }
  \label{sec_splitFE_derivate}

  To discretize the split variational form \eqref{equ_varsplit_topo} and \eqref{equ_varsplit_metri} 
  with Galerkin's method, we approximate the variables 
  $u^{(1)}, h^{(0)},\widetilde h^{(1)}, \widetilde u^{(0)}$ by suitable discrete $k$-forms 
  $u_h^{(1)}, h_h^{(0)},\widetilde h_h^{(1)}, \widetilde u_h^{(0)}$, respectively.  
  Following the above general ideas of the split FE method, the discrete $k$-forms with their straight 
  and twisted FE spaces that approximate the corresponding continuous spaces 
  should satisfy the following diagram:
   \begin{equation}\label{equ_spaces_cd}
   \begin{CD}  
     h_h^{(0)} \quad \in \quad \LOh               @>\mathbf{d}>>           \LIh  \quad \ni \quad u_h^{(1)}  \\
              \hspace*{+1.6cm}    @V  \starHLOh/ \starHLIh  \, VV            \hspace*{-1.4cm}               @VV \, \starULOh/ \starULIh  V          \\
       \widetilde h_h^{(1)}   \quad \in \quad \LtwIh      @<\mathbf{d}<<  \LtwOh  \quad \ni \quad \widetilde u_h^{(0)}  \\ 
    \end{CD}
   \end{equation}
   In the following we will provide definitions for these spaces and mappings.

  The preceding diagram is satisfied by the following set of FE spaces:
  \begin{itemize}
   \item $\LOh = \VOh$ and $\LtwOh = \VOh$;
   \item $\LIh = \{\omega_h^{(1)} = \omega_h^{(1)}(x) dx \ : \ \omega_h^{(1)}(x) \in \VIh  \}$; 
   \item $\LtwIh = \{\widetilde \omega_h^{(1)} = \widetilde\omega_h^{(1)}(x) dx \ : \ \widetilde\omega_h^{(1)}(x) \in \VIh  \}$;
  \end{itemize}
  with FE spaces $\VOh$ and $\VIh$ from Sect.~\ref{sec_standard} with piecewise linear basis
  $\{ \fctVOtrial \}_{\indxVOtrial = 1}^{N}$ and piecewise constant basis $\{ \fctVItrial \}_{\indxVItrial = 1}^{N}$, respectively.  
  Note that the superindices in $V_h^i,i=0,1,$ correlate with the polynomial degree 
  whereas the superindices in $\Lambda_h^i,i=0,1,$ denote the degree of the differential forms.
  Choosing other, higher order FE spaces is also possible as long as the relationship between these 
  spaces satisfy diagram~\eqref{equ_spaces_cd}. For this manuscript however, we use the lowest order FE 
  spaces possible. To establish time-dependent $k$-forms, we proceed as above for the mixed FE methods.
  
  Approximating the continuous by discrete $k$-forms taken from the above FE spaces and substituting 
  into \eqref{equ_varsplit_topo} and \eqref{equ_varsplit_metri} give the {\em discrete split variational 
  form} of the split 1D wave equations: 
  find $(u_h^{(1)}(x,t), h_h^{(0)}(x,t),\widetilde h_h^{(1)}(x,t),\widetilde u_h^{(0)}(x,t)) \in  
  (\LIh;0,T) \times (\LOh;0,T) \times (\LtwIh;0,T) \times (\LtwOh;0,T)$
  such that 
  \begin{equation}
    \begin{split}\label{equ_var_form_split}
    \int_L \left( \frac{\partial u_h^{(1)}(x,t)}{\partial t} + g \frac{\partial h_h^{(0)}(x,t)}{\partial x} \right) \hat \chi(x) dx &= 0  \, ,
                 \quad  \forall \hat \chi(x) \in \LIhtest \subset \LI,  \qquad \\
    \int_L \left( \frac{\partial \widetilde h_h^{(1)}(x,t)}{\partial t} + H \frac{\partial \widetilde u_h^{(0)}(x,t)}{\partial x}  \right) \hat \chi(x) dx  & = 0 \, ,
                 \quad  \forall \hat \chi(x) \in \LIhtest \subset \LI, 
   \end{split}
  \end{equation}
  for any test function $\hat \chi(x)$ of the test space $\LIhtest \subseteq \LIh$ subject to the discrete variational metric equations
  \begin{align}
  \widetilde u_h^{(0)} =  \starULih u_h^{(1)} \ \ \text{by} 
         \int_L    \widetilde u_h^{(0)}(x,t)   \hat\tau^i(x)  dx & =  \int_L    u_h^{(1)}(x,t)  \hat\tau^i(x)    dx  \, , 
         \        \forall \hat\tau^i(x) \in \hat\Lambda_h^{i} \subset \Lambda^{i}, i = 0,1, \label{equ_proj_1}\\ 
  \widetilde h_h^{(1)} =  \starHLjh  h_h^{(0)} \ \ \text{by} 
          \int_L    \widetilde h_h^{(1)}(x,t)  \hat\tau^{j}(x)   dx & =  \int_L    h_h^{(0)}(x,t)  \hat\tau^{j}(x)    dx  \, , 
          \        \forall \hat\tau^{j}(x) \in \hat\Lambda_h^{{j}}\subset \Lambda^{{j}}, {j} = 0,1,\label{equ_proj_2}
  \end{align}
  for any test functions $ \hat\tau^{i,j}(x)$ of the test spaces $\hat\Lambda_h^{{i,j}} \subseteq \Lambda_h^{{i,j}}$.  
  Equations~\eqref{equ_proj_1} and \eqref{equ_proj_2} define the discrete Hodge-star operators.
  We use the notation $\starULih$ and $\starHLjh$ for $i,j = 0,1$ to indicate that 
  $\starULOh/\starHLOh$ project on piecewise linear (P1) and $\starULIh/\starHLIh$ on piecewise constant (P0) spaces.

  Equations~\eqref{equ_var_form_split} are  trivial projections of the topological equations into 
  the test space $\LIhtest$ of piecewise constant test functions.   
  The discrete exterior derivative $\bf d$ is a surjective map from piecewise linear to piecewise constant functions and 
  can be represented by the metric-free coincidence matrix $\stiffmxen$ (see below). 
  The discrete topological equations are therefore exact up to the 
  errors that occur by the projections into the piecewise constant spaces.

  Equations~\eqref{equ_proj_1} and \eqref{equ_proj_2} are nontrivial Galerkin projections between different 
  spaces and provide approximations of the Hodge-star operator. For solving the system, 
  the discrete Hodge-star operator -- unlike $\bf d$ -- has to be inverted, 
  which requires special treatment in case $\starHLOh$ or $\starHLIh$ has a non-trivial kernel 
  (cf. Sect.~\eqref{sec_split_CN}). 
  Consequently, the discrete metric closure equations introduce additional errors into the systems,
  similarly to those occurring by partial integrations when formulating the wave equations weakly. 
%

  Note that a projection of the topological equations into the twisted test space $\LtwIhtest$ 
  would result in the same variational form because all minus signs carried 
  by the twisted forms would compensate. However, a projection into a combination
  of straight and twisted spaces is not allowed as it would change the sign of the 
  original wave equations. This statement holds also for the projection of the 
  metric equations. 
  
  
  In the next section we will derive suitable matrix-vector representations of the split schemes. 
  Following our ideas of treating topological and metric equations separately, 
  let us first consider the discrete topological equations in Sect.~\ref{sec_disc_top}
  and then the discrete metric equations in Sect.~\ref{sec_disc_metr} 
  that close the system of equations.

  \subsubsection{Discrete topological equations}
  \label{sec_disc_top}
  
  We expand the four variables by means of the above introduced FE spaces: 
  \begin{itemize}
   \item momentum pair (straight forms): \  
       ${u}_h^{(1)}(x,t) = \sum_{\indxVItrial = 1}^{N} \uVItrial(t) \fctVItrial$, \ 
       ${h}_h^{(0)}(x,t) = \sum_{\indxVOtrial = 1}^{N} \hVOtrial(t) \fctVOtrial$;
   \item continuity pair (twisted forms): \quad
       ${\widetilde{h}_h^{(1)}}(x,t) = \sum_{\indxVItrial = 1}^{N} \htwVItrial(t) \fctVItrial$, \ 
       ${\widetilde{u}_h^{(0)}}(x,t) = \sum_{\indxVOtrial = 1}^{N} \utwVOtrial(t) \fctVOtrial$; 
  \end{itemize}
  and substitute them into \eqref{equ_var_form_split}. Varying the test functions 
  $\hat \chi(x)$ over the basis $\{ \fctVItest  \}_{\indxVItest = 1}^{N}$ of $\LIhtest$,
  we obtain 
  \begin{equation}\label{equ_splitFE_top}
  \begin{split}
          \sum_{\indxVItrial = 1}^{N} \partial_t \uVItrial(t)  \int_L \fctVItrial \fctVItest dx      
      + g \sum_{\indxVOtrial = 1}^{N} \hVOtrial(t)  \int_L \frac{d\fctVOtrial}{dx} \fctVItest dx 
   & =     0  \, ,  \ \ \text{\ for\ } \indxVItest = 1,\dots,N  ,\\
          \sum_{\indxVItrial = 1}^{N} \partial_t \htwVItrial(t) \int_L \fctVItrial \fctVItest  dx 
      + H \sum_{\indxVOtrial = 1}^{N} \utwVOtrial(t) \int_L \frac{d\fctVOtrial}{dx} \fctVItest dx
   & =     0   \, ,    \ \ \text{\ for\ } \indxVItest = 1,\dots,N  .
  \end{split}
  \end{equation}
 
  The preceding equations are similar in form to the second equation of \eqref{equ_P0P1}.
  They can therefore be written in matrix-vector form analogously to the continuity 
  equation in~\eqref{equ_matrix_stbl} by using the matrices $\massmxee$ 
  with metric-dependent coefficients~\eqref{equ_mssmxee_coeff} and $\stiffmxen$ with 
  metric-free coefficients~\eqref{equ_stiffxne_coeff}. 
  Defining the vector $\metrivcte = (\Delta x_1,\dots \Delta x_m, \dots \Delta x_{N})$
  containing metric information of the mesh, we note that $\massmxee = {\bf Id}\cdot\metrivcte$.
  The metric coefficients combined with the coefficients for velocity and height constitute  
  the discrete 1-forms ${\bf u}_e^{(1)}$ and ${\bf \widetilde  h}_e^{(1)}$ that approximate
  the respective continuous 1-forms. 
  The discrete topological (metric-free) equations then read 
  \begin{align}\label{equ_matrix_splitFE}
     \frac{\partial {\bf u}_e^{(1)}}{\partial t}             + g \stiffmxen             {\bf h}_n^{(0)} = 0 , \quad 
     \frac{\partial \widetilde {\bf h}_e^{(1)}}{\partial t}  + H \stiffmxen  \widetilde {\bf u}_n^{(0)} = 0 ,
  \end{align}
  using the following definitions:
  \begin{itemize}
   \item ${\bf u}_e^{(1)} = (u_1(t) \Delta x_1 ,\dots,u_m(t) \Delta x_m,\dots  u_{N}(t) \Delta x_{N})$ approximates the 1-form $u^{(1)}\in (\LI;0,T)$;
   \item ${\bf \widetilde  h}_e^{(1)}= (\widetilde h_1(t) \Delta x_1,\dots \widetilde h_m(t) \Delta x_m,\dots   \widetilde h_{N}(t) \Delta x_{N})$ 
         approximates the 1-form $\widetilde h^{(1)} \in (\LtwI;0,T)$;
   \item ${\bf \widetilde  u}_n^{(0)}= (\widetilde u_1(t) ,\dots \widetilde u_l(t), \dots \widetilde u_{N}(t) )$ approximates the 0-form $\widetilde u^{(0)} \in (\LtwO;0,T)$;
   \item ${\bf h}_n^{(0)}= (h_1(t), \dots h_l(t),\dots h_{N}(t) )$ approximates the 0-form $h^{(0)}  \in (\LO;0,T)$.
  \end{itemize}
  
  Because the stiffness matrix $\stiffmxen$ is a metric-free approximation of the exterior derivative $\bf d$ 
  applied in both straight and twisted sequences of FE spaces, the discrete topological momentum and continuity 
  equations~\eqref{equ_matrix_splitFE} provide metric-free approximations of the corresponding continuous
  topological equations in \eqref{equ_split_wave}.

  \subsubsection{Discrete metric closure equations}
  \label{sec_disc_metr}
  
   \newcommand{\GPIu}{\color{black}{\rm GP1}_{u}\color{black}}  
   \newcommand{\GPOu}{\color{black}{\rm GP0}_{u}\color{black}}  
   \newcommand{\GPIh}{\color{black}{\rm GP1}_{h}\color{black}}  
   \newcommand{\GPOh}{\color{black}{\rm GP0}_{h} \color{black}}  
  
  \newcommand{\GPI}[1]{\color{black} {\rm GP1}_{{#1}} \color{black}}
  \newcommand{\GPO}[1]{\color{black} {\rm GP0}_{{#1}} \color{black}}

  For the given choice of FE spaces, there exist {\bf four} realizations of the 
  discrete metric equations \eqref{equ_proj_1} and \eqref{equ_proj_2}. We will classify 
  these realizations into three groups depending on how accurately they approximate 
  the continuous metric closure equations: 
  (i) {\em high accuracy closure} using the pair $(\starULOh,\starHLOh)$,
  (ii) {\em low accuracy closure} using  the pair $(\starULIh,\starHLIh)$,
  and (iii) {\em medium accuracy closure} using either $(\starULOh,\starHLIh)$ 
  or $(\starULIh,\starHLOh)$. 
  As the discrete Hodge-star operators are realized by nontrivial Galerkin 
  projections (GP) onto either the piecewise constant ($\GPO{}$) or piecewise linear ($\GPI{}$) space,
  we denote these schemes also by: (i) $\GPIu$--$\,\GPIh$, (ii) $\GPOu$--$\,\GPOh$, and (iii)
  $\GPIu$--$\,\GPOh$ or $\GPOu$--$\,\GPIh$, respectively, in analogy to the conventional 
  notation for mixed P1--P1 and P1--P0 schemes.

  \paragraph{High accuracy closure ($\GPIu$--$\,\GPIh$).} 
  The most accurate approximation of the metric closure equations~\eqref{equ_varsplit_metri} 
  is achieved by using the discrete Hodge-stars $\starULOh$ and $\starHLOh$ 
  in \eqref{equ_proj_1} and \eqref{equ_proj_2}, respectively. Hence,
  \begin{equation}\label{equ_splitFE_m}
    \begin{split}
     \sum_{\indxVOtrial = 1}^{N} \left\{ \begin{array}{c}    \utwVOtrial(t)   \\  \hVOtrial(t)    \end{array} \right\} 
       \int_L \fctVOtrial \fctVOtest  dx 
   =  \sum_{\indxVItrial = 1}^{N} \left\{ \begin{array}{c}    \uVItrial(t)   \\ \htwVItrial(t)  \end{array} \right\}  
        \int_L \fctVItrial \fctVOtest  dx \, , \quad \text{\ for\ } \indxVOtest = 1,\dots,N  .
    \end{split}
  \end{equation}
    
  To write equations~\eqref{equ_splitFE_m} in matrix-vector form, we note that the terms 
  on the left hand side are simply the metric-dependent coefficients~\eqref{equ_mssmxnn_coeff}
  of $\massmxnn$,  while those on the right form the $(N\times N)$ mass-matrix
  $\massmxne$ with metric-dependent coefficients 
  \begin{equation}\label{equ_mssmxne_coeff}
   M_{\indxVItrial \indxVOtest} = \int_L \fctVItrial \fctVOtest dx.
  \end{equation}
  This matrix can be written equivalently as $\massmxne = \prjcmxne (\metrivcte)^T $, 
  where $\prjcmxne$ is a metric-free $(N\times N)$ matrix representing a projection 
  operator $\LIh \rightarrow \LtwOh $ and $\LtwIh \rightarrow \LOh $. The explicit 
  representation of this projection (resp. mass-matrix) corresponds to determining node values from 
  averaging the neighboring cell values (resp. area-weighted cell values) (cf. Appendix~\ref{append_int_values}).
  We obtain the metric-dependent equations
  \begin{equation}\label{equ_mtx_high}
        \massmxnn {\bf \widetilde  u}_n^{(0)}   = \prjcmxne        {\bf u }_e^{(1)}  \, ,  \qquad\qquad
        \massmxnn {\bf h }_n^{(0)}               = \prjcmxne    {\bf \widetilde  h}_e^{(1)} \, ,
  \end{equation}
  in which we have combined the elements of $\metrivcte$ denoting the area weights of each element 
  with the coefficients for velocity and height to the discrete 1-forms ${\bf u }_e^{(1)}$ and 
  ${\bf \widetilde  h}_e^{(1)}$, respectively.

  \paragraph{Low accuracy closure ($\GPOu$--$\,\GPOh$).} 
  The most inaccurate approximation of the metric closure equations~\eqref{equ_varsplit_metri}
  results from using $\starULIh$ and $\starHLIh$ in \eqref{equ_proj_1} and \eqref{equ_proj_2}, respectively:
  \begin{equation}\label{equ_splitFE_m3}
  \begin{split}
      \sum_{\indxVOtrial = 1}^{N} \left\{ \begin{array}{c}    \utwVOtrial(t)   \\  \hVOtrial(t)    \end{array}  \right\} 
        \int_L \fctVOtrial  \fctVItest dx 
    = \sum_{\indxVItrial = 1}^{N} \left\{ \begin{array}{c}    \uVItrial(t)   \\ \htwVItrial(t)  \end{array}  \right\}  
        \int_L \fctVItrial  \fctVItest  dx \, , \  \text{\ for\ } \indxVItest = 1,\dots,N  .
   \end{split}
   \end{equation}
   In terms of matrix-vector formulation, these discrete metric-dependent equations read
   \begin{equation}\label{equ_mtx_low}
       \massmxen   {\bf \widetilde  u}_n^{(0)}   =   {{\bf Id}^{ee}}  {\bf u }_e^{(1)}  \, ,  \qquad\qquad 
        \massmxen  {\bf h }_n^{(0)}              =   {{\bf Id}^{ee}}  {\bf \widetilde  h_e}^{(1)} \, ,
   \end{equation}
   in which the metric-dependent $(N\times N)$ mass-matrix $\massmxen$ can be easily computed by $\massmxen = (\massmxne)^T$.
   For the terms on the right hand side, we exploit the fact that $\massmxee = {\bf Id}^{ee}\,(\metrivcte)^T$, 
   which again allows us to combine the area weights of $\metrivcte$ with the 
   corresponding coefficients for velocity or height to the 1-forms 
   ${\bf u }_e^{(1)}$ resp. ${\bf \widetilde  h_e}^{(1)}$ for each element.

  \paragraph{Medium accuracy closure ($\GPIu$--$\,\GPOh/\GPOu$--$\,\GPIh$).} 
  An approximation to the metric closure equations~\eqref{equ_varsplit_metri} at an intermediate accuracy 
  results from applying the approximations $\starULOh$ and $\starHLIh$ 
  in \eqref{equ_proj_1} and \eqref{equ_proj_2}, respectively:
  \begin{equation}\label{equ_splitFE_m6}
    \begin{split}
        \sum_{\indxVOtrial = 1}^{N} \utwVOtrial(t)  \int_L \fctVOtrial  \fctVOtest  dx 
    & = \sum_{\indxVItrial = 1}^{N} \uVItrial(t)    \int_L \fctVItrial  \fctVOtest  dx  \, , \quad \text{\ for\ } \indxVOtest = 1,\dots,N,  \\ 
        \sum_{\indxVOtrial = 1}^{N} \hVOtrial(t)    \int_L \fctVOtrial  \fctVItest  dx 
    & = \sum_{\indxVItrial = 1}^{N} \htwVItrial(t)  \int_L \fctVItrial  \fctVItest  dx \, , \ \  \text{\ for\ } \indxVItest = 1,\dots,N . 
   \end{split}
  \end{equation}
  Using the matrices defined above, these discrete metric-dependent equations 
  yield the matrix-vector form: 
  \begin{equation}\label{equ_mtx_med1}
      \massmxnn  {\bf \widetilde  u}_n^{(0)}   =   \prjcmxne {\bf u }_e^{(1)}  \, ,  \qquad\qquad 
      \massmxen  {\bf h }_n^{(0)}              =   {{\bf Id}^{ee}}   {\bf \widetilde  h}_e^{(1)} \, .
  \end{equation}
  
  A similarly accurate approximation of the metric closure equations~\eqref{equ_varsplit_metri} follows from 
  applying the approximations $\starULIh$ and $\starHLOh$ in 
  \eqref{equ_proj_1} and \eqref{equ_proj_2}, respectively:
  \begin{equation}\label{equ_splitFE_m7}
    \begin{split}
       \sum_{\indxVOtrial = 1}^{N} \utwVOtrial(t) \int_L \fctVOtrial \fctVItest  dx 
  & =  \sum_{\indxVItrial = 1}^{N} \uVItrial(t)   \int_L \fctVItrial \fctVItest  dx \, , \ \ \text{\ for\ } \indxVItest = 1,\dots,N,	  \\
       \sum_{\indxVOtrial = 1}^{N} \hVOtrial(t)   \int_L \fctVOtrial \fctVOtest  dx 
  & =  \sum_{\indxVItrial = 1}^{N} \htwVItrial(t) \int_L \fctVItrial \fctVOtest  dx \, , \, \quad \text{\ for\ } \indxVOtest = 1,\dots,N .
   \end{split}
  \end{equation}
  In matrix-vector form, they read 
  \begin{equation}\label{equ_mtx_med2}
       \massmxen  {\bf \widetilde  u}_n^{(0)}   = {{\bf Id}^{ee}}  {\bf u }_e^{(1)}  \, ,  \qquad\qquad 
       \massmxnn {\bf h }_n^{(0)}               =       \prjcmxne {\bf \widetilde  h}_e^{(1)} \, .
  \end{equation}

  Equations~\eqref{equ_mtx_high}, \eqref{equ_mtx_low}, and \eqref{equ_mtx_med1}/\eqref{equ_mtx_med2}
  provide discrete metric-dependent approximations to the metric equations in \eqref{equ_split_wave}.
  More precisely, we can determine the discrete Hodge-star operators with respect to the two possible 
  accuracy levels by
  \begin{equation}
   \begin{split}\label{equ_def_disc_hodge}
        \widetilde u^{(0)}  = \widetilde\star u^{(1)} \quad & \approx \quad  {\bf \widetilde  u}_n^{(0)} =  \underbrace{ (\massmxnn)^{-1} \prjcmxne             }_{=:\starULOh}  {\bf u }_e^{(1)} 
        \quad \text{or} \quad                                                {\bf \widetilde  u}_n^{(0)} =  \underbrace{ (\massmxen)^{-1} {{\bf Id}^{ee}}  }_{=:\starULIh}  {\bf u }_e^{(1)} ,\\
        \widetilde h^{(1)}  = \widetilde\star h^{(0)} \quad & \approx \quad  {\bf \widetilde  h}_e^{(1)} =  \underbrace{ (\prjcmxne)^{-1} \massmxnn             }_{=:\starHLOh}  {\bf h }_n^{(0)}
        \quad \text{or} \quad                                                {\bf \widetilde  h}_e^{(1)} =  \underbrace{ {({\bf Id}^{ee})^{-1} } \massmxen         }_{=:\starHLIh}  {\bf h }_n^{(0)},
   \end{split}
   \end{equation}
   which allows us to write diagram~\eqref{equ_spaces_cd} in terms of matrix-vector formulation as
   \begin{equation}\label{equ_spaces_cd_matrices}
   \begin{CD}  
     {\bf h }_n^{(0)}    \quad \in \quad \LOh               @>\mathbf{\stiffmxen}>>           \LIh  \quad \ni \quad  {\bf u }_e^{(1)}  \\
             \hspace*{+1.6cm}    @V { (\prjcmxne)^{-1} \massmxnn / ({\bf Id}^{ee})^{-1}  \massmxen} VV     
             \hspace*{-1.4cm}   @VV { (\massmxnn)^{-1} \prjcmxne / (\massmxen)^{-1} {{\bf Id}^{ee}} }   V          \\
       {\bf \widetilde  h}_e^{(1)}   \quad \in \quad \LtwIh      @<\mathbf{\stiffmxen}<<  \LtwOh  \quad \ni \quad {\bf \widetilde  u}_n^{(0)}  \\ 
    \end{CD}
   \end{equation}  
  using the matrix representations of $\bf d = \stiffmxen$. 
  
  Together with \eqref{equ_matrix_splitFE} each of these realizations of the metric equations forms 
  a closed system of semi-discrete matrix-vector equations that approximate the split wave equations. 
  In the next section, we will discuss the time-discretization of these sets of equations 
  and introduce a solution method for the resulting systems of algebraic equations.

  \subsection{Time discretization and solving algorithm.}
  \label{sec_split_CN}
  
  As for the mixed FE schemes, we use the symmetric, implicit Crank-Nicolson time discretization.
  For a time $t \in [0,T]$ with time step size $\Delta t$, each fully discretized split FE scheme reads 
  \begin{equation}\label{equ_matrix_eq_CN_split_full}
    \begin{split}
	 {\bf u}^{(1)}_{t+1} & =              {\bf u}^{(1)}_{t} +  \frac{1}{2}\Delta t g {\bf D}^{en} {\bf h}^{(0)}_{t} 
	 + \frac{1}{2}\Delta t g {\bf D}^{en} {\bf h}^{(0)}_{t+1} \, , \\ 
	\widetilde {\bf h}^{(1)}_{t+1} & =    \widetilde {\bf h}^{(1)}_{t} +  \frac{1}{2}\Delta t H {\bf D}^{en} \widetilde {\bf u}^{(0)}_{t} 
				    +  \frac{1}{2}\Delta t H {\bf D}^{en} \widetilde {\bf u}^{(0)}_{t+1} \, ,
    \end{split}
  \end{equation}
  subject to one of the projection pairs \eqref{equ_mtx_high}, \eqref{equ_mtx_med1}/\eqref{equ_mtx_med2},
  or \eqref{equ_mtx_low}. These metric equations do not depend on time, hence they need no special time discretization.
  For better readability, we skipped the subindices $e,n$.

  \paragraph{Solving algorithm.}
  The algorithm to solve the implicit equations is given by the {\it fixed point iteration}:
  \begin{enumerate}
   \item start loop over $k$ with initial guess ($k=0$) from values of previous timestep $t$: 
         ${\bf u}^{(1)}_{*,k=0} = {\bf u}^{(1)}_{t}$ and  ${\bf \widetilde h}^{(1)}_{*,k=0} = {\bf \widetilde h}^{(1)}_{t}$;
   \item project ${\bf u}^{(1)}_{*,k}$ onto ${\bf \widetilde u}^{(0)}_{*,k}$ by $\starULih$ for the appropriate
         test space $\hat\Lambda_h^{i}, i=0,1$;
   \item calculate height ${\bf \widetilde h}^{(1)}_{*,k+1}$:
     \begin{equation}\notag
        \widetilde {\bf h}^{(1)}_{*,k+1}  =  D_{t}   +  \frac{1}{2}\Delta t H {\bf D}^{en} \widetilde {\bf u}^{(0)}_{*,k} 
                                         \quad   \text{with} \quad  D_{t} := \widetilde {\bf h}^{(1)}_{t} +  \frac{1}{2}\Delta t H {\bf D}^{en} \widetilde {\bf u}^{(0)}_{t} \, ; 
     \end{equation}
    \item project ${\bf \widetilde h}^{(1)}_{*,k+1}$ onto ${\bf h}^{(0)}_{*,k+1}$ by $(\starHLjh)^{-1}$ for the appropriate
          test space $\hat\Lambda_h^{{j}}, {j}=0,1$;
    \item calculate velocity ${\bf u}^{(1)}_{*,k+1}$:
     \begin{equation}\notag
      {\bf u}^{(1)}_{*,k+1}  =  F_{t}   +  \frac{1}{2}\Delta t g {\bf D}^{en}  {\bf h}^{(0)}_{*,k+1} 
                                       \quad   \text{with} \quad  F_{t} := {\bf u}^{(1)}_{t} +  \frac{1}{2}\Delta t H {\bf D}^{en}  {\bf h}^{(0)}_{t} \, ; 
     \end{equation}
   \item set $k+1 = k$ and stop loop over $k$ if $||{\bf u}^{(1)}_{*,k+1}\vspace{-0.1em}   -  {\bf u}^{(1)}_{*,k} || 
            + || \widetilde {\bf h}^{(1)}_{*,k+1} -\widetilde {\bf h}^{(1)}_{*,k} || < \epsilon$ for a small positive $\epsilon$.
  \end{enumerate}
  In case of convergence, i.e. ${\bf u}^{(1)}_{*,k+1}\rightarrow{\bf u}^{(1)}_{t+1} $  
  and $\widetilde {\bf h}^{(1)}_{*,k+1}\rightarrow \widetilde {\bf h}^{(1)}_{t+1} $, 
  this algorithm indeed solves equations~\eqref{equ_matrix_eq_CN_split_full}. 
  Note that ${\bf u}^{(1)}_{*,k+1}$ could also be computed before ${\bf \widetilde h}^{(1)}_{*,k+1}$.

  As described in Sect.~\eqref{sec_standard00}, we have to restrict the time step size to 
  guarantee that the fastest waves are well resolved and that the iterative solver converges. 
  The corresponding CFL numbers, given by \eqref{cflnumber}, are determined numerically for 
  all meshes and test cases studied. We use
  \begin{itemize}
  \item $\mu \leq 1.15$ for the $\GPIu$--$\,\GPIh$ scheme (similar to the P1--P1 scheme);
  \item $\mu \leq 0.57$ for the $\GPIu$--$\,\GPOh/\GPOu$--$\,\GPIh$ schemes (similar to the P1--P0 scheme);
  \item $\mu \leq \frac{1.15}{N/2}$ for the $\GPOu$--$\,\GPOh$ scheme.
  \end{itemize}
  Note that in case of the $\GPOu$--$\,\GPOh$ scheme the CFL condition scales with the grid resolution (see explanation below).

  \paragraph{Treatment of possible non-trivial kernels of the discrete Hodge-stars.}

  Within each time step of the fixed point iteration algorithm, we compute
  an increment to the time update of ${\bf u }_e^{(1)}$ and map it to ${\bf \widetilde  u}_n^{(0)}$ 
  in order to compute an increment to the time update of ${\bf \widetilde  h}_e^{(1)}$;
  then we map the latter to ${\bf h }_n^{(0)}$ to determine the next more accurate increment to
  ${\bf u }_e^{(1)}$, until the required accuracy is achieved.
  This procedure corresponds to the mappings in diagram~\eqref{equ_spaces_cd_matrices}, where 
  the directions of mappings are indicated correctly by the arrows, except for the mappings between 
  height fields. As indicated in step 4) of the algorithm above, the inverse $(\starHLjh )^{-1}, {j}=0,1$ 
  is required to map ${\bf \widetilde  h}_e^{(1)}$ to ${\bf h }_n^{(0)}$.
  Moreover, since we want to study all possible combinations of metric closure equations, 
  we have to invert all matrices present in \eqref{equ_def_disc_hodge}. 
    
  In general it is not clear if the inverse matrices exist nor if they are computable, and their general
  treatment is an object of current research but not in the scope of this manuscript. Here, 
  we only consider the periodic 1D mesh of Fig.~\ref{fig_mesh} with $N$ elements. 
  In this case $\massmxnn$ has full rank and is always invertible. In case $N$ is odd,
  also $\prjcmxne$ and $\massmxen$ have full rank and are invertible, which allows us to 
  calculate $\starULIh,\starULOh,\starHLIh$, and $\starHLOh $ and their inverse. 
  However, in case $N$ is even, $\prjcmxne$ and $\massmxen$ have 
  one-dimensional kernels (null-vectors). 
 
  In order to fix this deficiency for even $N$, we consider the map 
  $(\starHLIh )^{-1}:  {\bf \widetilde  h}_e^{(1)} \mapsto  {\bf h}_n^{(0)} $
  with null-vector ${\bf K}_{0} $ of length $N$. Then we solve the augmented 
  problem: given ${\bf \widetilde  h}_e^{(1)} \in \LtwIh$, 
  find ${\bf h }_n^{(0)}\in \LOh$ such that
  \begin{equation} \label{equ_augmentedsystem}
  \left( \begin{matrix} 
             \massmxen           & { \bf K}_{0}   \\
             ( {\bf K}_{0}  )^T   &   0 
          \end{matrix} \right) 
   \left( \begin{matrix} 
             {\bf h }_n^{(0)}\\
             0 
          \end{matrix} \right) = 
   \left( \begin{matrix} 
             {\bf \widetilde  h}_e^{(1)} \\
             0 
          \end{matrix} \right) . 
  \end{equation}
  Analogously, we treat all other cases in which we have to use the inverse 
  of a matrix that has a non-trivial kernel.

  For the projection steps 2) and 4) of the algorithm, we can choose different spaces. In the 
  following paragraph and in Sect.~\ref{sec_num_analysis}, we will study how this choice 
  influences the performance of the resulting numerical scheme, in particular regarding
  it's accuracy, stability, and convergence properties.

 \subsection{Discussion about spatial stability}
 \label{sec_discuss_stability}
  
  Analogously to the mixed FE cases, we study the stability properties of the split FE
  schemes by investigating the discrete dispersion relations. In particular, we 
  investigate the impact of the choice of projection accuracy of the metric closure equations 
  \eqref{equ_proj_1} and \eqref{equ_proj_2} on the dispersion relations. 
  We present the results of an analytic derivation of the discrete dispersion relations for each 
  scheme (for details see Appendix~\ref{app_disc_disp_rel}).
    
\begin{figure}[t]\centering 
     \includegraphics[scale=1.2]{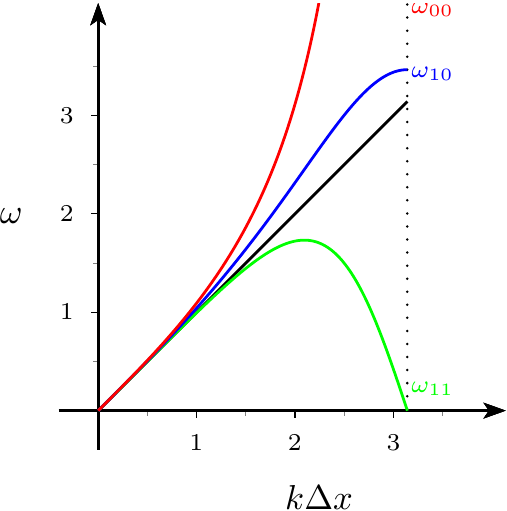}
    \caption{Dispersion relations: analytic (black) for $c=\sqrt{gH}$, ${\omega_{11}}$ in 
    green for the Hodge-star pair $(\starULOh,\starHLOh)$,
    ${\omega_{10}}$ in blue for $(\starULOh,\starHLIh)$ or $(\starULIh,\starHLOh$),
    and ${\omega_{00}}$ in red for $(\starULIh,\starHLIh$). }     
   \label{fig_disp_relation}
  \end{figure}

  For the split $\GPIu$--$\,\GPIh$ scheme, consisting of the topological momentum and 
  continuity equations~\eqref{equ_matrix_splitFE} and the metric closure equations
  \eqref{equ_mtx_high}, the angular frequency $\omega_{11} := \omega(k)$ 
  (the subscripts will indicate the order of the Hodge-star operators, in this case $\starULOh$ and $\starHLOh$)
  satisfies the discrete dispersion relation
  \begin{equation}\label{equ_disprel_split00}
	c_d = \frac{\omega_{11}   }{k} = \pm \sqrt{gH} \frac{\sin (k\frac{\Delta x}{2})}{k \frac{\Delta x}{2}} \frac{ 3\cos (k \frac{\Delta x}{2}) }{ (2 + \cos (k\Delta x))  }.
  \end{equation}
  One realizes that $c_d \rightarrow c = \sqrt{gH}$ for $k \rightarrow 0$. This can 
  be seen by expanding $\sin (k\frac{\Delta x}{2})$ in a Taylor series around zero. 
  Because of the double angle formula (i.e. $\sin (k \Delta x) = 2\sin(k\frac{\Delta x}{2})\cos(k\frac{\Delta x}{2}  )$),
  this dispersion relation equals the one in \eqref{equ_disprel_P1P1} for the P1--P1 scheme; 
  hence, it has a zero solution at $k= \frac{\pi}{\Delta x}$ (cf. Fig.~\ref{fig_disp_relation}) and it is therefore unstable.
  

  For the split $\GPIu$--$\,\GPOh$ and $\GPOu$--$\,\GPIh$ schemes which consist of the topological momentum and 
  continuity equations~\eqref{equ_matrix_splitFE} and metric closure equations
  \eqref{equ_mtx_med1} or \eqref{equ_mtx_med2}, 
  the discrete dispersion relation of the angular frequency $\omega_{10}  :=\omega(k)$
  for both combinations of metric closure equations reads
  \begin{equation}\label{equ_disprel_split01}
           c_d = \frac{{\omega_{10}}}{k} = \pm \sqrt{gH} \frac{\sin (k\frac{\Delta x}{2})}{k \frac{\Delta x}{2}} \left[ \frac{ 3  }{ (2 + \cos (k\Delta x))  } \right]^{\frac{1}{2}},
  \end{equation}
  with $c_d \rightarrow c = \sqrt{gH}$ for $k \rightarrow 0$.
  The preceding equation provides exactly the same dispersion relation as in \eqref{equ_disprel_P0P1} 
  for the stable P1--P0 FE scheme. Hence, both intermediate accuracy split schemes are stable as they expose no 
  spurious mode at $k= \frac{\pi}{\Delta x}$.

  For the split $\GPOu$--$\,\GPOh$ scheme, which consists of the topological momentum and 
  continuity equations~\eqref{equ_matrix_splitFE} and metric closure equations
  \eqref{equ_mtx_low}, the angular frequency ${\omega_{00}} := \omega(k)$ satisfies the discrete dispersion relation 
  \begin{equation}\label{equ_disprel_split11}
           c_d = \frac{{\omega_{00}} }{k} = \pm \sqrt{gH} \frac{\tan (k\frac{\Delta x}{2})}{k \frac{\Delta x}{2}} ,
  \end{equation}
  with $c_d \rightarrow c = \sqrt{gH}$ for $k \rightarrow 0$. 
  This dispersion relation has no second root. However, with increasing $k$ the wave speeds exceed
  the analytical solution significantly and grow infinitely for the smallest wave $k= \frac{\pi}{\Delta x}$. 
  We identify these fast traveling small scale waves as spurious modes, since
  they cause small scale noise on the entire mesh after some short simulation 
  time. 
  
  In fact, the occurrence of these spurious modes are the reason why the CFL number 
  $\mu \leq \frac{1.15}{N/2}$ of the iterative solver (cf. Sect.~\ref{sec_split_CN}) 
  is very small and even decreases proportionally to the element size of the mesh. 
  With $\Delta x \rightarrow 0$ and considering the short wave $k\rightarrow \frac{\pi}{\Delta x}$ 
  we have that $\frac{k\Delta x}{2} \rightarrow \frac{\pi}{2}$ and hence the discrete wave 
  speed $c_d \rightarrow \infty$ in \eqref{equ_disprel_split11}. So, with increasing mesh 
  resolution the wave speed increases and the CFL condition becomes more restrictive.
  The suggested $\mu \leq \frac{1.15}{N/2}$ for the $\GPOu$--$\,\GPOh$ scheme provides 
  an upper bound for the CFL number and assures convergence of the iterative solver
  for all meshes used and all test cases studied.

  Figure~\ref{fig_disp_relation} shows the dispersion relations plotted for the three different 
  choices of Hodge-star projections:  
  in green the high accuracy ($\GPIu$--$\,\GPIh$), 
  in blue the intermediate accuracy ($\GPIu$--$\,\GPOh$/$\GPOu$--$\,\GPIh$), 
  and in red the low accuracy ($\GPOu$--$\,\GPOh$) scheme, while 
  the analytic dispersion relation is shown in black.

 \section{Numerical analysis}
  \label{sec_num_analysis}

 The objective of this section is threefold. 
 First, we study if the mixed and split FE discrete schemes are structure-preserving in the sense that 
 they preserve the first principles of mass and momentum conservation.
 Second, we study the convergence behavior with respect to analytical 
 solutions of the 1D wave equations. Third, we investigate if the schemes 
 represent the corresponding discrete dispersion relations correctly. 
 
 As discussed in sections \ref{sec_standard} and \ref{sec_discuss_stability}, 
 the discrete dispersion relations differ from the analytical ones, in particular 
 for small wave numbers. In order to avoid these error sources when investigating
 the convergence properties, we perform simulations with a single sine wave
 (corresponding to wave number $k=1$) in test case (TC) 1. 
 Limiting the convergence tests to long waves will keep the errors introduced by 
 the discrete dispersion relation small. This will allow us to test convergence of 
 the schemes for integration times up to 5 cycles 
 (which means the wave propagates 5 times through the periodic domain),
 despite of the influence of the wave dispersion on the shape of the wave package. 
 
 In a second test case (TC 2), we use a Gaussian distribution function with intermediate width, 
 which consists of a superposition of all possible wave numbers. Here, the errors in the discrete 
 dispersion relations will lead to a separation of the Gaussian wave package 
 into its single wave components. As this effect grows in time, we will study 
 the convergence properties of the schemes only for short integration times (less than 1 cycle).
 
 Finally in TC 3, we use the Gaussian distribution function with very small width such that the 
 high frequency waves have large magnitude. This allows us to study the dispersion behavior 
 of each scheme numerically and to compare it with the analytic results of Sections~\ref{sec_standard}
 and \ref{sec_discuss_stability}.

 \paragraph{Test cases and initialization.}
 
 For all numerical simulations, we use a uniform mesh with periodic boundaries, as 
 shown in Fig.~\ref{fig_mesh}, for a domain of length $L = 1000\,$m. Simulation times 
 are measured in full periods (cycles) with time $T= \frac{L}{c}$ for 1 cycle, 
 with respect to the wave speed $c = \sqrt{gH}$. 
 To keep the temporal error for all mesh resolutions in the same order,
 we use in all simulations the same fixed time step size of $\Delta t = 6.3102 \cdot 10^{-04}\,$s for the 
 high and medium accuracy, and $\Delta t =  \frac{1}{200} 6.3102\cdot 10^{-04}\,$s 
 for the low accuracy schemes. These values follow from the CFL number presented in Sect.~\ref{sec_split_CN}
 and from the grid spacing of the highest resolved mesh applied.

 We use the following analytical solutions of the 1D wave equations with 
 $x \in [0,L]$ and $t\in \mathbbm R$ as initializations at time $t=0$ 
 and to determine the convergence rates: 
 \medskip
 
 \noindent Analytical solution for TC 1; {\em single sine wave}:
 \begin{equation}\label{equ_h_u_sinus_ana}
 \begin{split}
    h(x,t)  &  = H  +  \frac{\Delta H}{2} \sin \left( \frac{2\pi}{L}\big(x - c  t \big) \right) 
                  +  \frac{\Delta H}{2} \sin \left( \frac{2\pi}{L}\big(x + c  t \big) \right) , \\
    u(x,t)  & =   \frac{c \Delta H}{2H} \sin \left( \frac{2\pi}{L}\big(x - c  t \big) \right) 
              -  \frac{c \Delta H}{2H} \sin \left( \frac{2\pi}{L}\big(x + c  t \big) \right) ,             
 \end{split}
 \end{equation}
 with parameters: $c = \sqrt{gH}$, $H = 1000\,$m, $\Delta H = 75\,$, 
 and gravitational acceleration $g = 9.81\, \rm{m}/\rm{s}^2$.

  \medskip
  \noindent Analytical solution for TC 2/TC 3; {\em Gaussian distribution with intermediate/small width}:
  \begin{equation}\label{equ_h_u_gauss_ana}
  \begin{split}
   h(x,t) & = H     +    \frac{\Delta H}{2}     e ^{ -\Big( \frac{\Delta w}{2\pi}\sin\big(\frac{\pi}{L} ( x - ct - x_c )\big)\Big)^2 }  
                    +    \frac{\Delta H}{2}     e ^{ -\Big( \frac{\Delta w}{2\pi}\sin\big(\frac{\pi}{L} ( x + ct - x_c )\big)\Big)^2}  , \\
   u(x,t) & =       +    \frac{c \Delta H}{2H}  e ^{ -\Big( \frac{\Delta w}{2\pi}\sin\big(\frac{\pi}{L} ( x - ct - x_c )\big)\Big)^2}  
                    -    \frac{c \Delta H}{2H}  e ^{ -\Big( \frac{\Delta w}{2\pi}\sin\big(\frac{\pi}{L} ( x + ct - x_c )\big)\Big)^2} ,
  \end{split}               
   \end{equation}
   with parameters: $\Delta w = 40$, $x_c = \frac{1}{2}L$. For TC 3, we use the preceding analytical
   solution with very small width by setting $\Delta w = 1000$. 
   
   A direct calculation shows that both analytical expressions are indeed solutions of
   the 1D wave equations~\eqref{equ_standard_wave} or \eqref{equ_split_wave}. To initialize the discrete schemes,
   we project the functions $h(x,0)$ and $u(x,0)$ onto either the piecewise constant
   or piecewise linear FE spaces (analogously to Sect.~\ref{sec_standard} and \ref{sec_split}).

   \paragraph{Structure-preserving nature of the split schemes.}
   
   As shown in Sect.~\ref{sec_splitFE_derivate}, the split FE schemes are structure-preserving in 
   the sense that the corresponding discrete sets of equations preserve the splitting into 
   topological and metric parts. Here, we illustrate that these discrete equations also
   fulfill the first principles of mass and momentum conservation. 
   
   In the continuous case with solutions of the wave equations for velocity $u(x,t)$ and height $h(x,t)$,
   we define the mass $m(t)$ of the water column and it's momentum $p(t)$ by 
   \begin{equation}\label{equ_massmomt}
    m(t) := \int_L h(x,t) dx , \qquad p(t) := \int_L h(x,t) u(x,t) dx .
   \end{equation}
   In fact, the first principles of mass and momentum conservation allow for the derivation 
   of momentum and continuity equations (see e.g.~\cite{Bauer2016}). Therefore, we require our 
   structure preserving schemes to conserve these principles discretely.
      
\begin{figure}[t] \centering
   \begin{tabular}{cc} 
    \hspace{-0.5cm}{\includegraphics[scale=0.52]{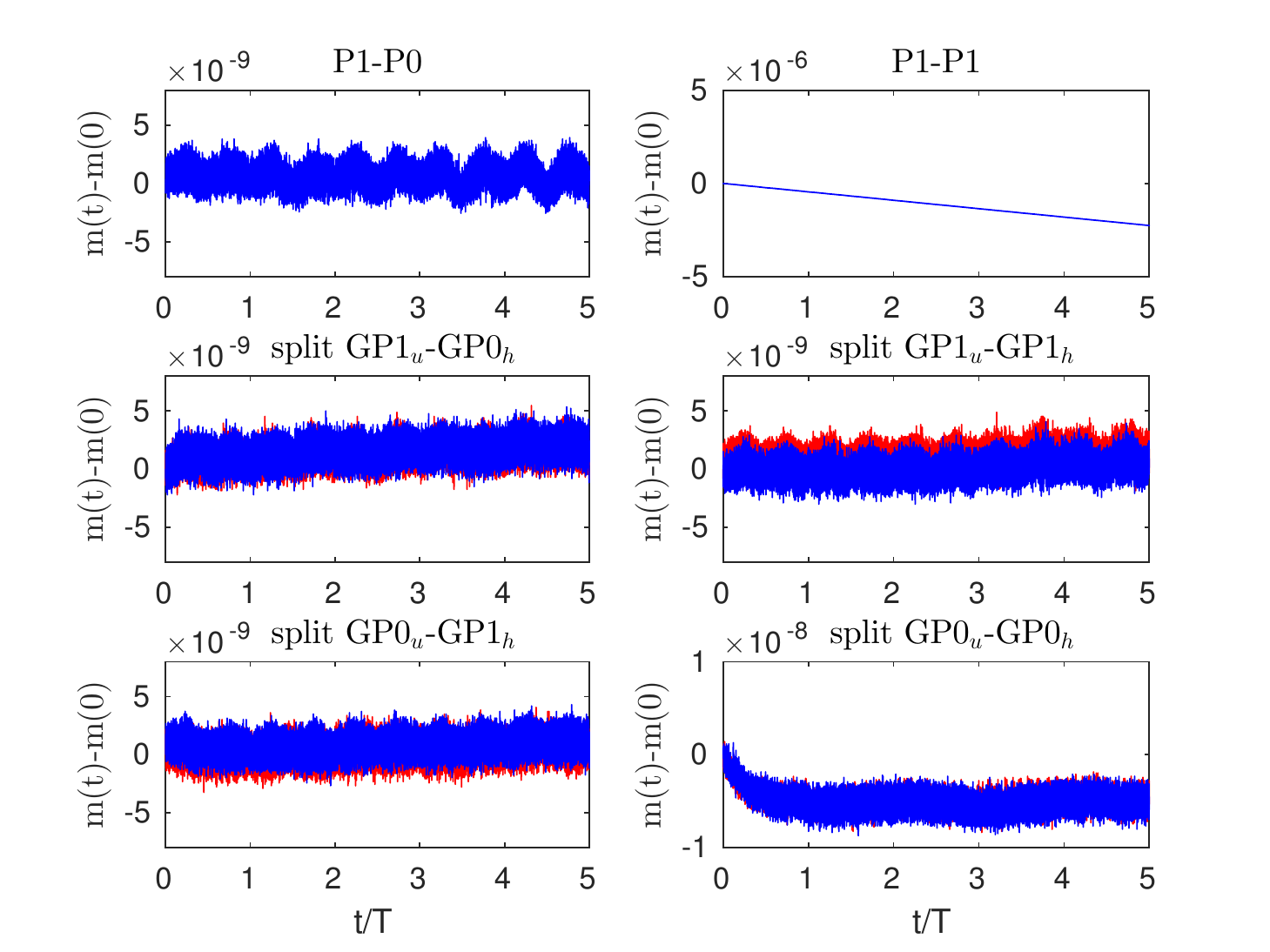}} & 
    \hspace{-0.5cm}{\includegraphics[scale=0.52]{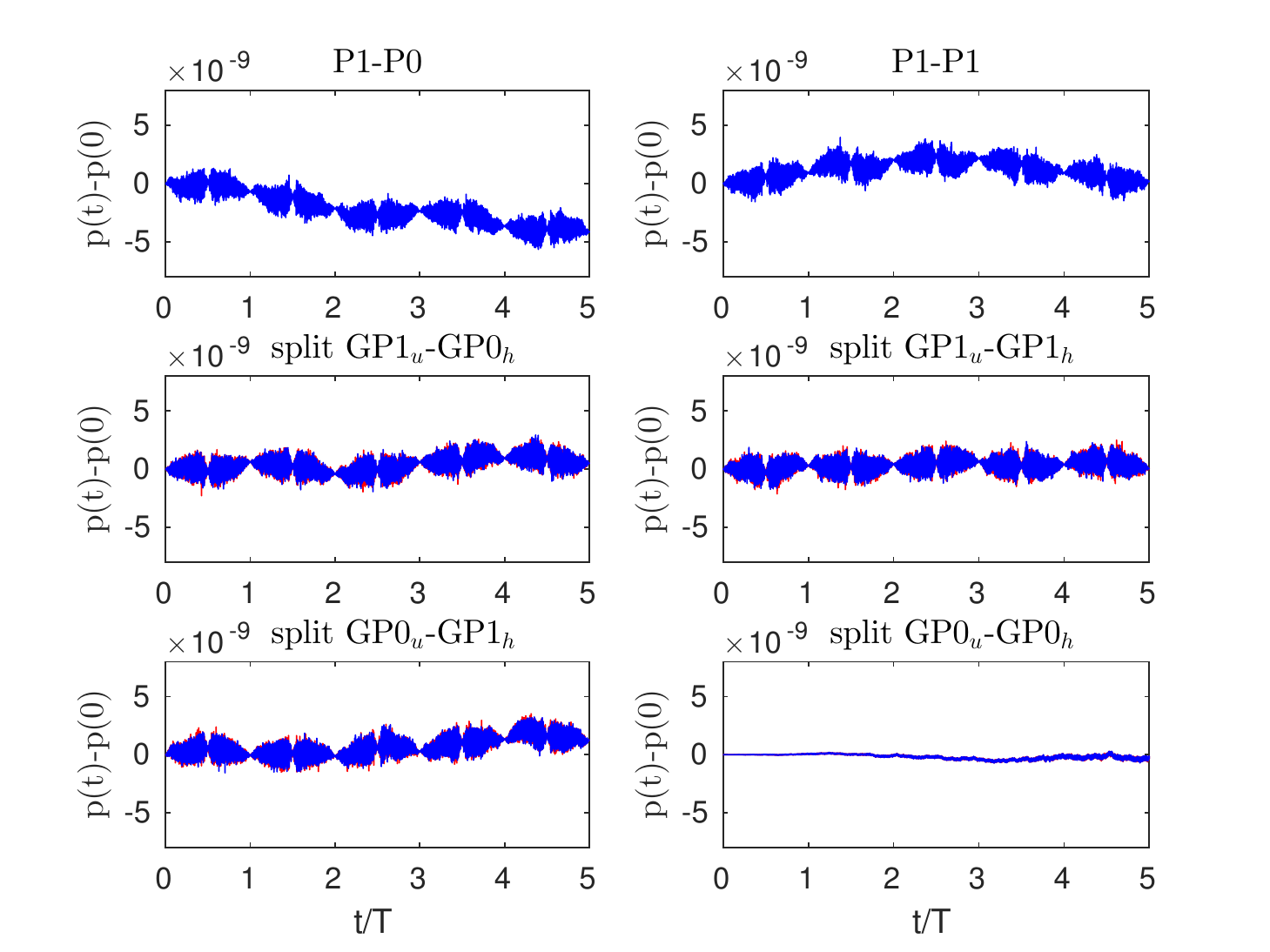}} 
   \end{tabular}
    \caption{Conservation of mass $m(t)$ (left block) and of momentum $p(t)$ (right block) 
    of all schemes for an integration time of $t = 5T$. For the split schemes, the red lines 
    show error values calculated with piecewise constant fields,
    the blue lines with piecewise linear fields.
    }
   \label{fig_conserve}
   \end{figure}

   In Fig.~\ref{fig_conserve}, we present relative errors for mass and momentum 
   with respect to the corresponding initial values for all schemes. The results 
   are determined using TC 2 for an integration time of $t = 5T$ on a mesh with 
   $1024$ elements for the time step sizes mentioned above. 
   For other mesh resolutions and time step sizes, these values remain conserved 
   at the same order of accuracy. More precisely, all schemes expose a mass conservation 
   error (left block) at the order of $10^{-9}$ with the unstable P1--P1 FE scheme 
   being the only exception with an error of approx. $10^{-6}$. The error in the 
   momentum (right block) is of the order of $10^{-9}$ for all schemes.

   \paragraph{Results for TC 1 and TC 2.}
      
\begin{figure}[t] \centering
   \begin{tabular}{cc} 
    \hspace{-0.5cm}{\includegraphics[scale=0.5]{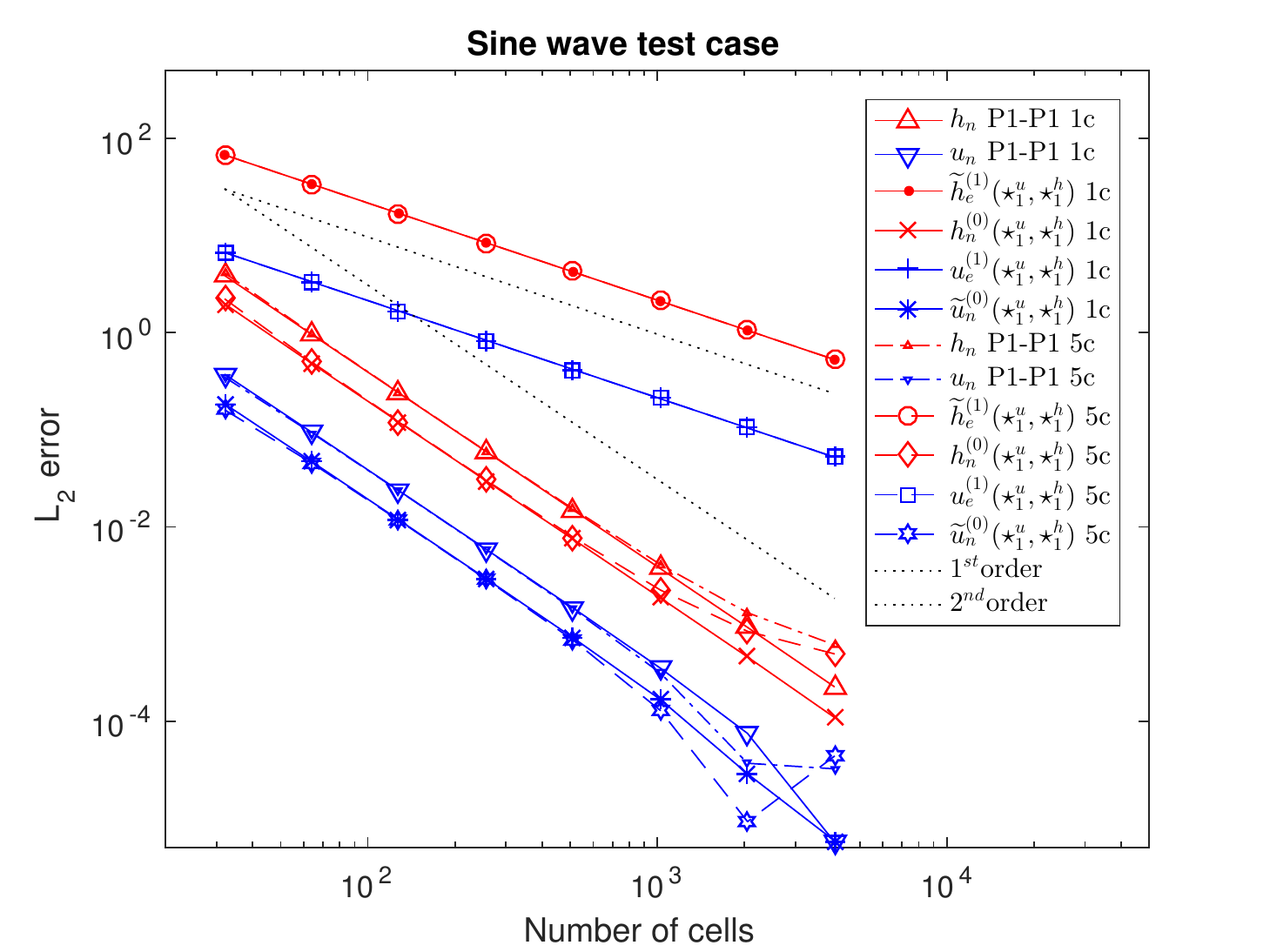}}  & 
    \hspace{-0.5cm}{\includegraphics[scale=0.5]{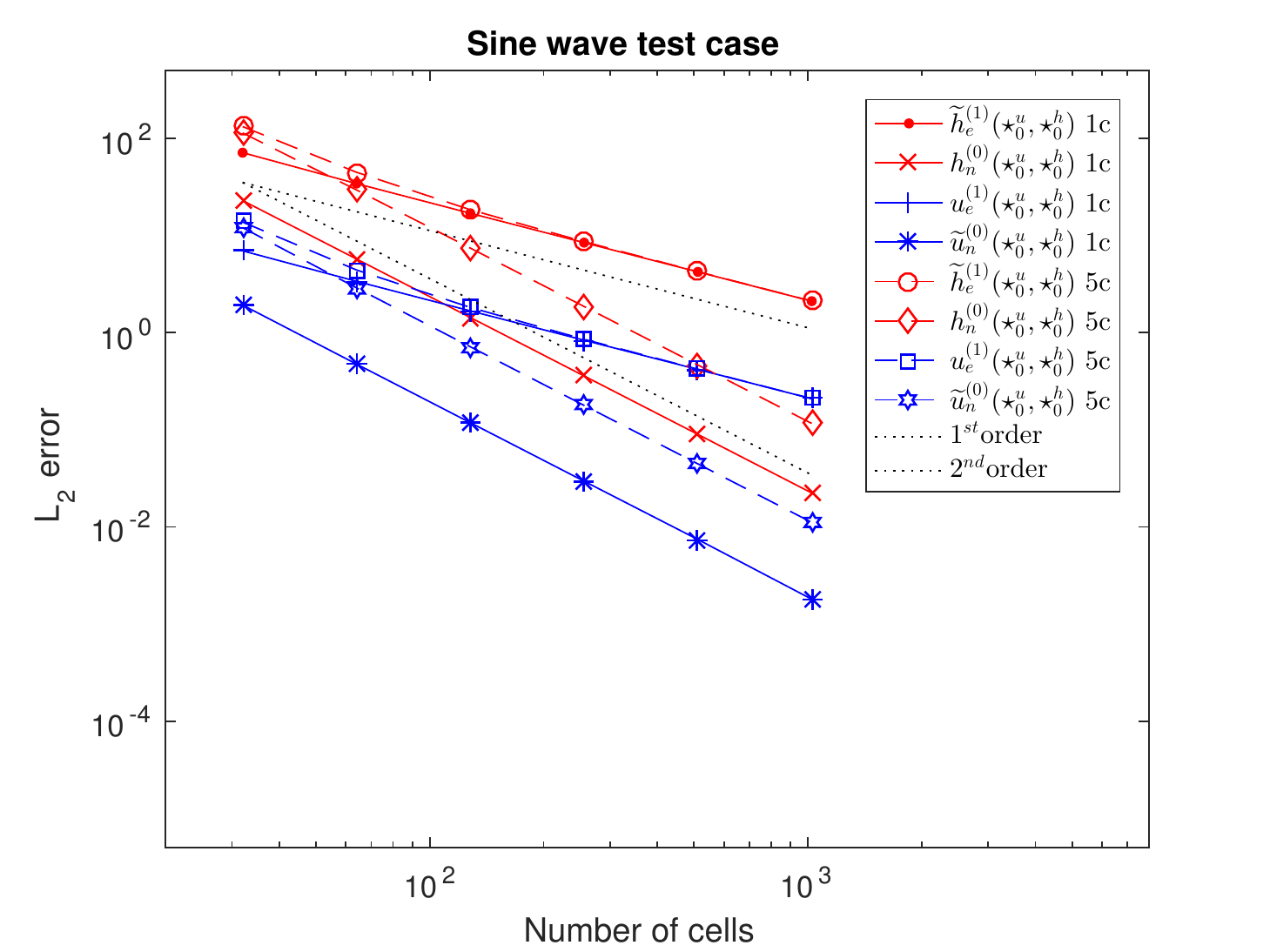}}  \\
    \hspace{-0.5cm}{\includegraphics[scale=0.5]{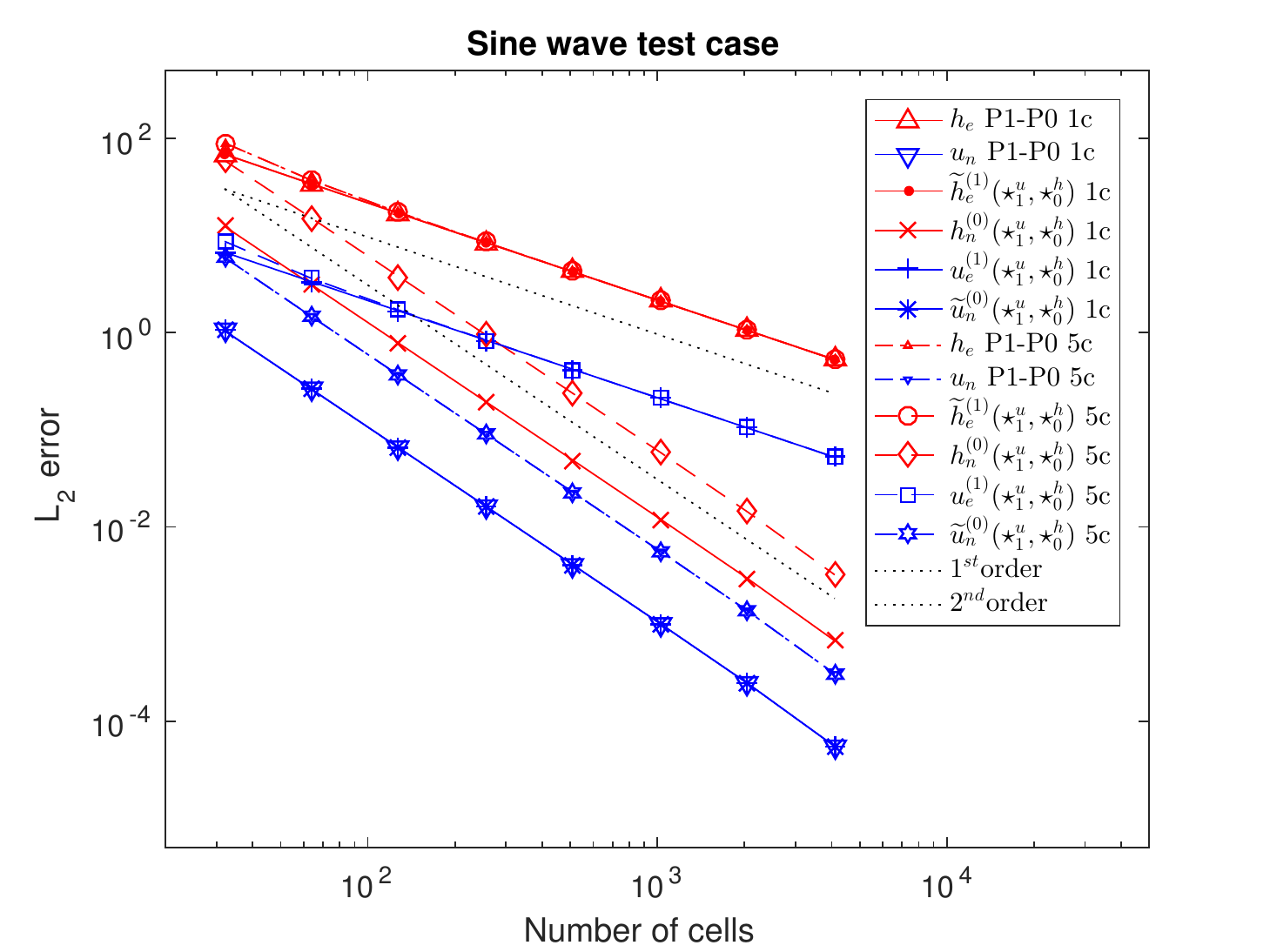}} & 
    \hspace{-0.5cm}{\includegraphics[scale=0.5]{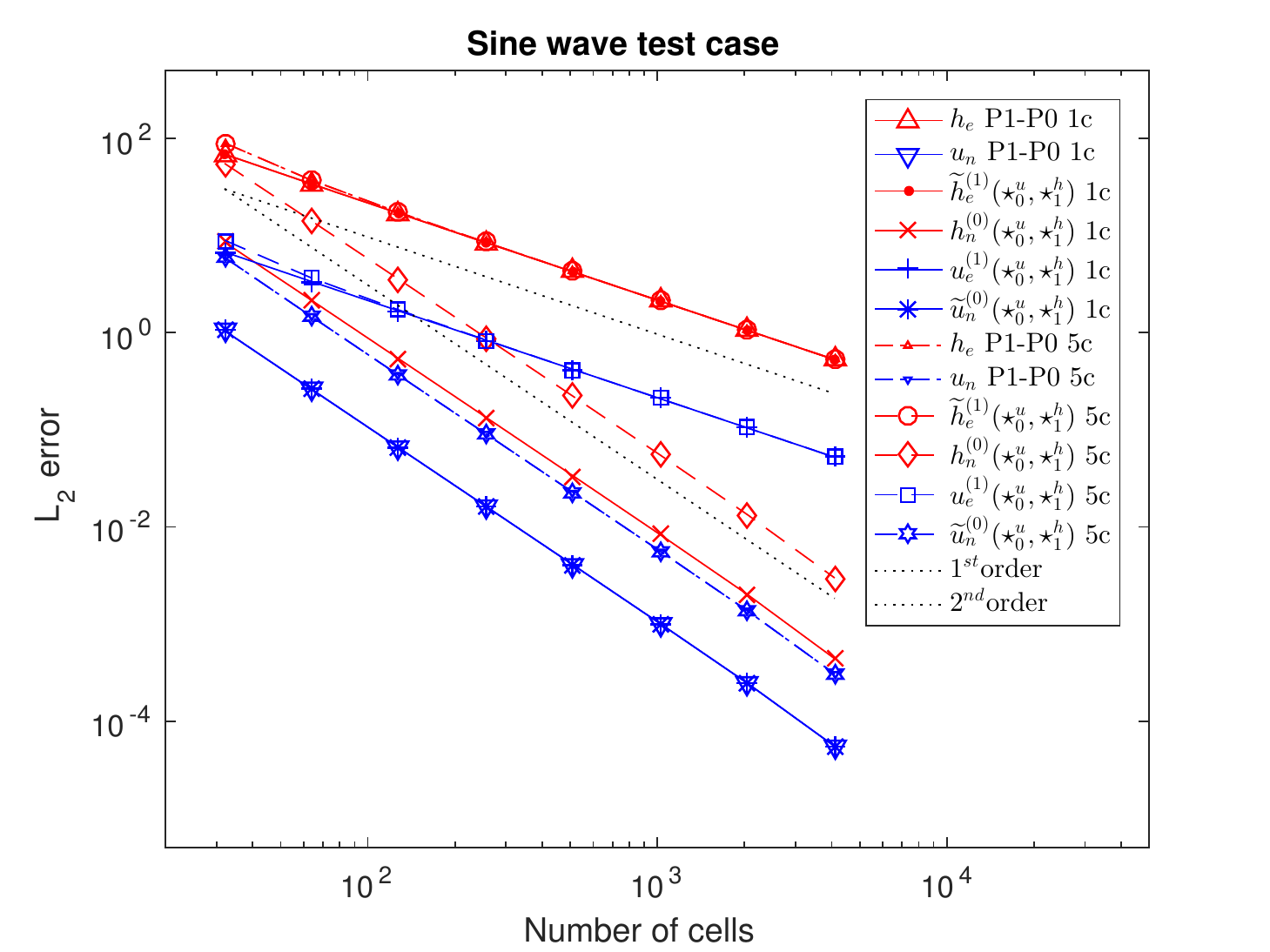}} 
   \end{tabular}
    \caption{Convergence rates of solutions of P1--P1 and $\GPIu$--$\,\GPIh$ (upper left), 
    of $\GPOu$--$\,\GPOh$ (upper right), of P1--P0 and $\GPIu$--$\,\GPOh$ (lower left), 
    and of P1--P0 and $\GPOu$--$\,\GPIh$ (lower right) against analytical sine wave 
    solution~\eqref{equ_h_u_sinus_ana} for integration times $t = 0.875T$ and $t = 4.875T$.
    }
   \label{fig_cnv_sinana}
   \end{figure}

\begin{figure}[t] \centering
  \begin{tabular}{cc} 
    \hspace{-0.5cm}{\includegraphics[scale=0.5]{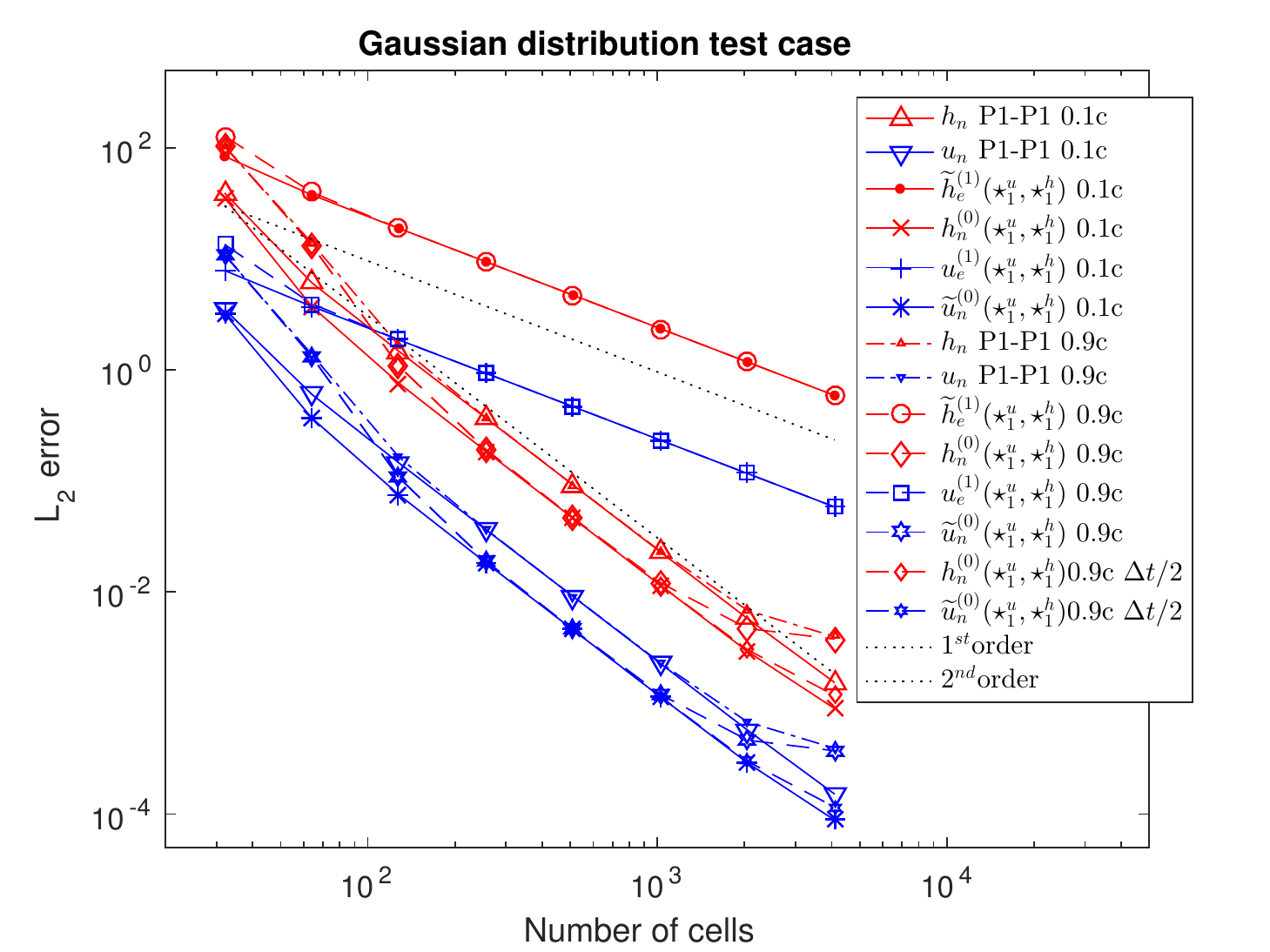}}  & 
    \hspace{-0.5cm}{\includegraphics[scale=0.5]{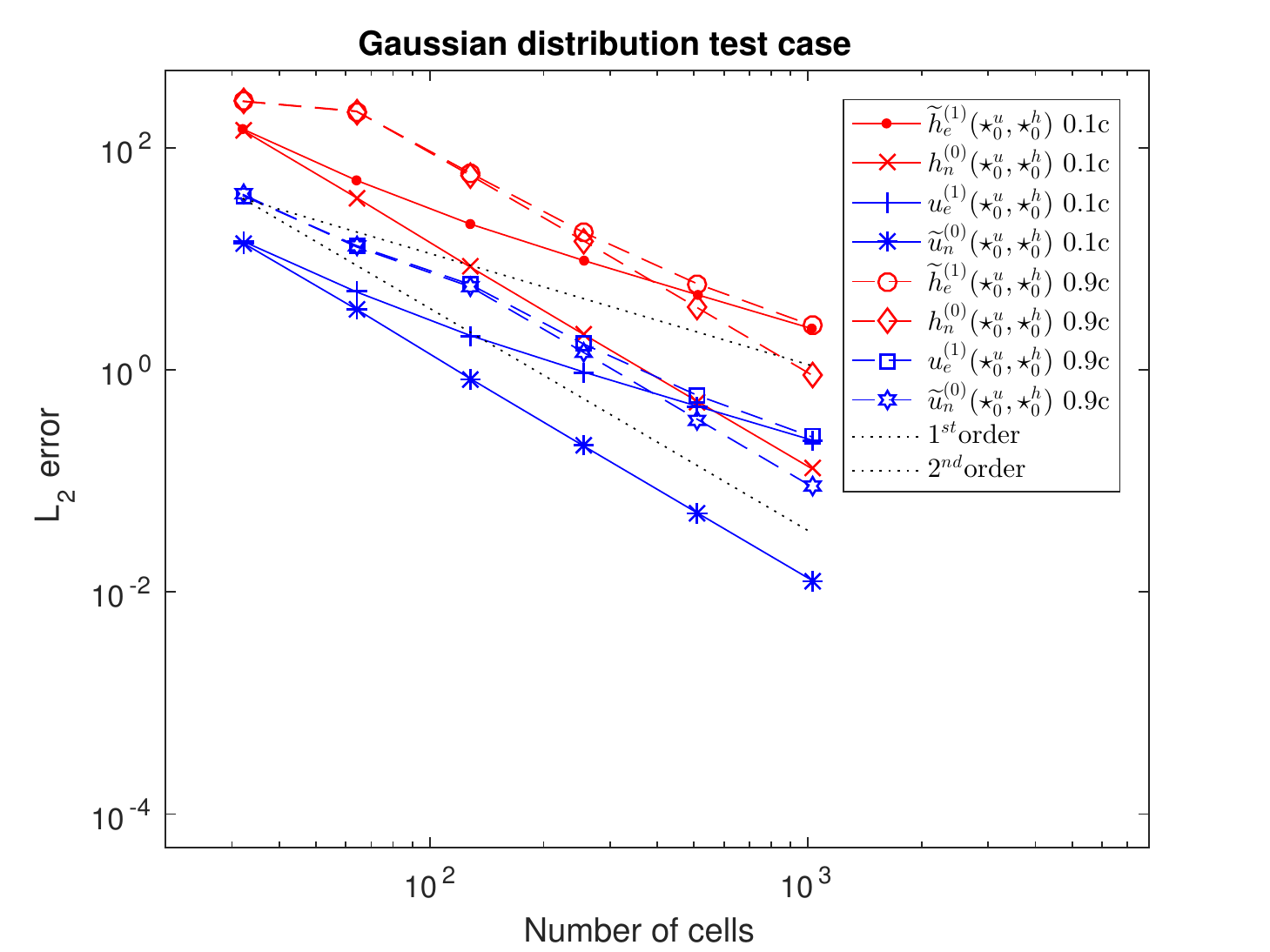}}  \\
    \hspace{-0.5cm}{\includegraphics[scale=0.5]{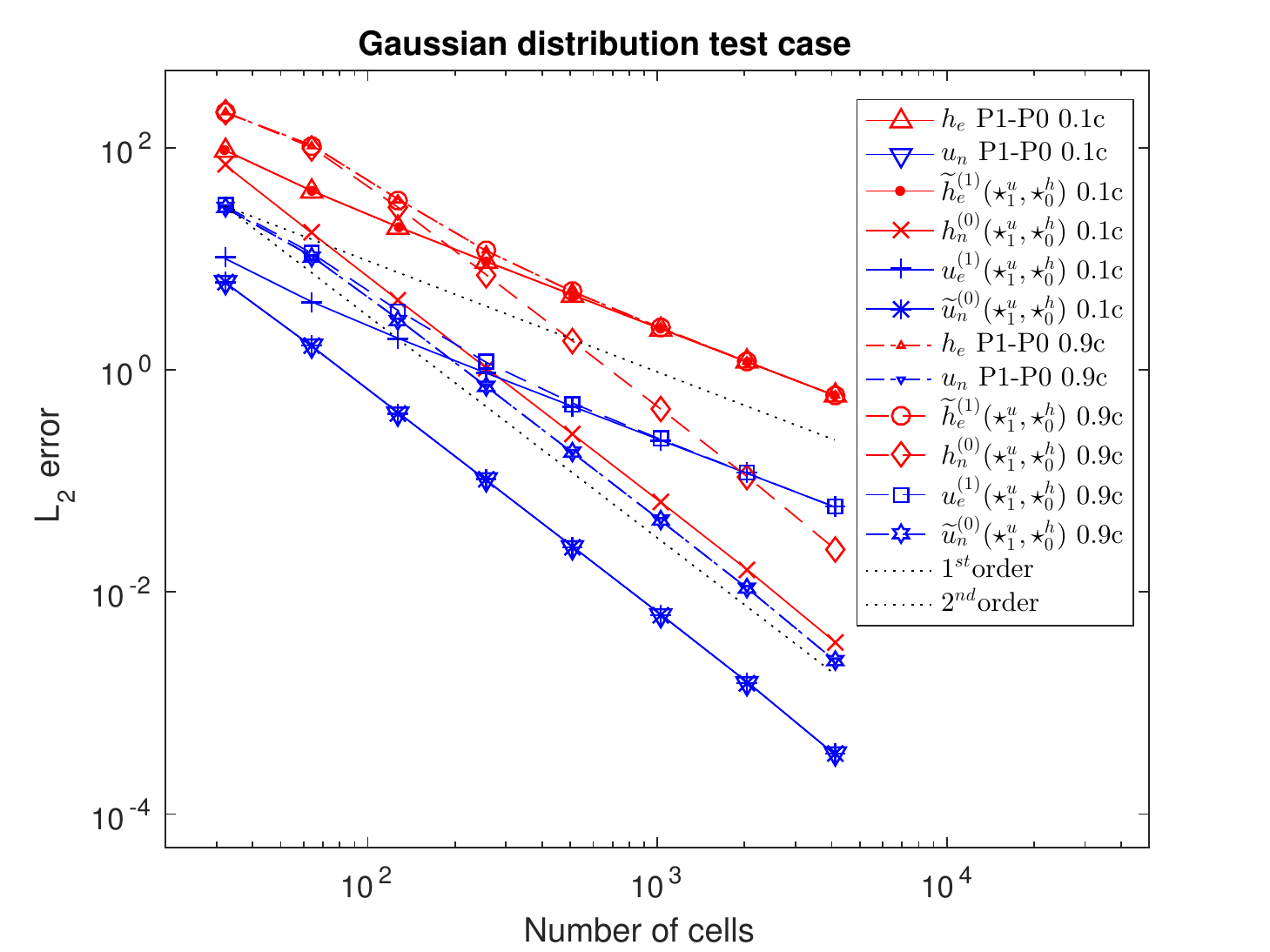}} & 
    \hspace{-0.5cm}{\includegraphics[scale=0.5]{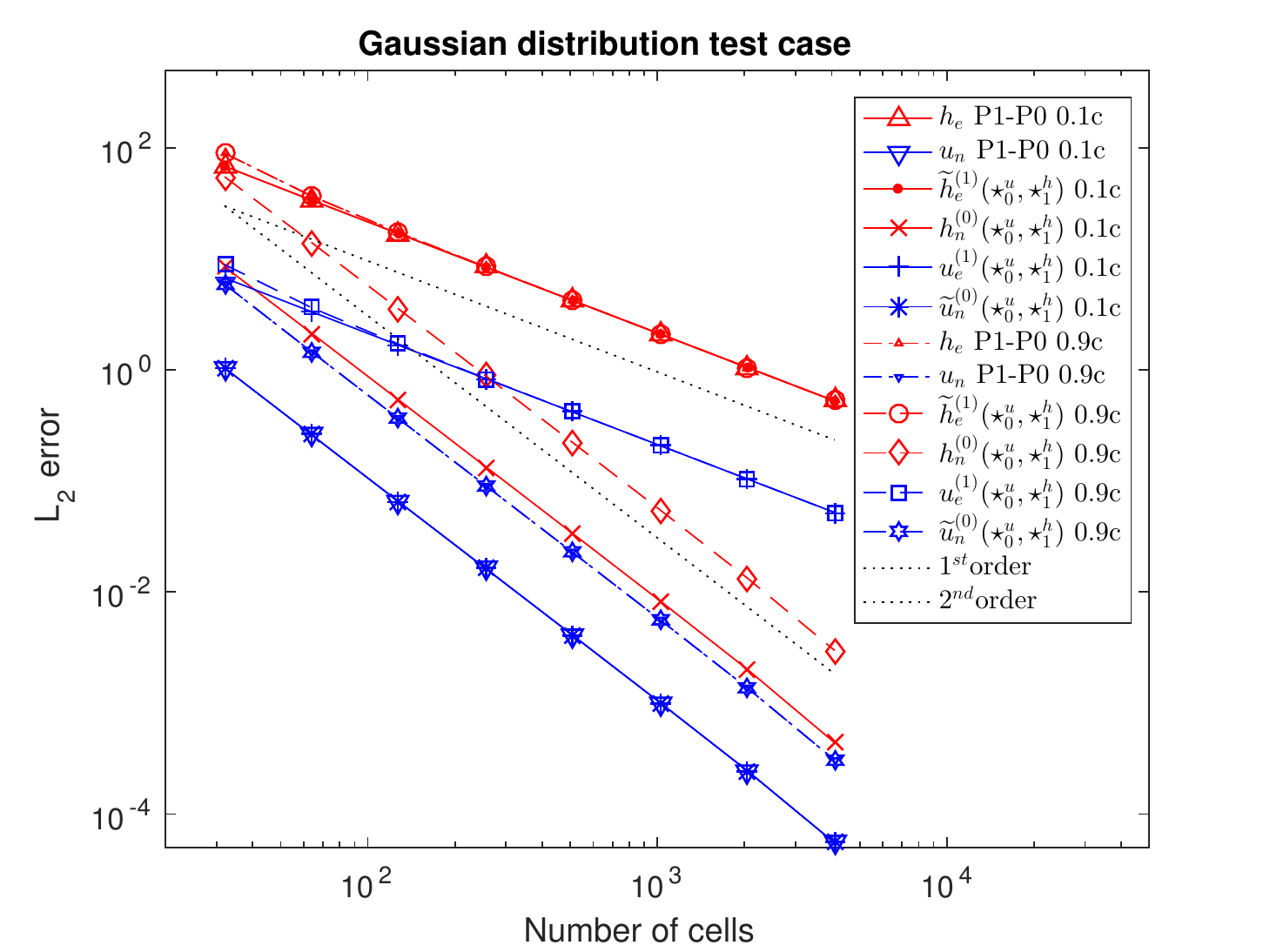}} 
  \end{tabular}
    \caption{Convergence rates of solutions of P1--P1 and $\GPIu$--$\,\GPIh$ 
    (upper left), of $\GPOu$--$\,\GPOh$ (upper right), of P1--P0 and $\GPIu$--$\,\GPOh$ (lower left), 
    and of P1--P0 and $\GPOu$--$\,\GPIh$ (lower right) 
    against analytical Gaussian wave solution~\eqref{equ_h_u_gauss_ana}
    for integration times $t = 0.125T$ and $t = 0.875T$. 
    }
  \label{fig_cnv_gaussdist01c09c}
  \end{figure}

   Figure~\ref{fig_cnv_sinana} shows convergence towards the analytical sine wave 
   solutions~\eqref{equ_h_u_sinus_ana} for integration times $t = 0.875T$ and $t = 4.875T$.
   Similarly, in Figure~\ref{fig_cnv_gaussdist01c09c} convergences to the analytical 
   Gaussian-shaped solutions~\eqref{equ_h_u_gauss_ana} for integration times 
   $t = 0.125T$ and $t = 0.875T$ are given. To determine convergence rates valid for any smooth solution 
   given in the form of equations~\eqref{equ_h_u_sinus_ana} and \eqref{equ_h_u_gauss_ana}, 
   we use values for $t$ that are not multiples 
   of $T/4$ such that neither height nor velocity solutions 
   coincide with constant functions.

   The reduction in integration time from about 5 cycles for TC 1 down to about 1 cycle for TC 2 
   is a consequence of the errors caused by the discrete dispersion relations. These
   lead to a deformation of the Gaussian-shaped wave package and to the development of 
   high oscillatory waves, which prevents us from determining convergence rates in case of long integration times. 
   
   Considering the absolute error values of all split schemes for both test cases,
   we notice that the low accuracy ($\GPOu$--$\,\GPOh)$ scheme has in general 
   significantly higher error values than both medium accuracy ($\GPIu$--$\,\GPOh/\GPOu$--$\, \GPIh$)
   schemes, which in turn have higher error values than the high accuracy ($\GPIu$--$\,\GPIh$)
   one. This can be expected as the low accuracy scheme applies two low order 
   closure equations, the medium schemes a combination of low and higher order, and the high
   accuracy scheme two higher order closure equations. Hence, these absolute error values indeed
   reflect how accurately the two metric closure equations are discretized. 
   
   Comparing the split with the mixed schemes, we observe a general agreement of the 
   convergence rates of the corresponding solutions. In more detail, those fields represented by 
   piecewise constant basis functions show the expected $1^{st}$-order convergence rate 
   for all schemes. Moreover, the absolute error values of corresponding mixed and split schemes 
   agree very well. As we are not aware of a mixed FE counterpart for the low accuracy case, 
   we cannot compare it to an established FE method.

   The fields represented by piecewise linear bases show the expected $2^{nd}$-order convergence rates. 
   Here, error values calculated with the split $\GPIu$--$\,\GPIh$ scheme are a little smaller then 
   those corresponding to solutions of the P1--P1 scheme,
   but the solutions obtained from each scheme for short and long integration times agree very well. 
   In contrast, the $L_2$ error values for the $\GPOu$--$\,\GPOh$ scheme for the 
   fields with piecewise linear representation increase with integration time. This property is shared 
   by both $\GPIu$--$\,\GPOh$ and $\GPOu$--$\,\GPIh$ schemes and the mixed P1--P0 scheme, while  
   the absolute error values of mixed and split schemes agree very well. 
   When the $L_2$ errors between mixed and split schemes show small differences in low resolution meshes, 
   they usually vanish with higher resolution.

   In case of mixed P1--P1 and split $\GPIu$--$\,\GPIh$ schemes, we notice for high resolutions of about $2048$ and $4096$ 
   elements with error values smaller than $10^{-4}$ a discrepancy from the $2^{nd}$-order convergence rate.  
   At such high resolution the spatial errors are in the order of the temporal errors leading
   to a flattening of the convergence curve, or even to some jumps as in Fig.~\ref{fig_cnv_sinana} (upper-left).
   This is confirmed by halving the time step size, which reduces the time step errors by a factor of $4$, which leads to the 
   improved convergence rates as shown in (Fig.~\ref{fig_cnv_gaussdist01c09c}, upper-left)). 
   The factor $4$ corresponds well to the expected $2^{nd}$-order convergence rate of the 
   CN time discretization scheme, introduced in Sect.~\ref{sec_split_CN}.

   Both $\GPOu$--$\,\GPOh$ and $\GPIu$--$\,\GPOh/\GPOu$--$\,\GPIh$ schemes exhibit straight convergence lines of 
   second order because here the temporal errors are smaller than the spatial errors. 
   The $\GPOu$--$\,\GPOh$ scheme exposes spatial errors one order of magnitude larger than all 
   other cases with resolution up to $1024$ elements, even with very small time steps (see explanation below).
   For the $\GPIu$--$\,\GPOh/\GPOu$--$\,\GPIh$ schemes, there is no flattening in the convergence curves
   as these schemes approximate the continuous dispersion relation better than the other schemes 
   (cf. Fig.~\ref{fig_disp_relation}), which keeps the error in the time integration
   lower than in the unstable cases.

  \paragraph{Results of TC 3.}
  
  Applying the Gaussian distribution~\eqref{equ_h_u_gauss_ana} with very small width (TC 3), 
  we illustrate the schemes' representation of waves and their dispersions. To this end, 
  we let the simulations run up to $t = 0.1T$ and study how the Gaussian-shaped wave packages 
  change their forms during simulations. These changes are a consequence of the different
  discrete dispersion relations for the schemes studied in this manuscript. 
  Note that the short simulation time is justified by the design of TC 3 such that 
  the numerical behavior in the short wave regime 
  (right side in figure \ref{fig_disp_relation}) is exposed.
  
\begin{figure}[h] \centering
  \begin{tabular}{ccc} 
   \begin{overpic}[scale=0.5,unit=1mm]
    {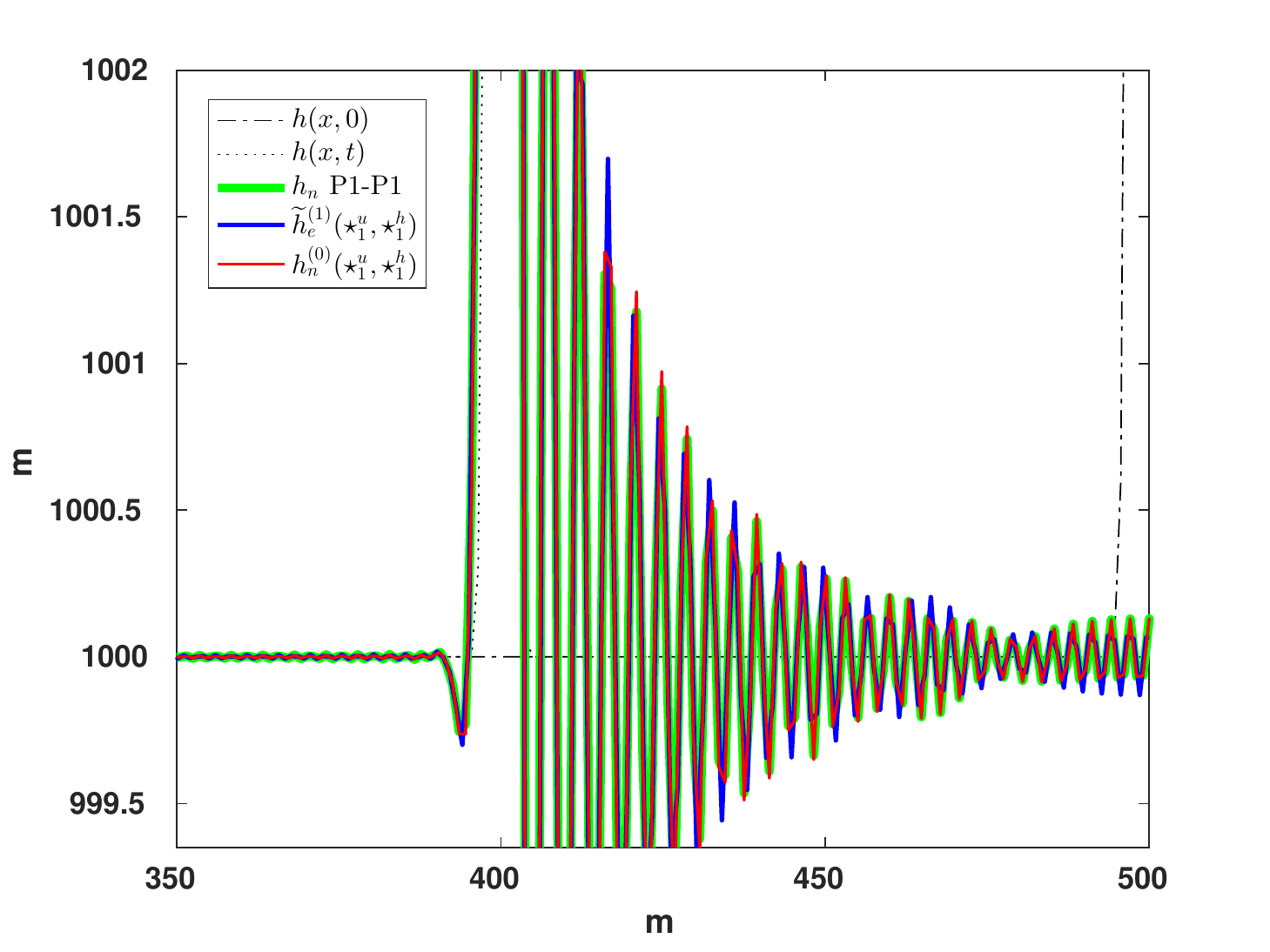}
     \put(40.2,27.8){\includegraphics[scale=.25]{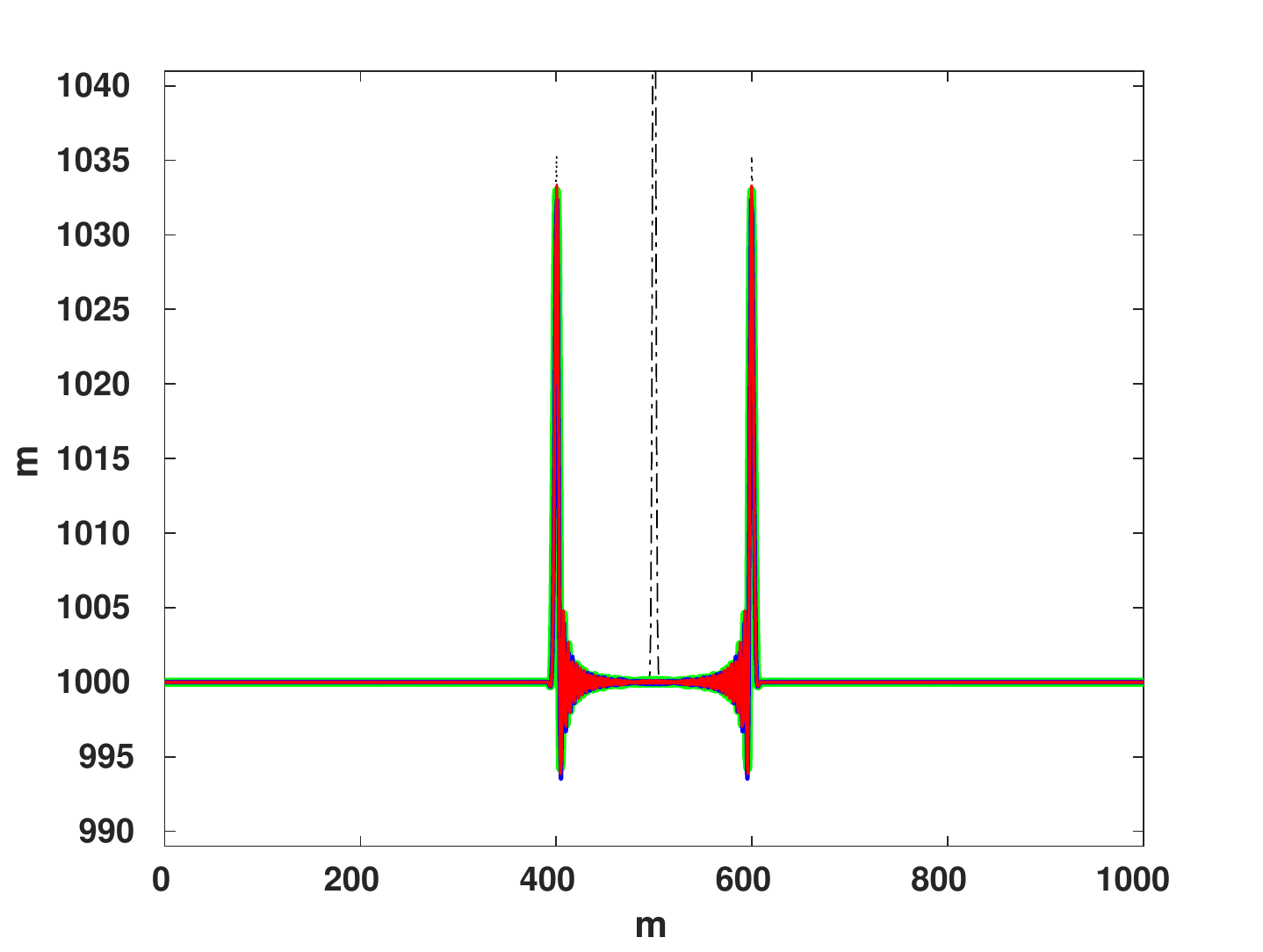}}
   \end{overpic}
   \begin{overpic}[scale=0.5,unit=1mm]
    {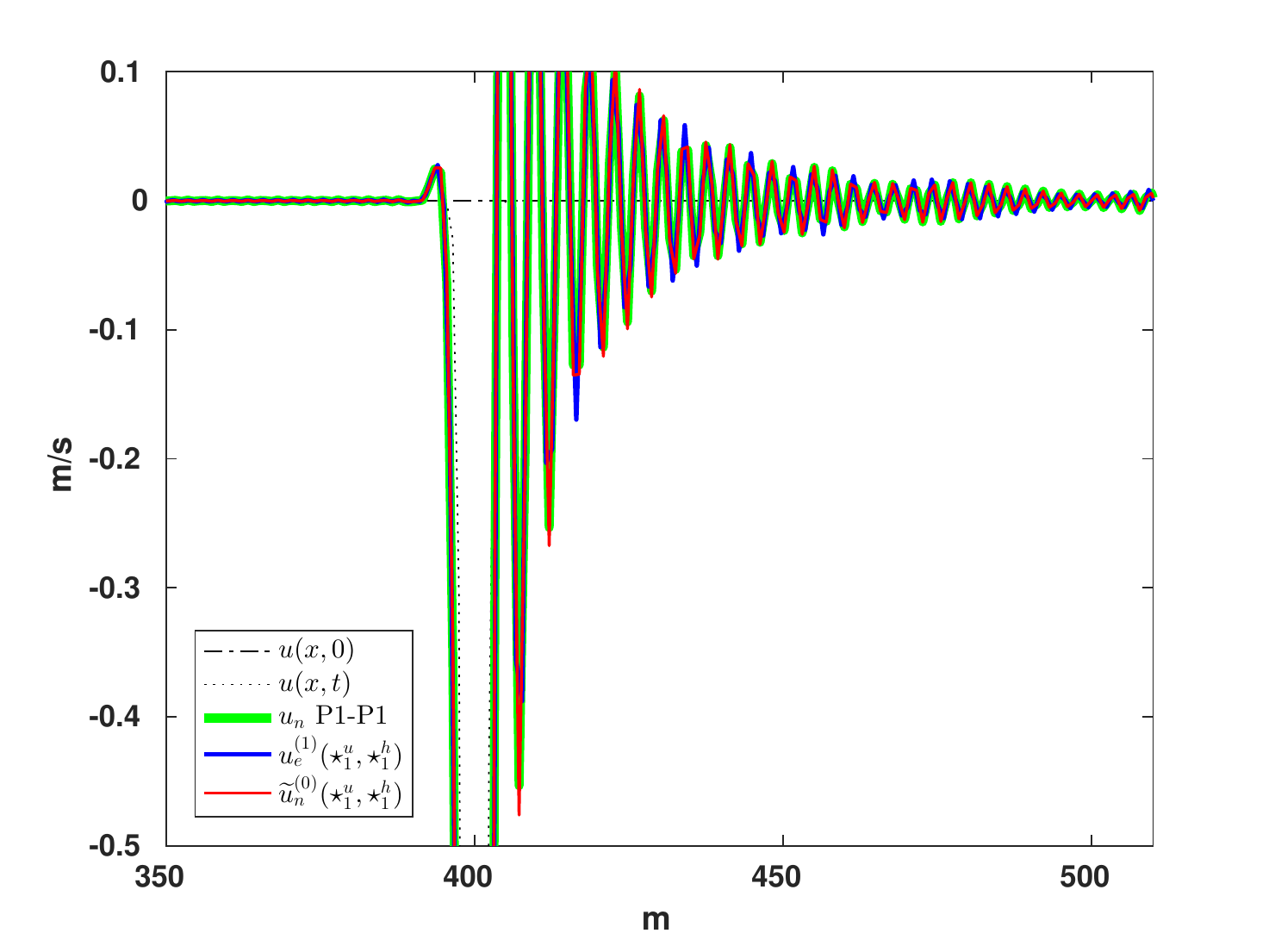}
      \put(34.2,7.2){\includegraphics[scale=.27]{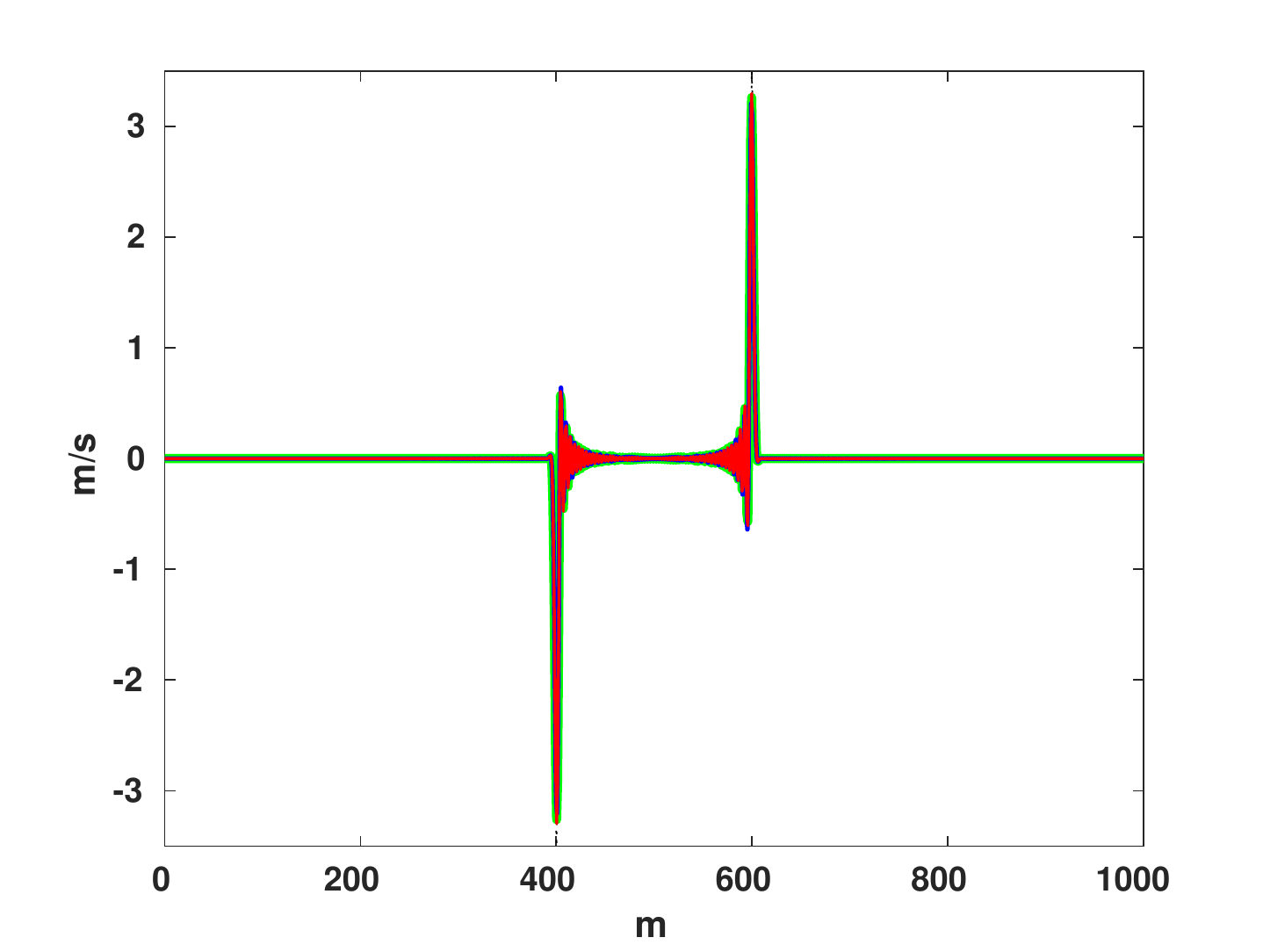}}
   \end{overpic}
   \end{tabular}
   \caption{Illustration of the effects of the discrete dispersion relation~\eqref{equ_disprel_split00} 
    of the unstable P1--P1 and $\GPIu$--$\,\GPIh$ schemes on Gaussian-shaped height and velocity peaks 
    after a simulation time of $t=0.1T$.}  
    \label{fig_disprel_spltv0h0}
  \end{figure}

  Considering first the unstable mixed P1--P1 and split $\GPIu$--$\,\GPIh$ schemes in 
  Figure~\ref{fig_disprel_spltv0h0}, we notice that the piecewise linear solutions of
  both mixed and split schemes are very similar. At time $t = 0.1T$,
  the left wave package travels to the left, the right one to the right. 
  The front position on both sides agree well with the analytical solution
  and there are no forerunning fast waves. However, trailing oscillations develop, 
  which change the shape of the Gaussian distribution, rendering it lower and wider. 
  Small scale oscillations occupy the entire region between the two fronts.

  As predicted by the discrete dispersion relation~\eqref{equ_disprel_split00} (or by \eqref{equ_disprel_P1P1}), 
  the wave with shortest wave length $k = \frac{\pi}{\Delta x}$ has zero phase velocity. 
  This spurious mode that is not transported away from the center of the domain but oscillates
  on the grid scale level, is visible in Figure~\ref{fig_disprel_spltv0h0}. 
  In general, as derived theoretically, 
  our simulations show that the phase velocity of all waves are bounded from above by 
  $c_d \leq c$, and with increasing wave numbers $k$, the wave speeds get slower and 
  slower up to the shortest wave with zero phase velocity (cf. Fig.~\ref{fig_disp_relation}).

\begin{figure}[h] \centering
  \begin{tabular}{ccc} 
   \begin{overpic}[scale=0.5,unit=1mm]
    {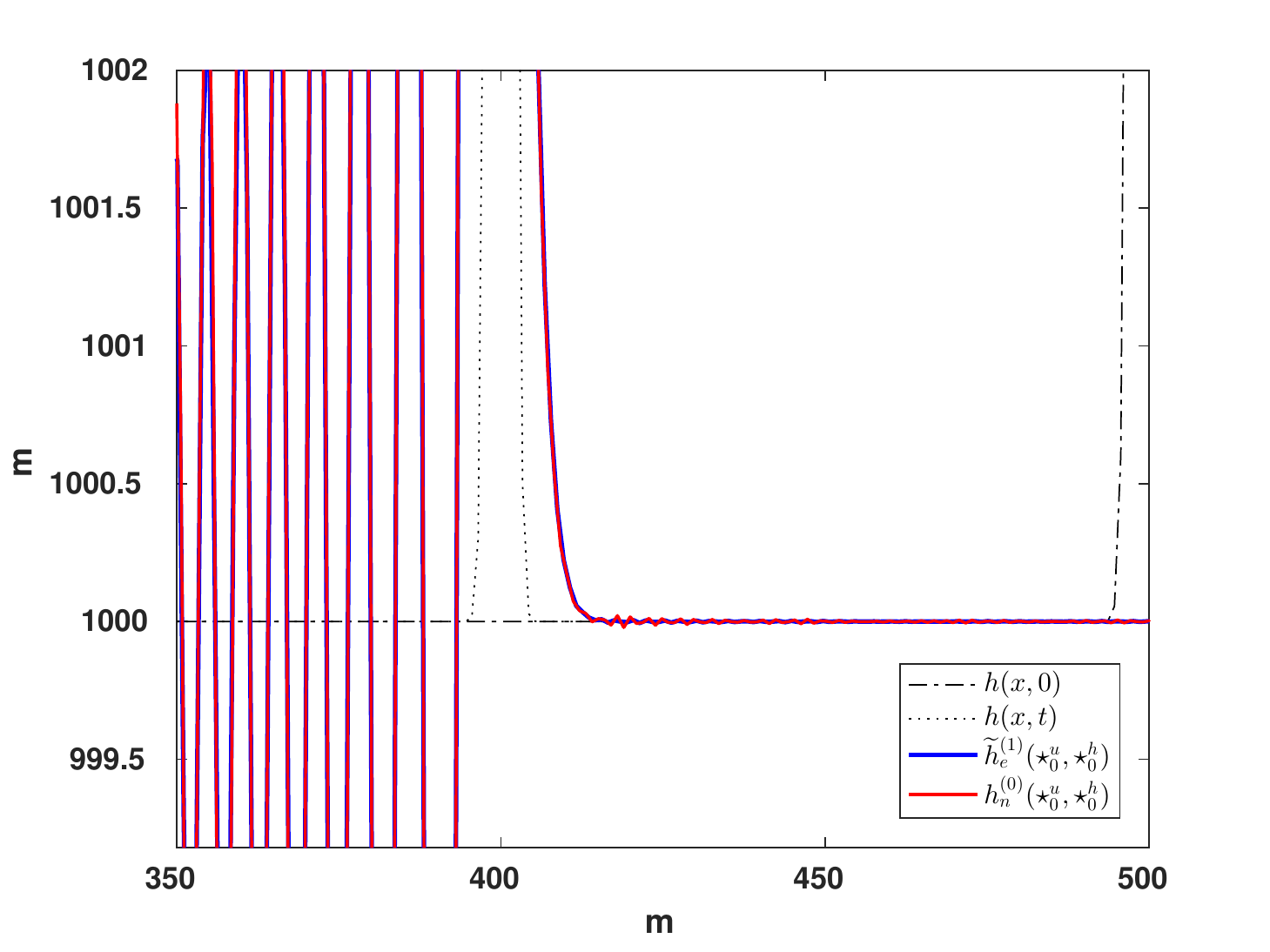}
     \put(33,26.2){\includegraphics[scale=.27]{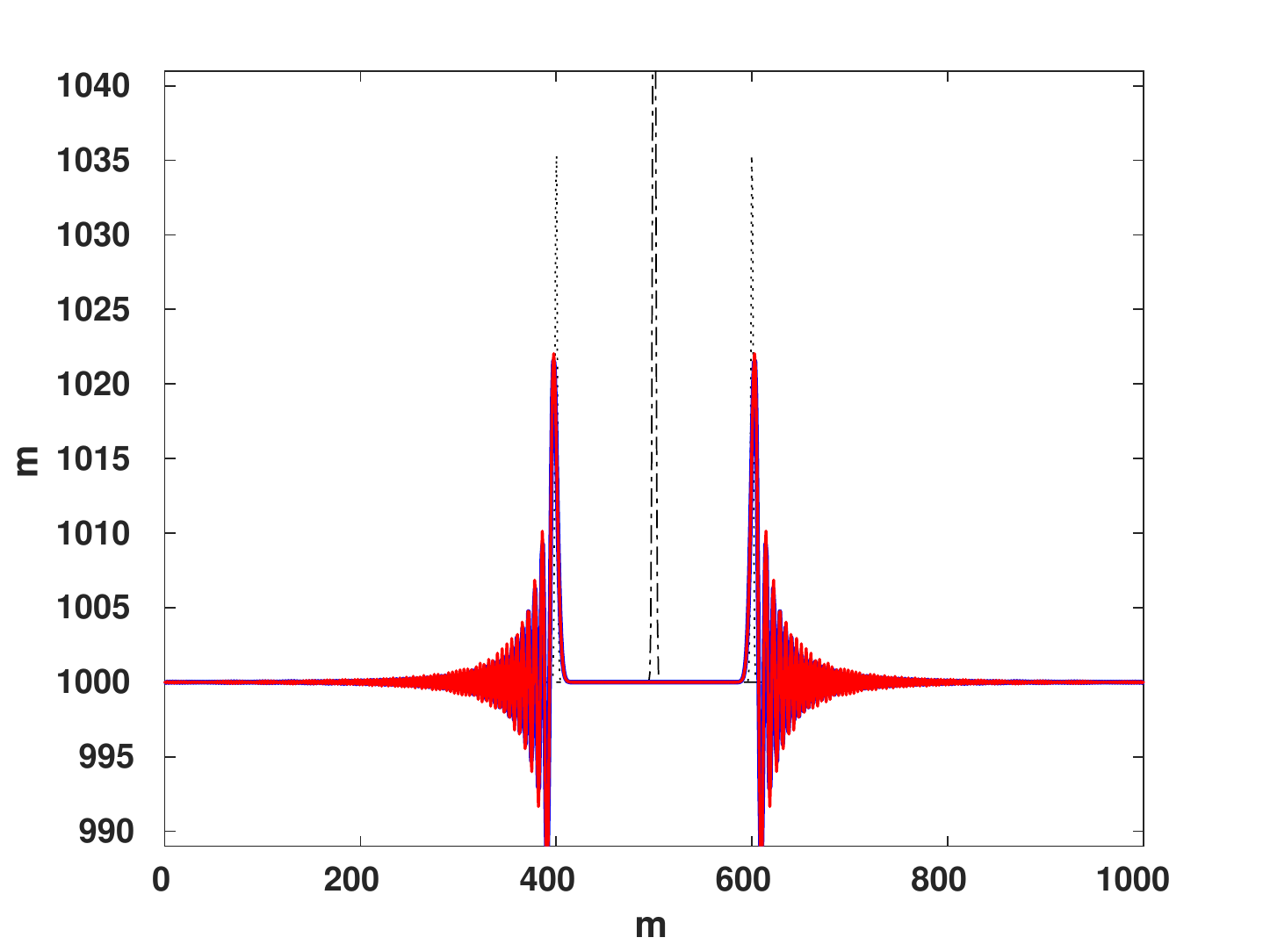}}
   \end{overpic}
   \begin{overpic}[scale=0.5,unit=1mm]
    {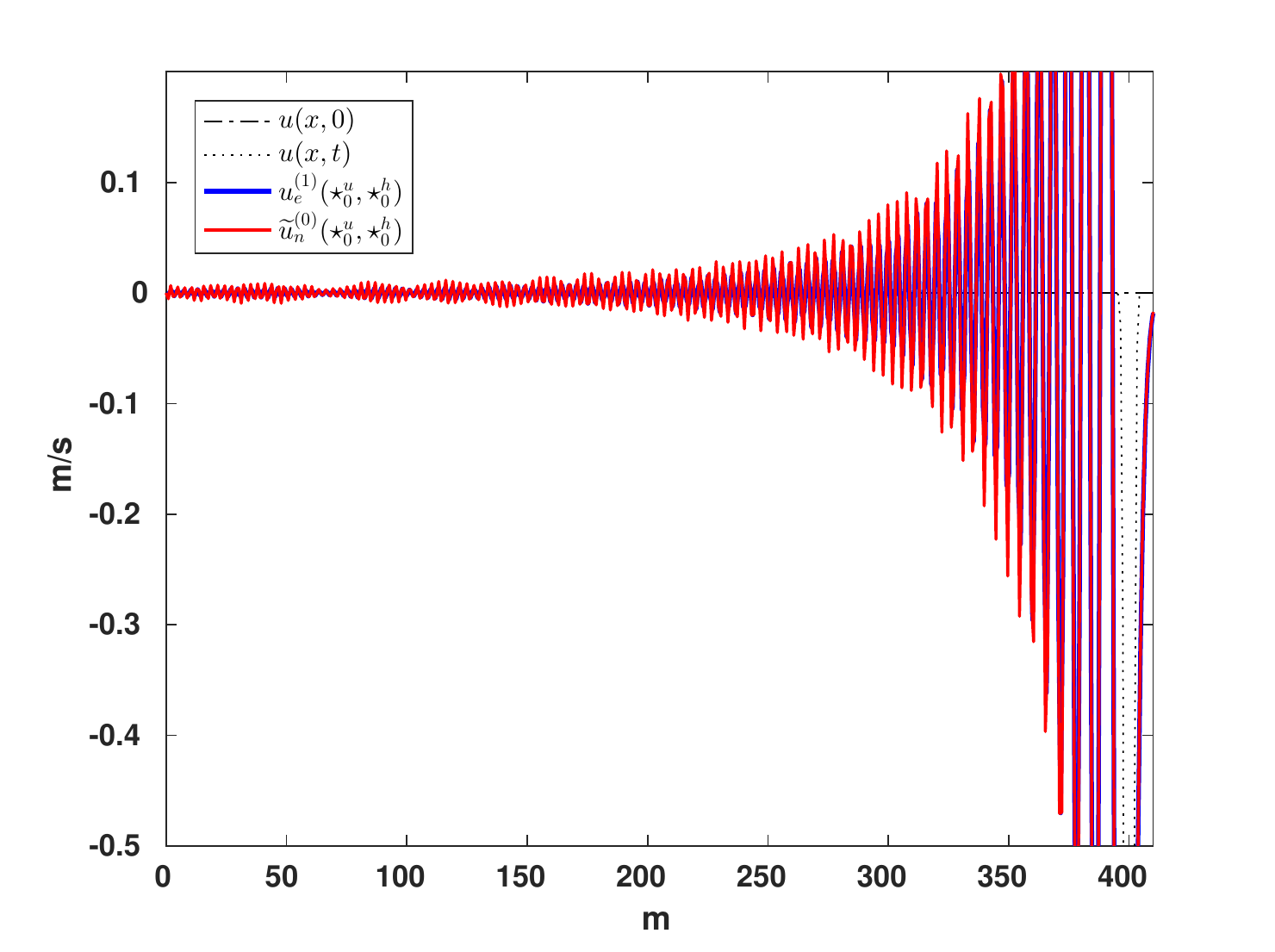}
      \put(9.9,6.3){\includegraphics[scale=.25]{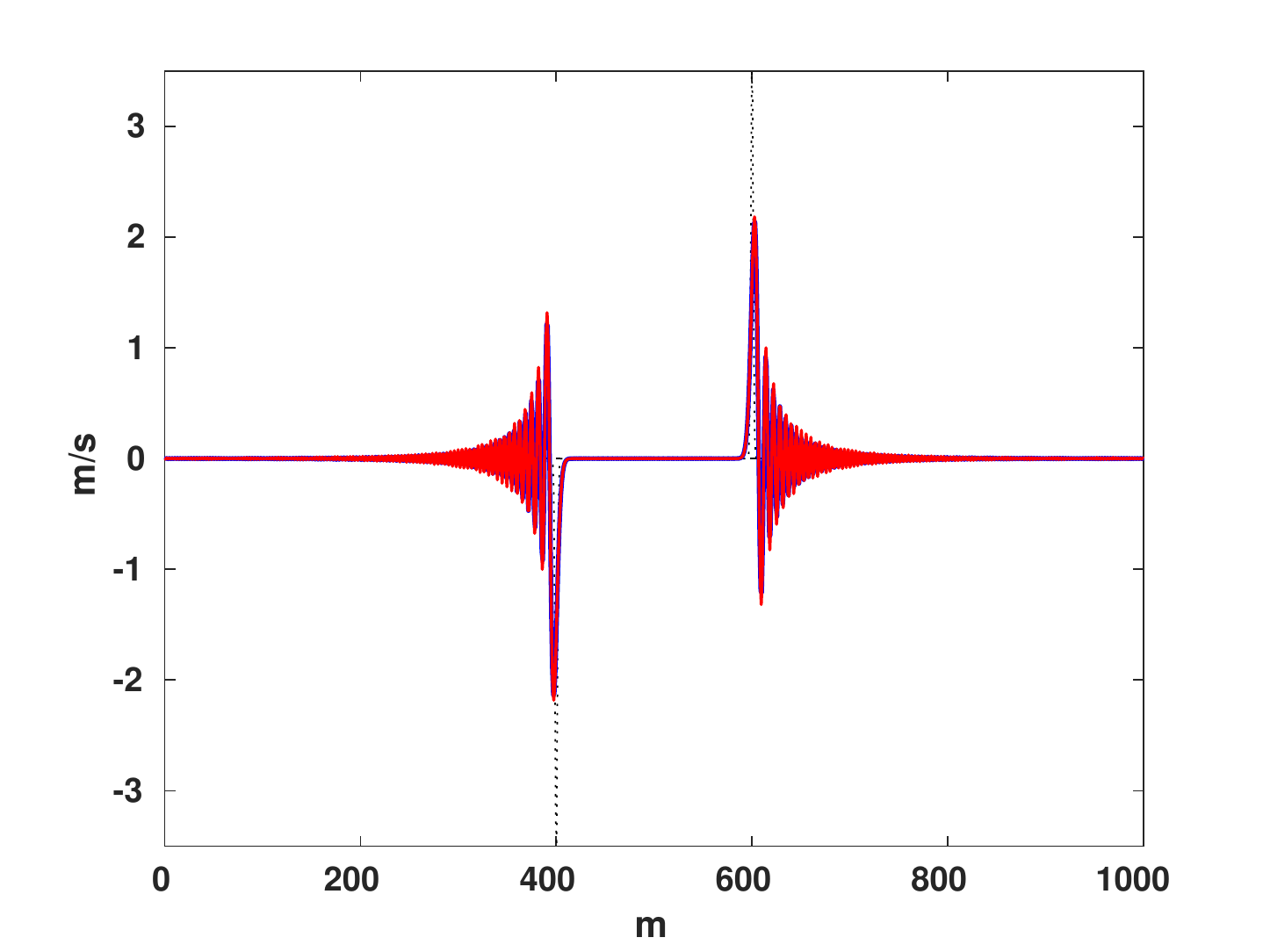}}
   \end{overpic}
   \end{tabular}
    \caption{Illustration of the effects of the discrete dispersion relation~\eqref{equ_disprel_split11} 
    of the $\GPOu$--$\,\GPOh$ scheme on Gaussian-shaped height and velocity peaks 
    after a simulation time of $t=0.1T$.}
  \label{fig_disprel_spltv1h1}
  \end{figure}

  Looking at the split $\GPOu$--$\,\GPOh$ scheme (cf. Figure~\ref{fig_disprel_spltv1h1}), 
  the wave packages show a similar general behavior as in the previous case. However, the differing
  discrete dispersion relation~\eqref{equ_disprel_split11} leads to differently 
  deformed Gaussian-shaped wave packages compared to the previous unstable schemes.
  As above, the peaks become smaller and wider. However, in the $\GPOu$--$\,\GPOh$ scheme
  the central part between the packages is free from high frequency modes, 
  which indicates that the scheme does not support a standing spurious mode. 
  On the other hand the high frequency spurious waves
  occur outside the central part, hence the scheme supports a fast traveling spurious mode. 
  In fact, all waves travel with speeds equal to or larger than $c$, the analytically correct 
  wave speed, and increase with wave number $k$ up to very large wave speeds. 
  This is why we find waves occupying the entire outer region, 
  and why our time step size is so much smaller when compared to the other schemes 
  (cf. Sect.~\ref{sec_discuss_stability}).
  This observation agrees very well with the discrete dispersion relation for the split 
  $\GPOu$--$\,\GPOh$ scheme shown in Fig.~\ref{fig_disp_relation}.

\begin{figure}[h] \centering
  \begin{tabular}{ccc} 
   \begin{overpic}[scale=0.5,unit=1mm]
    {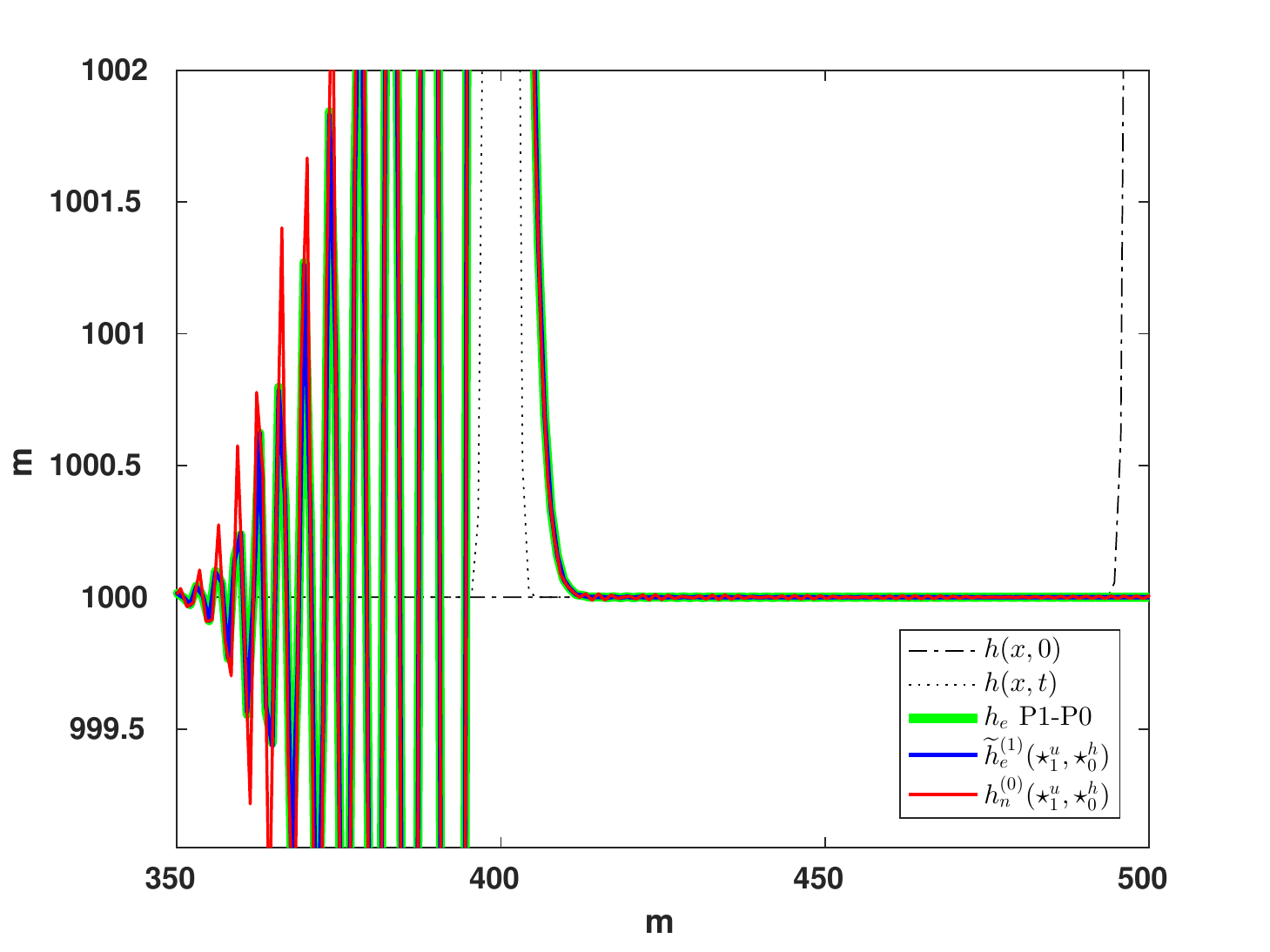}
     \put(33,26.2){\includegraphics[scale=.27]{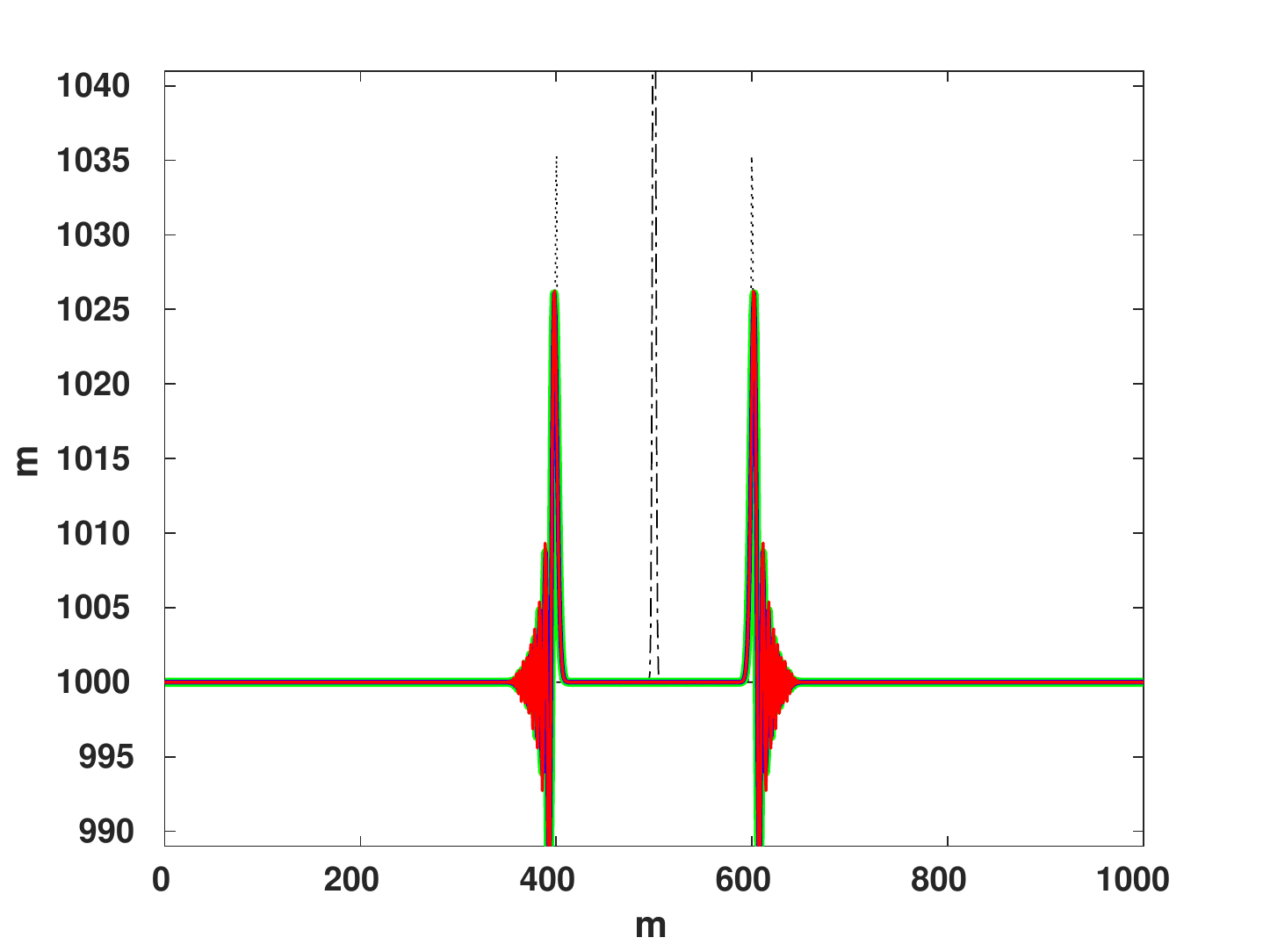}}
   \end{overpic}
   \begin{overpic}[scale=0.5,unit=1mm]
    {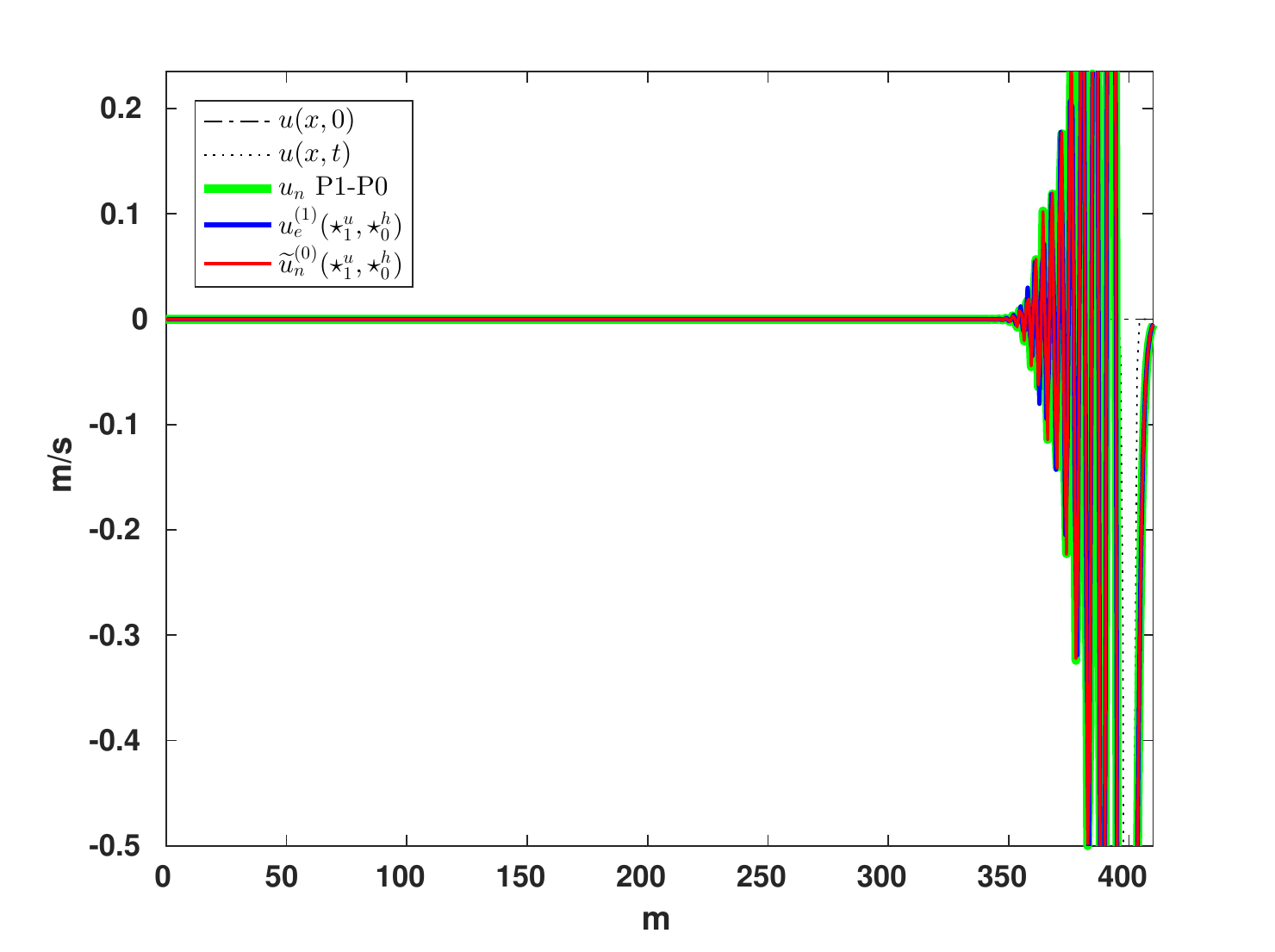}
      \put(14.2,6.6){\includegraphics[scale=.27]{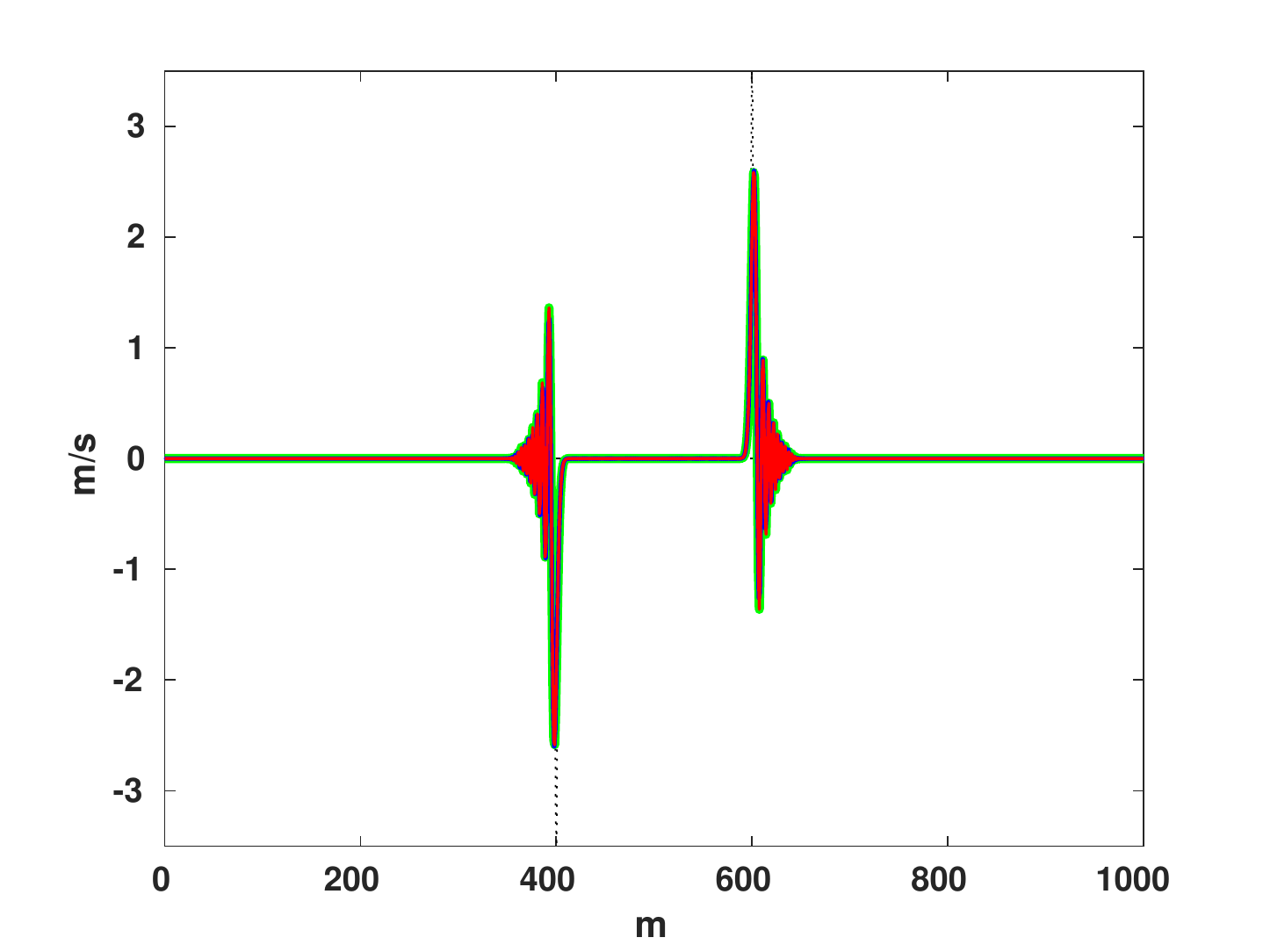}}
   \end{overpic}
   \end{tabular}
   \caption{Illustration of the effects of the discrete dispersion relation~\eqref{equ_disprel_split01} 
    of the stable P1--P0 and $\GPIu$--$\,\GPOh$ schemes on Gaussian-shaped height and velocity peaks 
    after a simulation time of $t=0.1T$.}
  \label{fig_disprel_spltv1h0}
  
  \end{figure}

\begin{figure}[h] \centering
  \begin{tabular}{ccc} 
   \begin{overpic}[scale=0.5,unit=1mm]
    {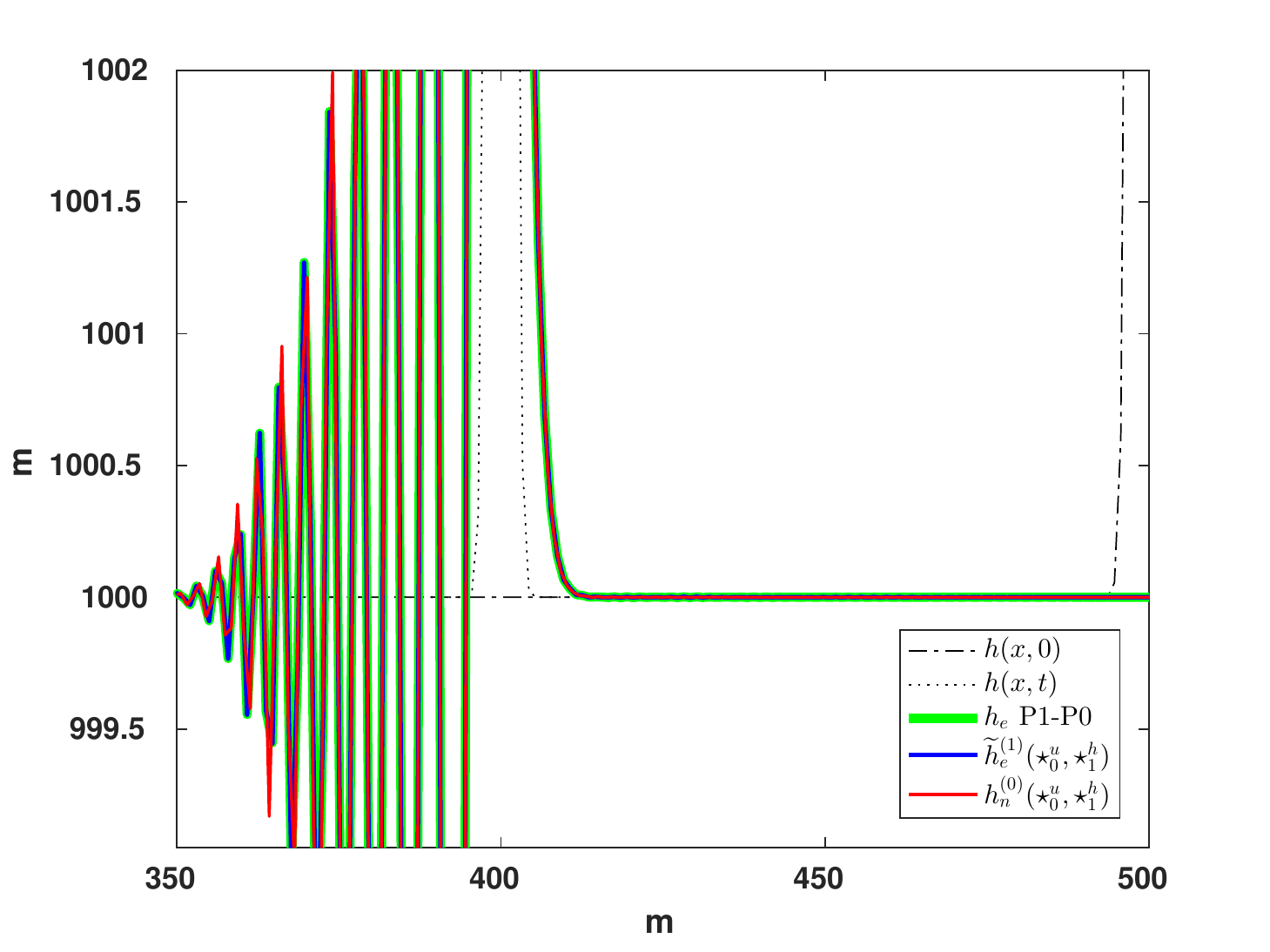}
     \put(33,26.2){\includegraphics[scale=.27]{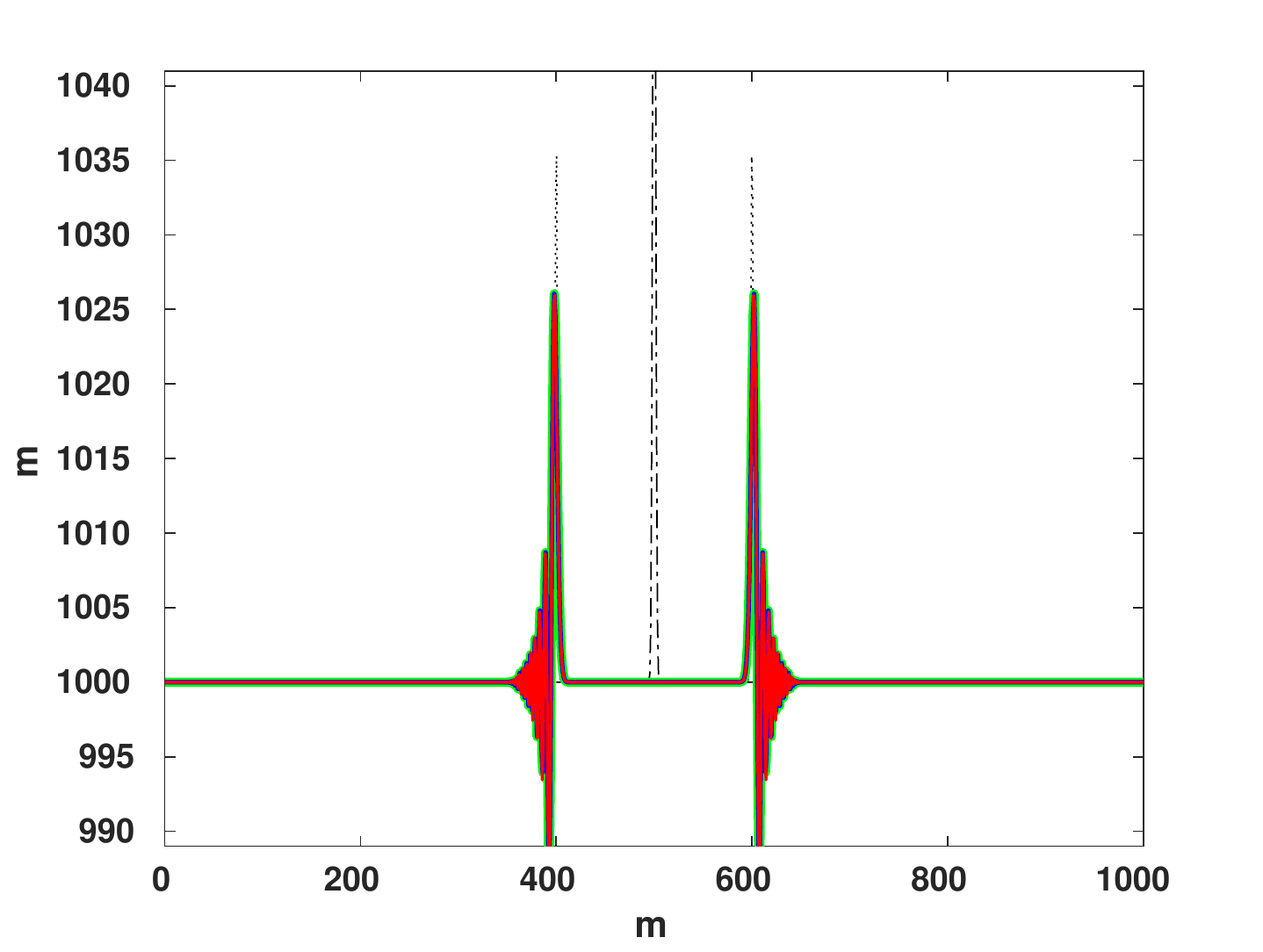}}
   \end{overpic}
   \begin{overpic}[scale=0.5,unit=1mm]
    {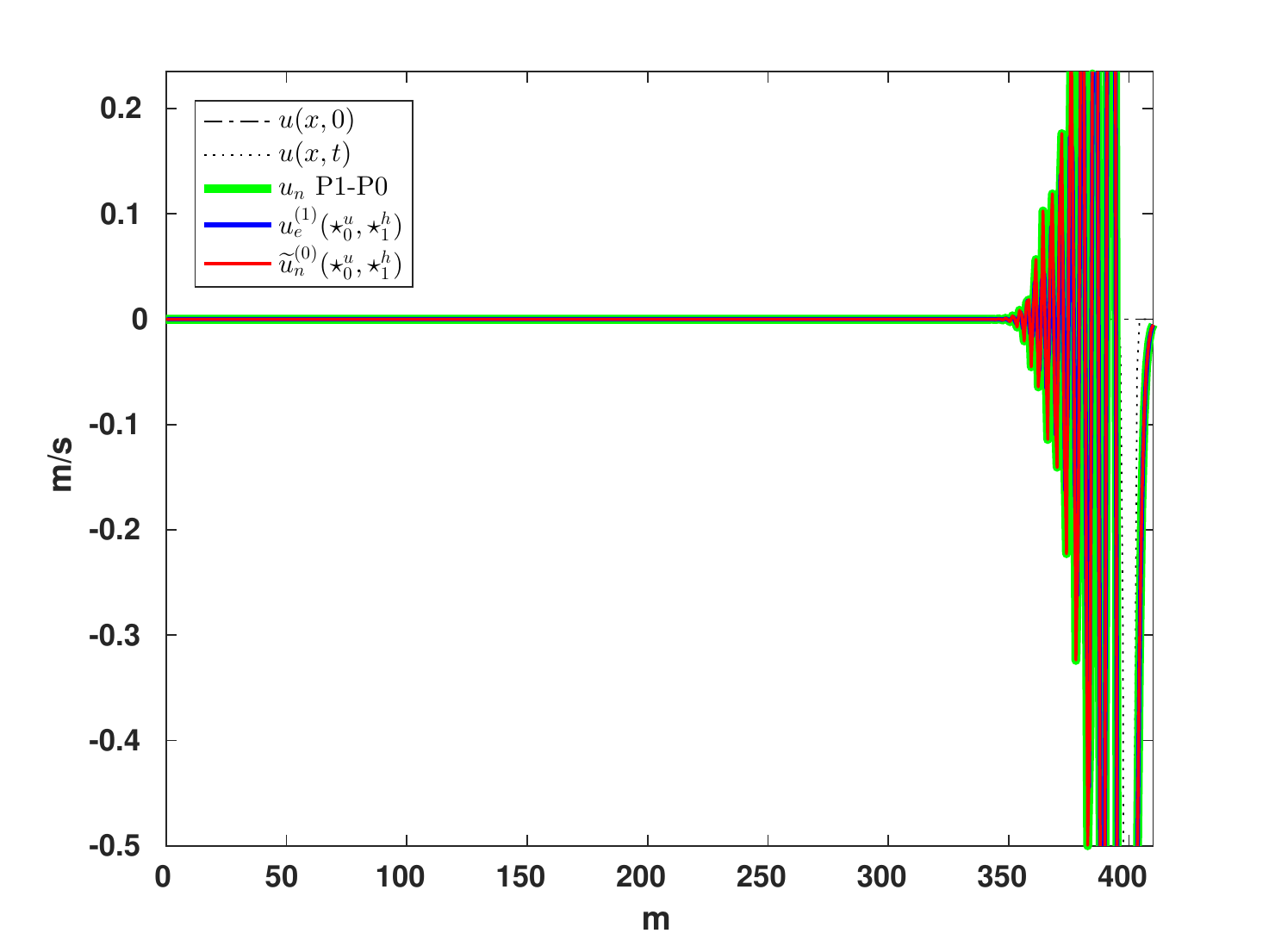}
      \put(14.2,6.6){\includegraphics[scale=.27]{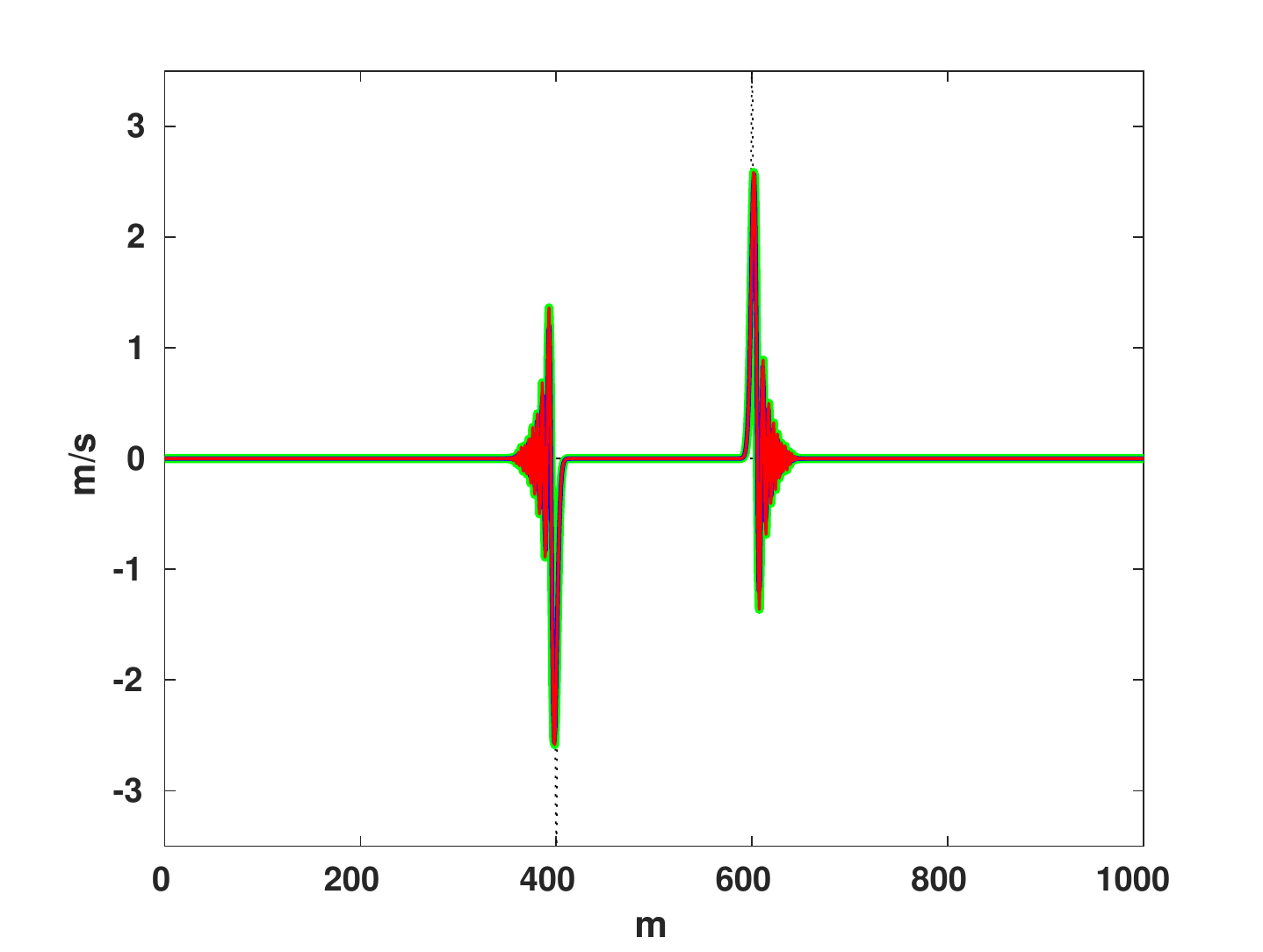}}
   \end{overpic}
   \end{tabular}
   \caption{Illustration of the effects of the discrete dispersion relation~\eqref{equ_disprel_split01} 
    of the stable P1--P0 and $\GPOu$--$\,\GPIh$ schemes on Gaussian-shaped height and velocity peaks 
    after a simulation time of $t=0.1T$.}
  \label{fig_disprel_spltv0h1}
  
  \end{figure}

  Let us finally study simulations of TC 3 performed with the mixed P1--P0, the split 
  $\GPIu$--$\,\GPOh$, and the split $\GPOu$--$\,\GPIh$ schemes.  First of all 
  we note that the solutions of the P1--P0 and of both split schemes agree very well. 
  This can be inferred from Figures~\ref{fig_disprel_spltv1h0} and \ref{fig_disprel_spltv0h1}, 
  in which the solutions of the P1--P0 scheme (green lines) 
  and those of the split schemes (blue and red lines) show only tiny differences.
  Again, propagation directions are correct, while the peaks become smaller and 
  wider. In between the peaks, no high frequency waves occupy the central region, similarly to the split 
  $\GPOu$--$\,\GPOh$ scheme. In contrast to the latter however, the high frequency waves 
  do not occupy the entire outer region but only some parts directly in front of the traveling 
  wave packages. 
  This reflects very well the theoretical findings described by the discrete dispersion 
  relation~\eqref{equ_disprel_split01}, which states that all waves travel equal to, or faster then 
  $c$, but not faster than about $1.1c$ (cf. Fig.~\ref{fig_disp_relation}).

  \section{Summary, Conclusions, Outlook}
  \label{sec_summary}

  We developed a new finite element (FE) discretization framework based on the 
  split 1D wave equations~\eqref{equ_split_wave} with the long term objective to 
  apply this framework to derive structure-preserving discretizations 
  for Geophysical Fluid Dynamics (GFD) \cite{Bauer2016}. 
  The splitting of these covariant equations into topological momentum/continuity equations 
  and metric-dependent closure equations relies on the introduction of additional differential 
  forms (DF). These establish pairs with the original DFs, which are connected by the Hodge-star operator. 
  By providing proper FE spaces for all DFs such that the differential operators 
  (here gradient and divergence) hold in strong form, our split FE framework provides 
  discretizations that preserve geometrical structure. In particular, the discrete momentum 
  and continuity equations are metric-free and hold in strong form (no partial integration is required). 
  They are therefore exact up to the errors introduced by the trivial projections into the correspond FE 
  spaces. Additional errors are however introduced by the discrete metric-dependent closure equations   
  that provide approximations to the Hodge-star operator by nontrivial Galerkin projections (GP).

  Requirements from structure preservation together with our objective to develop a 
  lowest-order scheme for the split 1D wave equations impose certain conditions on 
  the choice of FE spaces that approximate the corresponding DFs. 
  In detail, both topological equations are projected onto piecewise constant spaces, 
  whereas each of the two metric equations is projected onto both piecewise constant 
  {\em and} piecewise linear spaces. Although this gives four possible realizations, we find in 
  principle three classes of schemes: a {\em high accuracy scheme} ($\GPI{}$--$\,\GPI{}$), 
  in which the Galerkin projections map onto piecewise linear (GP1) spaces, 
  a {\em low accuracy scheme} ($\GPO{}$--$\,\GPO{}$), in which the Galerkin projections 
  map onto piecewise constant (GP0) spaces, and two {\em medium accuracy schemes} 
  ($\GPI{}$--$\,\GPO{}$), in which one Galerkin projection maps onto 
  piecewise linear, the other to piecewise constant space.
  
  A comparison to conventional methods, the unstable P1--P1 scheme that applies piecewise linear 
  (P1) spaces to approximate both velocity and height fields, and the stable P1--P0 scheme that 
  approximates the velocity by a piecewise linear and the height by a piecewise
  constant (P0) field, was performed with respect to conservation of mass and momentum, dispersion relations, 
  and convergence behavior. Analytical derivations showed that the discrete dispersion relations 
  of the unstable P1--P1 and $\GPI{}$--$\,\GPI{}$ schemes as well as of the stable P1--P0 
  and both $\GPI{}$--$\,\GPO{}$ schemes coincide, in spite of the differences in how they are derived.
  In particular, the $\GPI{}$--$\,\GPI{}$ scheme and the P1--P1 scheme share the problem 
  of permitting spurious modes, whereas the stable P1--P0 scheme and the two $\GPI{}$--$\,\GPO{}$
  schemes are free of such modes. 
      
  The piecewise linear fields of P1--P1 and $\GPI{}$--$\,\GPI{}$ schemes show more or less 
  the same convergence rates of second order. Applying four instead 
  of only two FE spaces, we illustrated that the piecewise 
  constant fields, present only in the split schemes, showed the expected first order 
  convergence rate. Similarly, the corresponding fields of the P1--P0 and both $\GPI{}$--$\,\GPO{}$
  schemes agree in their convergence behavior. Both piecewise constant height fields 
  show first order and both piecewise linear velocity fields second order convergence 
  rates. Utilizing a piecewise linear approximation to the height field, the 
  split $\GPI{}$--$\,\GPO{}$ schemes give for both velocity {\em and} height fields
  second order convergence rates. This combination of projections can 
  therefore be interpreted as stable version of the conventional unstable P1--P1 scheme. 
    
  We further verified the structure-preserving nature of the split schemes fulfilling 
  the first principles of mass and momentum conservation discretely. In more detail, 
  the corresponding relative errors show fluctuations on the order of about $10^{-9}$ 
  without indicating a trend. When compared to the conventional schemes, 
  this conservation behavior is only shared by the stable P1--P0 scheme, 
  whereas the unstable P1--P1 scheme shows a tendency to lose mass, for the given 
  simulation time on the order of $10^{-6}$.

  Besides the $\GPI{}$--$\,\GPI{}$ and $\GPI{}$--$\,\GPO{}$ schemes that have a counterpart 
  with the P1--P1 and P1--P0 schemes in literature, we found a split $\GPO{}$--$\,\GPO{}$
  scheme, for which we are not aware of a corresponding standard formulation. 
  This low accuracy scheme preserves the first principles of mass and momentum 
  conservation in the same way as the other split schemes. However, it does not 
  provide a very accurate approximation to the dispersion relation, exposing
  very fast traveling short waves. Despite of the latter drawback, the finding 
  of this new scheme illustrates the general character of our split FE discretization 
  framework. By allowing a larger variety of combinations of different FE spaces 
  compared to conventional mixed FE methods, this framework provides a toolbox 
  to find and study new combinations of finite element spaces. In this vein, 
  it is subject to current and future work to investigate various low- and higher-order 
  FE space combinations to find new structure-preserving discretizations of the 
  nonlinear rotating shallow-water equations, and of the equations of GFD in general.

  \section{Acknowledgements}

  This project has received funding from the European Union's Horizon 2020
  research and innovation programme under the Marie Sk\l odowska-Curie grant agreement No 657016. 
  Additionally, the second author acknowledges support by the DFG excellence cluster CliSAP (EXC177).


 
 \newcommand{\fctVOtriall}[1]{\color{black}  \phi_{{#1}}(x)     \color{black}}
 \newcommand{\fctVItriall}[1]{\color{black}  \chi_{{#1}}(x)     \color{black}}

  \appendix

  \section{Explicit representation of mass and stiffness matrices}
  \label{append_int_values}

  Here, we explicitly represent the mass and stiffness matrices used
  in Sect.~\ref{sec_standard} and Sect.~\ref{sec_split}. In particular, 
  we show the coefficients which follow by evaluating the integrals of the variational 
  formulations. Imposing periodic boundary conditions, the presented matrices 
  correspond to a mesh with $N$ independent 
  DoFs for nodes and elements, as illustrated in Figure~\ref{fig_mesh}. 
  
  The basis functions for the piecewise linear FE space are given by the 
  Lagrange functions $\fctVOtrial, \ \forall l = 1,\dots N$, and those for the piecewise constant FE space
  by the step functions $\fctVItrial, \ \forall m=1,\dots N,$ (cf. Figure~\ref{fig_mesh}).
  They are defined by 
  \begin{equation}
   \fctVOtrial = \begin{cases}
                  1 - \frac{x - x_l     }{\Delta x_m       } &    \text{if} \ x \in [x_l    , x_{l+1}] ,\\
                     \frac{ x - x_{l-1} }{\Delta x_{m - 1} } &    \text{if} \ x \in [x_{l-1}, x_{l}  ] ,\\
                  0 			         &  \text{else};
                 \end{cases}
   \qquad 
   \fctVItrial = \begin{cases}
                  1    &  \text{if} \  x \in [x_l    , x_{l+1}] ,\\
                  0    &  \text{else}.
                 \end{cases}
  \end{equation}  
   
  Note that, for all $\indxVOtrial = 1,\dots N$, there is
  $\int_L \fctVOtriall{l-1} \fctVOtrial dx = \frac{1}{6} \Delta x_{\indxVItrial-1}$,
  $\int_L \fctVOtriall{l  } \fctVOtrial dx = \frac{1}{3} \Delta x_{\indxVItrial-1}+\frac{1}{3} \Delta x_{\indxVItrial} $,  
  $\int_L \fctVOtriall{l+1} \fctVOtrial dx = \frac{1}{6} \Delta x_{\indxVItrial}$,
  and zero else. Hence, the evaluation of coefficients~\eqref{equ_mssmxnn_coeff} 
  gives the metric-dependent mass-matrix
  \begin{small}
  \begin{equation}
   \massmxnn =  \left( \begin{matrix} 
         \frac{1}{3}\Delta x_{N} + \frac{1}{3}\Delta x_{1} &  \frac{1}{6}\Delta x_{1}                           &  0                           & \ldots &   \frac{1}{6}\Delta x_{N}   \\ 
         \frac{1}{6}\Delta x_{1}                           &  \frac{1}{3}\Delta x_{1} + \frac{1}{3}\Delta x_{2} &  \frac{1}{6}\Delta x_{2}     & \ldots   &      0                      \\ 
         \vdots                                            &  \ddots                                            &  \ddots                      &  \ddots &\vdots\\ 
         0                                                 &  \ldots  & \frac{1}{6}\Delta x_{N-2} &  \frac{1}{3}\Delta x_{N-2} + \frac{1}{3}\Delta x_{N-1} &  \frac{1}{6}\Delta x_{N}    \\
         \frac{1}{6}\Delta x_{1} &   \ldots & 0 & \frac{1}{6}\Delta x_{N-1} & \frac{1}{3}\Delta x_{N-1} + \frac{1}{3}\Delta x_{N}     
        \end{matrix} \right) .
  \end{equation}  
  \end{small}
  
  \noindent
  Because of the periodic boundary, we identify here and throughout the entire manuscript 
  $\fctVOtriall{0}$ with $\fctVOtriall{N}$ and $\fctVOtriall{N+1}$ with $\fctVOtriall{1}$.
  Moreover, for all $\indxVOtrial = 1,\dots N$, there is
   $\int_L \frac{d\fctVOtriall{l \mp 1}}{dx} \fctVOtrial dx = \mp \frac{1}{2}$, and zero else. 
  Therefore, the evaluation of coefficients~\eqref{equ_stiffxnn_coeff} gives the metric-free stiffness matrix
  \begin{equation}
    \stiffmxnn =  \left( \begin{matrix} 
         0              &  \frac{1}{2}   &   0           &   \ldots      &  \frac{1}{2}  \\ 
         -\frac{1}{2}   & 	0        &  \frac{1}{2}  &   \dots       &   0           \\ 
         \vdots         &  \ddots        &  \ddots       &   \ddots      &   \vdots       \\ 
         0              & \ldots         &  -\frac{1}{2} &   0           &  \frac{1}{2}     \\   
         \frac{1}{2}    & \ldots         &   0           & -\frac{1}{2}  &  0      
        \end{matrix} \right)  .
  \end{equation}  
  
  Applying the values $\int_L \fctVItrial  \fctVItrial = \Delta x_\indxVItrial$ for all $m = 1,\dots N$, zero else,
  to evaluate coefficients~\eqref{equ_mssmxee_coeff} and 
  $\int_L \frac{d \fctVOtriall{l/l+1}}{dx} \fctVItrial dx = -1 / 1$ for all $m = 1,\dots N$, zero else, 
  to calculate coefficients~\eqref{equ_stiffxne_coeff}, 
  there follow the metric-dependent mass-matrix $\massmxee$ and the metric-free stiffness matrix $\stiffmxen$: 
  \begin{equation}\notag
    \massmxee = \left( \begin{matrix} 
          \Delta x_{1} &  0           &  \ldots    & 0        \\ 
          0            & \Delta x_{2} &  \ldots    & 0        \\ 
          \vdots       & \vdots       &  \ddots    & \vdots   \\ 
          0            & 0            &  \ldots    & \Delta x_{N}       
        \end{matrix} \right)  , 
        \qquad         
    \stiffmxen = \left( \begin{matrix} 
         -1     	&   1  		&   0           &  \ldots    &  0       \\ 
          0     	&  -1  		&   1           &  \ldots    & 0        \\ 
          \vdots       &   \vdots       &  \ddots    	& \ddots     & \vdots   \\ 
          0            &    0           &  \ldots  	& -1  	    & 1         \\
          1            &    0           &  \ldots  	& 0  	    & -1          
        \end{matrix} \right)  , 
  \end{equation}  
  with $\stiffmxne = (\stiffmxen)^T$.

  Finally, using the values $\int_L \fctVItriall{m-1/m} \fctVOtrial dx = \frac{1}{2}\Delta x_{m-1}/\frac{1}{2}\Delta x_{m}$ 
  for all $l = 1,\dots N$, zero else, the evaluation of coefficients~\eqref{equ_mssmxne_coeff} 
  gives the metric-dependent mass-matrix $\massmxne$ with corresponding metric-free projection matrix $\prjcmxne$:
  \begin{equation}\notag
    \massmxne = \left( \begin{matrix} 
          \frac{1}{2}\Delta x_1	&                      0           &  \ldots    &  \frac{1}{2}\Delta x_N     \vspace{0.3em}  \\
          \frac{1}{2}\Delta x_1	&  \frac{1}{2}\Delta x_2   		    &  \ldots    &  0      \\ 
          \vdots            &  \ddots    	& \ddots     & \vdots   \\ 
          0            &    \ldots           	&   \frac{1}{2}\Delta x_{N-1} 	    &  \frac{1}{2}\Delta x_{N}         
        \end{matrix} \right) , 
   \qquad 
    \prjcmxne = \left( \begin{matrix} 
          \frac{1}{2} 	&                      0           &  \ldots    &  \frac{1}{2}     \vspace{0.3em}  \\
          \frac{1}{2} 	&  \frac{1}{2}   		    &  \ldots    &  0      \\ 
          \vdots            &  \ddots    	& \ddots     & \vdots   \\ 
          0            &    \ldots           	&   \frac{1}{2}  	    &  \frac{1}{2}          
        \end{matrix} \right) .   
  \end{equation}

  \section{Derivation of the discrete dispersion relations}
  \label{app_disc_disp_rel}
  
  In this section we present a detailed derivation of the discrete dispersion relations, 
  presented in Sect.~\ref{sec_standard} for the mixed, and in Sect.~\ref{sec_split} 
  for the split schemes. For the mixed schemes we mainly follow the approach presented 
  in \cite{WaltersCarey1983}. For the split schemes, we adapt this approach such that it 
  suits to the split form of the wave equations.

  \subsection{The analytical case}
  
  To explain the method, we first determine the dispersion relation for the 
  analytical wave equations.
  As above, we consider the domain $L$ with coordinates 
  $x \in [0,L] \subset \mathbbm R$ and time variable $t \in \mathbbm R$. 
  We assume that the solutions can be formulated as periodic in space and time.
  In more detail, we seek periodic solutions of the equations~\eqref{equ_standard_wave}, 
  i.e. we decompose the solutions like
  \begin{equation}\label{equ_sol_rep}
   u(x,t) =   \hat{u}(x) e^{i\omega t }  =  \overline{u} e^{i\omega t - i k x}, \quad  h(x,t) =   \hat{h}(x) e^{i\omega t }  =  \overline{h} e^{i\omega t - i k x} , 
  \end{equation}
  in which $\overline{u},\overline{h}$ are constants. $\omega$ is the angular frequency
  and $k$ is the wave number.  
  
  Applying the first decomposition of \eqref{equ_sol_rep} with respect to time to evaluate 
  the time derivative of \eqref{equ_standard_wave}, 
  we obtain equations in {\em harmonic form}
  \begin{equation}\label{equ_harmform}
   i\omega \hat{u}(x) + g \frac{\partial \hat{h}(x)}{\partial x} = 0 , \quad i\omega  \hat{h}(x) + H \frac{\partial \hat{u}(x)}{\partial x} = 0 \, .
  \end{equation} 
  Using further the second decomposition of \eqref{equ_sol_rep} in wave vector space with wave vector $k$, 
  we find  
  \begin{equation}\label{equ_standard_wave_matrixform}
  \left( \begin{array}{cc} 
	i\omega & igk  \\ 
	iHk & i\omega
	  \end{array} \right) 
  \left( \begin{array}{c} 
	    \overline{u} \\ \overline{h} 
	  \end{array} \right) = 0 \, .
  \end{equation}
  For nontrivial solutions to exist, the determinant of the coefficient matrix in 
  \eqref{equ_standard_wave_matrixform} must vanish. 
  Then, the analytic dispersion relation for shallow-water waves reads
  \begin{equation}\label{equ_dsprel_an}
    c = \frac{\omega}{k} = \pm \sqrt{gH} \, ,
  \end{equation}
  in which $c$ describes the phase velocity. The positive and negative solutions stand for right and 
  left traveling waves, respectively.
  
  The preceding equation also describes the analytical dispersion relation of the split wave 
  equations~\eqref{equ_split_wave}. This fact follows directly by substituting in \eqref{equ_split_wave} 
  the metric into the topological equations, which leads to equations~\eqref{equ_standard_wave}. 
  An explicit treatment of the split equations avoiding such substitution, as presented
  in Sect.~\ref{sec_der_dspsplit} for the discrete case, would result in the dispersion
  relation~\eqref{equ_dsprel_an} too.

  \subsection{The discrete case}

  The discrete dispersion relations of the schemes derived in this manuscript 
  can be calculated analogously to the continuous case. In particular,
  we begin our derivations from equations in harmonic form (as in \eqref{equ_harmform}), 
  in which the equations are continuous in time and only discretized in space. 
  Representing the variables in terms of Fourier expansions on nodal and elemental unknowns,
  this approach results in discrete dispersion relations in which the magnitude of the 
  temporal frequency is expressed as a function of the wave number $k$, cf. \cite{LeRoux2007,LeRoux2008}.
  
  To simplify calculations, we assume for all the following cases that the domain with length $L$ is discretized by
  a uniform mesh whose $N$ elements have equilateral size $\Delta x = \Delta x_m = x_{l+1} - x_l = \frac{L}{N}$ 
  for all $m,l =1,\dots N$, (cf. Fig.~\ref{fig_mesh}).

  \subsubsection{The unstable P1--P1 FE scheme}

  The semi-discrete harmonic equations of the unstable P1--P1 scheme 
  are given by \eqref{equ_stnd_1} in Sect.~\ref{sec_standard00}, in which the partial time derivatives 
  $\partial_t \uVOtrial(t)$ and $\partial_t \hVOtrial(t)$ are replaced by 
  $i\omega \uVOtrial$ and $i\omega \hVOtrial$, respectively. 
  Expanding further all nodal values $\uVOtrial,\hVOtrial $ in wave vector form, i.e. 
  $\uVOtrial = \overline{u}e^{ikx_l}, \hVOtrial = \overline{h}e^{ikx_l} $ for constants $\overline{u}$ and $\overline{h}$, 
  we obtain 
  \begin{equation}\label{equ_disp_P1P1_1}
  \begin{split}
      & i\omega 
        \left\{ \begin{array}{c} \overline{u} \\ \overline{h} \end{array} \right\}
     \Big(e^{ikx_{l-1}} \int_L \fctVOtriall{l-1} \fctVOtrial  dx
        + e^{ikx_{l}}   \int_L \fctVOtriall{l}   \fctVOtrial  dx
        + e^{ikx_{l+1}} \int_L \fctVOtriall{l+1} \fctVOtrial  dx \Big) \\	
       & + \left\{ \begin{array}{c} g \\ H \end{array} \right\}
           \left\{ \begin{array}{c} \overline{h} \\ \overline{u} \end{array} \right\}
     \Big(e^{ikx_{l-1}} \int_L \frac{d \fctVOtriall{l-1}  }{dx} \fctVOtrial dx 
         +e^{ikx_{l+1}} \int_L \frac{d \fctVOtriall{l+1}  }{dx} \fctVOtrial dx \Big) =     0  \ \forall \indxVOtrial = 1,\dots N .
  \end{split}
  \end{equation}
  Dividing the preceding equations by $e^{ikx_{l}}$ and noting that 
  $\Delta x:= x_{l+1} -  x_{l} = x_{l} -  x_{l-1}$, we find
       \begin{equation}\label{equ_disp_P1P1_2}
	\begin{split}
       & i\omega 
         \left\{ \begin{array}{c} \overline{u} \\ \overline{h} \end{array} \right\}
	\Delta x \Big(  \frac{2}{3} 
	      + \frac{1}{6} e^{-ik\Delta x}   
	      + \frac{1}{6} e^{ ik\Delta x}   \Big) 	
        + \left\{ \begin{array}{c} g \\ H \end{array} \right\}
	 \left\{ \begin{array}{c} \overline{h} \\ \overline{u} \end{array} \right\}
	\Big(- \frac{1}{2} e^{-ik\Delta x} 
	     + \frac{1}{2} e^{ ik\Delta x} \Big) =     0 \, ,  \\
      & i\omega 
         \left\{ \begin{array}{c} \overline{u} \\ \overline{h} \end{array} \right\}
	\Delta x\Big(  \frac{2}{3}  + \frac{1}{3} \cos (k\Delta x)    \Big) 	
        + \left\{ \begin{array}{c} g \\ H \end{array} \right\}
	 \left\{ \begin{array}{c} \overline{h} \\ \overline{u} \end{array} \right\}
	\Big(i \sin (k\Delta x) \Big) =     0 , 
	\end{split}
      \end{equation}
    where we used the integral values of Appendix~\ref{append_int_values}. 
    Equivalently, these equations can be written in matrix-vector form:
    \begin{equation}\label{equ_disp_P1P1_Matrix}
       \left( \begin{array}{cc} 
         \frac {i\omega }{3} \Big(  2  +   \cos (k\Delta x)    \Big)    &  \frac{ig}{\Delta x}  \sin (k\Delta x)  \\ 
         \frac{iH}{\Delta x}  \sin (k\Delta x)     & \frac{ i\omega}{3} \Big(  2  +   \cos (k\Delta x)    \Big)
              \end{array} \right) 
        \left( \begin{array}{c} 
           \overline{u} \\ \overline{h} 
               \end{array} \right) = 0 . 
    \end{equation}
    
    Setting the determinant of this coefficient matrix to zero gives the discrete dispersion relation
    \begin{equation}\label{equ_disp_rel_P1P1}
       c_d = \frac{\omega}{k} = \pm \sqrt{gH}\frac{\sin (k\Delta x) }{k \Delta x} \frac{3 }{ 2 + \cos (k\Delta x)  } \, .
    \end{equation}
    For $k \rightarrow 0$, there follows $c_d \rightarrow c = \sqrt{gH}$, which can be verified by expanding 
    $\sin (k {\Delta x})$ in a Taylor series around zero. At the shortest wave length $k = \frac{\pi}{\Delta x}$,
    equation~\eqref{equ_disp_rel_P1P1} has a second root (spurious mode). Therefore, the P1--P1 scheme is unstable.

    \subsubsection{The stable P1--P0 FE scheme}

    In case of the stable P1--P0 scheme introduced in Sect.~\ref{sec_standard01},
    the semi-discrete harmonic equations are given by \eqref{equ_P0P1},
    with $\partial_t \uVOtrial(t)$ and $\partial_t \hVItrial(t)$ replaced by 
    $i\omega \uVOtrial$ and $i\omega \hVItrial$, respectively. 
    Here we expand the variables for all nodes as $\uVOtrial = \overline{u}e^{ikx_l}$ 
    and for all elements as $\hVItrial = \overline{h}e^{ikx_{m}} $.
    
    Then, the discrete momentum equation reads
    \begin{equation}\label{equ_disp_P0P1_1}
	\begin{split}
       & i\omega 
        \overline{u} 
	\Big( e^{ikx_{l-1}} \int_L  \fctVOtriall{l-1}  \fctVOtrial    dx
	    + e^{ikx_{l  }} \int_L  \fctVOtrial        \fctVOtrial    dx
	    + e^{ikx_{l+1}} \int_L  \fctVOtriall{l+1}  \fctVOtrial    dx \Big) \\	
	& - g \overline{h} 
	\Big( e^{ikx_{m-1}} \int_L  \fctVItriall{m-1}  \frac{d\fctVOtrial }{dx} dx 
	    + e^{ikx_{m  }} \int_L  \fctVItriall{m}    \frac{d\fctVOtrial }{dx} dx \Big) =     0  , \  \forall \ l = 1,\dots N,
	\end{split}
    \end{equation}
    and the continuity equation reads 
    \begin{equation}
      \begin{split}
        &   i\omega \overline{h}   e^{ikx_{m  }}      \int_L \fctVItrial  \fctVItrial  dx + \\
	 &  H  \overline {u}\Big(  e^{ikx_{l  }}  \int_L \frac{d \fctVOtriall{l  } }{dx}  \fctVItrial  dx 
	                       +   e^{ikx_{l+1}}  \int_L \frac{d \fctVOtriall{l+1} }{dx}  \fctVItrial  dx     \Big) =     0 , \  \forall \ m = 1,\dots N.
       \end{split}
    \end{equation}
    Dividing the former equation by $e^{ikx_{l}}$ and the latter by $e^{ikx_{m}}$,
    noting that $\frac{\Delta x}{2} = x_{l+1} -  x_{m} = x_{m} -  x_{l-1}$ on a uniform mesh, 
    and using the integral values of Appendix~\ref{append_int_values}, 
    the momentum equation becomes
    \begin{equation}\label{equ_disp_P0P1_2}
	\begin{split}
        i\omega    \overline{u}  
        \Delta x  \Big(  \frac{2}{3} 
	      + \frac{1}{6} e^{-ik\Delta x}   
	      + \frac{1}{6} e^{ ik\Delta x}   \Big) 	
        +   g  \overline{h} \Big( -  e^{-ik \frac{\Delta x}{ 2} }   +  e^{ik\frac{\Delta x}{ 2} } \Big) &  =     0   \, , \\
       i\omega   \overline{u} \Delta x \Big(  \frac{2}{3}  + \frac{1}{3} \cos (k\Delta x)    \Big) 	
        +   g \overline{h} 2 i  \sin \left(k \frac{\Delta x}{ 2}  \right)  &  =     0 \, ,
	\end{split}
    \end{equation}
    and the continuity equation becomes
    \begin{equation}\label{equ_disp_P0P1_3}
	\begin{split}
        i\omega    \overline{h}  
        \Delta x  + H   \overline{u} \Big( -  e^{-ik \frac{\Delta x}{ 2} }   +  e^{ik\frac{\Delta x}{ 2} } \Big) & =     0   \, , \\
       i\omega    \overline{h}  
        \Delta x  + H   \overline{u} 2 i  \sin \left(k \frac{\Delta x}{ 2}  \right) &  =     0   \, . \\
      	\end{split}
      \end{equation}
    In matrix-vector formulation, the preceding set of discrete equations yields 
    \begin{equation}\label{equ_disp_P0P1_Matrix}
	    \left( \begin{array}{cc} 
		  \frac {i\omega }{3} \left(  2  +   \cos (k\Delta x) \right)    &  \frac{i2g}{\Delta x}  \sin (k\frac{\Delta x}{ 2}) \vspace{0.3em} \\ 
		  \frac{i2H}{\Delta x}  \sin (k\frac{\Delta x}{ 2})     & i \omega 
		    \end{array} \right) 
	    \left( \begin{array}{c} 
		  \overline{u}\vspace{0.3em}  \\ \overline{h} 
		\end{array} \right) = 0 \, .
    \end{equation}
    
    Setting the determinant of this coefficient matrix to zero results in the discrete dispersion relation
    \begin{equation}\label{equ_disp_rel_P0P1}
       c_d =  \frac{\omega}{k} = \pm \sqrt{gH}  \frac{\sin (k \frac{\Delta x}{ 2} ) } { k \frac{\Delta x}{ 2}} 
                                            \left[  \frac{ 3}{ 2 + \cos (k\Delta x)  } \right]^{1/2} \, .
    \end{equation}
    Also here, $c_d \rightarrow c$ for $k \rightarrow 0$. However, in contrast to the unstable P1--P1 
    scheme, there is no second root at shortest wave length $k = \frac{\pi}{\Delta x}$. Instead, the 
    discrete dispersion relation approximates well the analytical one. 
    Therefore, we consider this scheme as stable.

  \subsubsection{The split FE schemes}
  \label{sec_der_dspsplit}

  For the split schemes derived in Sect.~\ref{sec_split}, the semi-discrete harmonic
  topological equations are given by \eqref{equ_splitFE_top},  
  with $\partial_t \uVItrial(t)$ and $\partial_t \htwVItrial(t)$ replaced by 
  $i\omega \uVItrial$ and $i\omega \htwVItrial$, respectively.
  Here we expand the coefficient of velocity and height in wave vector form as follows: at nodal values as 
  $ \utwVOtrial = \overline{\widetilde u}e^{ikx_l},  \hVOtrial = \overline{h}e^{ikx_l} $, 
  and at elemental values as $\htwVItrial   = \overline{\widetilde h}e^{ikx_{m}},  \uVItrial  = \overline{u}e^{ikx_{m} } $,
  with constants $\overline{\widetilde u},\overline{h},\overline{\widetilde h},$ and $\overline{u}$.
  We obtain
  \begin{equation}\label{equ_disp_splitP0P1_1}
  \begin{split}
       & i\omega 
         \left\{ \begin{array}{c} \overline{u} \\ \overline{\widetilde h} \end{array} \right\}
          \Big(e^{ikx_{m}} \int_L \fctVItrial \fctVItrial    dx  \Big) \\	
       & + \left\{ \begin{array}{c} g \\ H \end{array} \right\}
          \left\{ \begin{array}{c} \overline{h} \\ \overline{\widetilde u} \end{array} \right\}
          \Big(e^{ikx_{l  }} \int_L \frac{d \fctVOtriall{l  }  }{dx} \fctVItrial  dx 
              +e^{ikx_{l+1}} \int_L \frac{d \fctVOtriall{l+1}  }{dx} \fctVItrial  dx \Big) =     0 , \ \forall \indxVItrial =1,\dots N. 
  \end{split}
  \end{equation}
  Dividing the preceding equations by $e^{ikx_{m}}$ and using 
  $\frac{\Delta x}{2} = x_{l+1} -  x_{m} = x_{m} -  x_{l-1}$ and the 
  integral values of Appendix~\ref{append_int_values}, there follows
  \begin{equation}\label{equ_disp_splitP0P1_2}
  \begin{split}
	 i\omega 
         \left\{ \begin{array}{c} \overline{u} \\ \overline{\widetilde h} \end{array} \right\}
	 \Delta x   
	+ \left\{ \begin{array}{c} g \\ H \end{array} \right\}
	 \left\{ \begin{array}{c} \overline{h} \\ \overline{\widetilde u} \end{array} \right\}
	\Big(e^{ik \frac{\Delta x}{2} }  -  e^{-ik \frac{\Delta x}{2} }  \Big) &  =     0  \, , \\
	 i\omega 
         \left\{ \begin{array}{c} \overline{u} \\ \overline{\widetilde h} \end{array} \right\}
	+ \left\{ \begin{array}{c} g \\ H \end{array} \right\}
	 \left\{ \begin{array}{c} \overline{h} \\ \overline{\widetilde u} \end{array} \right\}
	\frac{2i}{\Delta x} \sin \left( k \frac{\Delta x}{2} \right)           & =     0 \, . \\
  \end{split}
  \end{equation}
  
  In Sect.~\ref{sec_disc_metr}, we presented four realizations of discrete metric equations 
  as closure conditions to the topological equations~\eqref{equ_matrix_splitFE}.
  Here we proceed similarly and derive discrete metric equations that close the 
  harmonic equations~\eqref{equ_disp_splitP0P1_2},
  i.e. we distinguish between high, medium, and low accuracy closures.

  \paragraph{High accuracy closure.}

  The discrete metric closure equations for the high accuracy case are given by 
  \eqref{equ_splitFE_m}. Representing the variables in wave vector form at 
  nodal and elemental values using the expansion as for the topological equations,
  we obtain 
  \begin{equation}\label{equ_split_MetP1P1_discrete}
  \begin{split}
	   &  \left\{ \begin{array}{c}  \overline{\widetilde u} \\ \overline{h }      \end{array} \right\} 
	    \Big(e^{ikx_{l-1}} \int_L \fctVOtriall{l-1} \fctVOtrial   dx
	       + e^{ikx_{l  }} \int_L \fctVOtriall{l  } \fctVOtrial   dx
	       + e^{ikx_{l+1}} \int_L \fctVOtriall{l+1} \fctVOtrial   dx \Big) \\	
	    & =    \left\{ \begin{array}{c} \overline{ u}          \\ \overline{\widetilde h }   \end{array} \right\}  
	    \Big(e^{ikx_{m-1}} \int_L \fctVItriall{m-1} \fctVOtrial   dx 
	        +e^{ikx_{m  }} \int_L \fctVItriall{m  } \fctVOtrial   dx \Big) \ \forall \indxVOtrial = 1,\dots N. 
  \end{split}
  \end{equation}
  Dividing the preceding equation by $e^{ikx_{l}}$, using $\frac{\Delta x}{2} = x_{l+1} -  x_{m} = x_{m} -  x_{l-1}$, 
  and applying the integral values of Appendix~\ref{append_int_values}, we find
  \begin{equation}\label{equ_split_MetP1P1_discrete2}
  \begin{split}
	   \left\{ \begin{array}{c}  \overline{\widetilde u} \\ \overline{h }      \end{array} \right\} 
	    \Delta x  \Big(  \frac{2}{3}  
		  + \frac{1}{6} e^{-ik\Delta x} 		    
		  + \frac{1}{6} e^{ ik\Delta x}   \Big) 
         & =  \left\{ \begin{array}{c} \overline{ u}          \\ \overline{\widetilde h }   \end{array} \right\}  
	     \frac{\Delta x}{2} \Big(  e^{-ik \frac{\Delta x}{ 2} }   +  e^{ik\frac{\Delta x}{ 2} } \Big) \, , \\
	   \left\{ \begin{array}{c}  \overline{\widetilde u} \\ \overline{h }      \end{array} \right\} 	  
	    \Big(  \frac{2}{3}  + \frac{1}{3} \cos (k\Delta x)    \Big) 	
	   & = \left\{ \begin{array}{c} \overline{ u}          \\ \overline{\widetilde h }   \end{array} \right\}  
	     \cos \left(k \frac{\Delta x}{ 2}  \right)  .
  \end{split}
  \end{equation}
  The combination of these latter with the discrete topological 
  equations~\eqref{equ_disp_splitP0P1_2} gives the 
  matrix-vector equations:
  \begin{equation}\label{equ_split_disp_P0P1_Matrix_v1}
	    \left( \begin{array}{cccc} 
	           i \omega & i \frac{2g}{\Delta x}  \sin (k \frac{\Delta x}{2})  & 0 & 0 \\ 
	           0 & 0 &  i \frac{2H}{\Delta x}  \sin (k \frac{\Delta x}{2}) & i \omega  \\ 
                  - \cos (k \frac{\Delta x}{2}) & 0 & \Big(  \frac{2}{3}  +   \frac{1}{3} \cos (k\Delta x)    \Big) & 0  \\  
                   0 & \Big(  \frac{2}{3}  +   \frac{1}{3} \cos (k\Delta x)    \Big) & 0 &  - \cos (k \frac{\Delta x}{2}) 		  
		    \end{array} \right) 
	    \left( \begin{array}{c} 
		  \overline{u} \vspace{0.4ex} \\ \overline{h} \vspace{0.4ex} \\ \overline{\widetilde u} \vspace{0.4ex}\\ \overline{\widetilde h} \\ 
		\end{array} \right) = 0 \, .
	  \end{equation}
	  
  To find nontrivial solutions of this equation, we set the determinant of this coefficient matrix to zero. 
  There follows the discrete dispersion relation
  \begin{equation}
           c_d = \frac{{\omega_{11}}}{k} = \pm \sqrt{gH} \frac{\sin (k\frac{\Delta x}{2})}{k \frac{\Delta x}{2}} 
                                      \frac{ 3\cos (k \frac{\Delta x}{2}) }{ (2 + \cos (k\Delta x))  } \, .
  \end{equation}
  For $k \rightarrow 0$, there follows $c_d \rightarrow c$. The high accuracy scheme $\GPIu$--$\,\GPIh$
  has a root solution at $k= \frac{\pi}{\Delta x}$ and is hence unstable.
  Because of the double angle formula (i.e. $\sin (k \Delta x) = 2\sin(k\frac{\Delta x}{2})\cos(k\frac{\Delta x}{2}  )$),
  this dispersion relation equals the one in \eqref{equ_disp_rel_P1P1} for the unstable P1--P1 scheme.

  \paragraph{Low accuracy closure.} 
  
  In the low accuracy case, the discrete metric equations are given by \eqref{equ_splitFE_m3},
  in which the variables at nodal and elemental values are expanded as for the topological equations.
  Then,	
  \begin{equation}\label{equ_split_MetP0P0_discrete}
	    \begin{split}
	   &  \left\{ \begin{array}{c}  \overline{\widetilde u} \\ \overline{h }      \end{array} \right\} 
	    \Big(  e^{ikx_{l  }}\int_L  \fctVOtriall{l  } \fctVItrial  dx 
		+  e^{ikx_{l+1}}\int_L  \fctVOtriall{l+1} \fctVItrial  dx  \Big)  \\	
	    & =    \left\{ \begin{array}{c} \overline{ u}          \\ \overline{\widetilde h }   \end{array} \right\}  
	            e^{ikx_{m}} \int_L  \fctVItrial       \fctVItrial  dx ,  \ \forall \indxVItrial = 1,\dots N. 
	    \end{split}
  \end{equation}
  Dividing the latter equations by $e^{ikx_{m}}$, applying $\frac{\Delta x}{2} = x_{l+1} -  x_{m} = x_{m} -  x_{l-1}$, and 
  using the integral values of Appendix~\ref{append_int_values}, there follows 
  \begin{equation}\label{equ_split_MetP0P0_discrete2}
  \begin{split}
	   \left\{ \begin{array}{c}  \overline{\widetilde u} \\ \overline{h }      \end{array} \right\} 
	   \frac{\Delta x}{2} \Big(  e^{-ik \frac{\Delta x}{ 2} }   +  e^{ik\frac{\Delta x}{ 2} } \Big) 
         & =  \left\{ \begin{array}{c} \overline{ u}          \\ \overline{\widetilde h }   \end{array} \right\} \, ,  \\
	   \left\{ \begin{array}{c}  \overline{\widetilde u} \\ \overline{h }      \end{array} \right\} 
	        \cos \left(k \frac{\Delta x}{ 2}  \right)   
	   & = \left\{ \begin{array}{c} \overline{ u}          \\ \overline{\widetilde h }   \end{array} \right\} \, .
  \end{split}
  \end{equation}
  Combining the preceding equations with \eqref{equ_disp_splitP0P1_2} gives the matrix-vector equations
  \begin{equation}\label{equ_split_disp_P0P1_Matrix_v3}
	    \left( \begin{array}{cccc} 
	           i \omega & i \frac{2g}{\Delta x}  \sin (k \frac{\Delta x}{2})  & 0 & 0 \\ 
	           0 & 0 &  i \frac{2H}{\Delta x}  \sin (k \frac{\Delta x}{2}) & i \omega  \\ 
	           1 & 0& - \cos (k \frac{\Delta x}{2})   &   0  		  \\	 
	           0 & - \cos (k \frac{\Delta x}{2})   &   0  & 1		  
		    \end{array} \right) 
	    \left( \begin{array}{c} 
		  \overline{u} \vspace{0.4ex} \\ \overline{h} \vspace{0.4ex} \\ \overline{\widetilde u} \vspace{0.4ex}\\ \overline{\widetilde h} \\ 
		\end{array} \right) = 0 \, .
  \end{equation}
  
  To find nontrivial solutions of this system, we set the determinant of this coefficient matrix to zero.
  This gives the discrete dispersion relation
  \begin{equation}
     c_d = \frac{{\omega_{00}} }{k} = \pm \sqrt{gH} \frac{\tan (k\frac{\Delta x}{2})}{k \frac{\Delta x}{2}} ,
  \end{equation}
  where $c_d \rightarrow c$ for $k \rightarrow 0$. At wave number $k= \frac{\pi}{\Delta x}$, there is no second 
  root but a singularity. That is, with increasing wave number $k\rightarrow \frac{\pi}{\Delta x}$, 
  the phase velocity tends to infinity. Therefore, we consider the low accuracy scheme $\GPOu$--$\,\GPOh$ as unstable.

  \paragraph{Medium accuracy closure.} 
  
  Similarly to Sect.~\ref{sec_disc_metr}, we discuss two medium accuracy closures.
  For the $\GPIu$--$\,\GPOh$ scheme corresponding to the discrete Hodge-star pair $(\starULOh,\starHLIh)$, 
  the discrete metric equations are given by \eqref{equ_splitFE_m6}.
  Using the expansions of the variables as above,
  the resulting closure equation for the velocities is given by \eqref{equ_split_MetP1P1_discrete},
  and for the height fields by \eqref{equ_split_MetP0P0_discrete}. 
  Combining the latter with the topological equations~\eqref{equ_disp_splitP0P1_2},
  we obtain the matrix-vector equations 
  \begin{equation}\label{equ_split_disp_P0P1_Matrix_v2}
	    \left( \begin{array}{cccc} 
	           i \omega & i \frac{2g}{\Delta x}  \sin (k \frac{\Delta x}{2})  & 0 & 0 \\ 
	           0 & 0 &  i \frac{2H}{\Delta x}  \sin (k \frac{\Delta x}{2}) & i \omega  \\ 
                  - \cos (k \frac{\Delta x}{2}) & 0 & \Big(  \frac{2}{3}  +   \frac{1}{3} \cos (k\Delta x)    \Big) & 0  \\  
                   0 & - \cos (k \frac{\Delta x}{2})  & 0 &   1  		  
		    \end{array} \right) 
	    \left( \begin{array}{c} 
		  \overline{u} \vspace{0.4ex} \\ \overline{h} \vspace{0.4ex} \\ \overline{\widetilde u} \vspace{0.4ex}\\ \overline{\widetilde h} \\ 
		\end{array} \right) = 0 \, .
  \end{equation}
  
  For the $\GPOu$--$\,\GPIh$ scheme corresponding to the discrete Hodge-star pair $(\starULIh,\starHLOh)$, 
  the discrete metric equations are given by \eqref{equ_splitFE_m7}. Expanding the variables 
  as above, the resulting closure equation for the velocities is now given by \eqref{equ_split_MetP0P0_discrete}
  and for the height fields by \eqref{equ_split_MetP1P1_discrete}. 
  The combination of the latter with the discrete topological equations yields
  \begin{equation}\label{equ_split_disp_P0P1_Matrix_v2_1}
	    \left( \begin{array}{cccc} 
	           i \omega & i \frac{2g}{\Delta x}  \sin (k \frac{\Delta x}{2})  & 0 & 0 \\ 
	           0 & 0 &  i \frac{2H}{\Delta x}  \sin (k \frac{\Delta x}{2}) & i \omega  \\  
	           1 & 0& - \cos (k \frac{\Delta x}{2})   &   0  		  \\
                   0 & \Big(  \frac{2}{3}  +   \frac{1}{3} \cos (k\Delta x)    \Big) & 0 & - \cos (k \frac{\Delta x}{2}) 
		    \end{array} \right) 
	    \left( \begin{array}{c} 
		  \overline{u} \vspace{0.4ex} \\ \overline{h} \vspace{0.4ex} \\ \overline{\widetilde u} \vspace{0.4ex}\\ \overline{\widetilde h} \\ 
		\end{array} \right) = 0 \, . 
  \end{equation}
  
  Setting the determinants of these coefficient matrices to zero, there follows for both cases the same 
  discrete dispersion relation
  \begin{equation}\label{equ_disp_med}
           c_d = \frac{{\omega_{10}} }{k} = \pm \sqrt{gH} \frac{\sin (k\frac{\Delta x}{2})}{k \frac{\Delta x}{2}} \left[ \frac{ 3  }{ (2 + \cos (k\Delta x))  } \right]^{\frac{1}{2}},
  \end{equation}
  with $c_d \rightarrow c$ for $k \rightarrow 0$. 
  For the wave number $k= \frac{\pi}{\Delta x}$, there is neither a second root nor a singularity.
  The dispersion relation~\eqref{equ_disp_med} equals exactly the one in \eqref{equ_disp_rel_P0P1}
  of the stable P1--P0 scheme. Hence, we consider both medium accuracy schemes as stable.

  \begin{rem}
   Considering the matrix-vector equations for all four cases, we realize that the first two lines remain 
   unchanged while the second two lines change with the choice of discrete metric equations. 
   Incorporating how the straight and twisted forms are connected, it is the choice of the metric equations
   that determine -- by these last two lines -- the schemes' discrete dispersion relations. 
  \end{rem}
      


\end{document}